\documentclass{amsart}

\RequirePackage{amsmath}
\RequirePackage{amsopn}
\RequirePackage{amsfonts}
\RequirePackage{amssymb}
\RequirePackage{amsthm}

\usepackage{varioref}
\usepackage[hyperindex]{hyperref}
\usepackage{amsrefs}
\usepackage[all]{xy}

\newtheorem{thm}{Theorem}[section]

\newtheorem{cor}[thm]{Corollary}
\newtheorem{lemma}[thm]{Lemma}
\newtheorem{prop}[thm]{Proposition}

\theoremstyle{definition}
\newtheorem{definition}[thm]{Definition}

\theoremstyle{remark}
\newtheorem{remark}[thm]{Remark}
\newtheorem{example}[thm]{Example}

\numberwithin{equation}{section}

\def\mathcs{C^{*}}
\newcommand{\cs}{\ensuremath{\mathcs}}

\DeclareMathSymbol{\rtimes}{\mathbin}{AMSb}{"6F}
\def\R{\mathbf{R}}
\def\C{\mathbf{C}}
\def\T{\mathbf{T}}
\def\Z{\mathbf{Z}}
\def\K{\mathcal{K}}
\let\ker\relax

\DeclareMathOperator{\Ind}{Ind}

\DeclareMathOperator{\Prim}{Prim}
\DeclareMathOperator*{\supp}{supp}
\DeclareMathOperator{\ker}{ker}

\DeclareMathOperator{\Aut}{Aut}

\DeclareMathOperator{\id}{id}
\def\set#1{\{\,#1\,\}}
\let\tensor=\otimes
\def\restr#1{|_{{#1}}}

\newbox\hidebox
\def\spechide#1{\setbox\hidebox=\hbox{$#1$}
\hbox to\wd\hidebox{$\box\hidebox^\wedge$\hss}}
\makeatletter
\def\labelenumi{\textnormal{(\@alph\c@enumi)}}
\def\theenumi{\@alph \c@enumi}
\def\labelenumii{\textnormal{(\@roman\c@enumii)}}
\def\theenumii{\@roman \c@enumii}
\newcount\charno
\def\alphapart#1{\charno=96
\advance\charno by#1\char\charno}
\def\partref#1{\textup{(}\textnormal{\alphapart{#1}}\textup{)}}
\makeatother
\def\<{\langle}
\def\>{\rangle}
\let\ipscriptstyle=\scriptscriptstyle
\def\lipsqueeze{{\mskip -3.0mu}}
\def\ripsqueeze{{\mskip -3.0mu}}
\def\ipcomma{\nobreak\mathrel{,}\nobreak}
\newbox\ipstrutbox
\setbox\ipstrutbox=\hbox{\vrule height8.5pt
width 0pt}
\def\ipstrut{\copy\ipstrutbox}
\def\lip#1<#2,#3>{\mathopen{\relax_{\ipstrut\ipscriptstyle{
#1}}\lipsqueeze
\langle} #2\ipcomma #3 \rangle}
\def\blip#1<#2,#3>{\mathopen{\relax_{\ipstrut
\ipscriptstyle{ #1}}\lipsqueeze\bigl\langle} #2\ipcomma #3 \bigr\rangle}
\def\rip#1<#2,#3>{\langle #2\ipcomma #3
\rangle_{\ripsqueeze\ipstrut\ipscriptstyle{#1}}}
\def\brip#1<#2,#3>{\bigl\langle #2\ipcomma #3
\bigr\rangle_{\ripsqueeze\ipstrut\ipscriptstyle{#1}}}
\def\angsqueeze{\mskip -6mu}
\def\smangsqueeze{\mskip -3.7mu}
\def\trip#1<#2,#3>{\langle\smangsqueeze\langle #2\ipcomma #3
\rangle\smangsqueeze\rangle_{\ripsqueeze\ipstrut\ipscriptstyle{#1}}}
\def\btrip#1<#2,#3>{\bigl\langle\angsqueeze\bigl\langle #2\ipcomma
#3
\bigr\rangle
\angsqueeze\bigr\rangle_{\ripsqueeze\ipstrut\ipscriptstyle{#1}}}
\def\tlip#1<#2,#3>{\mathopen{\relax_{\ipstrut\ipscriptstyle{
#1}}\lipsqueeze \langle\smangsqueeze\langle} #2\ipcomma #3
\rangle\smangsqueeze\rangle}
\def\btlip#1<#2,#3>{\mathopen{\relax_{\ipstrut\ipscriptstyle{
#1}}\lipsqueeze
\bigl\langle\angsqueeze\bigl\langle} #2\ipcomma #3
\bigr\rangle\angsqueeze\bigr\rangle}

\def\ip(#1|#2){(#1\mid #2)}
\def\bip(#1|#2){\bigl(#1 \mid #2\bigr)}
\def\Bip(#1|#2){\Bigl( #1 \bigm| #2 \Bigr)}
\DeclareMathOperator{\br}{Br}

\newcommand{\cox}{\ensuremath{C_{0}(X)}}

\IfFileExists{mathrsfs.sty}{\usepackage{mathrsfs}}
{\let\mathscr\mathcal} 
\newcommand{\ib}{im\-prim\-i\-tiv\-ity bi\-mod\-u\-le}

\newcommand{\sme}{\,\mathord{\mathop{\text{--}}\nolimits_{\relax}}\,}
\newcommand{\modulefont}[1]{\mathsf{#1}}
\newcommand{\X}{\modulefont{X}}

\newcommand{\bundlefont}[1]{\mathscr{#1}}
\newcommand{\A}{\bundlefont A}
\newcommand{\B}{\bundlefont B}
\newcommand{\E}{\bundlefont E}
\newcommand{\HH}{\bundlefont H}
\newcommand{\XX}{\bundlefont X}
\newcommand{\VV}{\bundlefont V}

\newcommand{\aga}{(\A,G,\alpha)}
\def\cp(#1,#2,#3){#1\rtimes_{#3}#2}
\newcommand{\acg}{\cp(\A,G,\alpha)}
\newcommand{\go}{G^{(0)}}
\newcommand{\ho}{H^{(0)}}
\def\sa_#1(#2;#3){\Gamma_{#1}(#2;#3)}

\def\bchip<#1,#2>{\tlip\cp(\B,H,\beta)<{#1},{#2}>}
\def\bbchip<#1,#2>{\btlip\cp(\B,H,\beta)<{#1},{#2}>}
\def\acgip<#1,#2>{\trip\cp(\A,G,\alpha)<{#1},{#2}>}
\def\bacgip<#1,#2>{\btrip\cp(\A,G,\alpha)<{#1},{#2}>}

\newcommand{\Es}{\E^{*}}
\newcommand{\ps}{p^{*}}
\newcommand{\dual}{{}^\flat}
\newcommand{\op}[1]{{#1}^{\text{\normalfont op}}}
\newcommand{\Xop}{\op{X}}
\renewcommand{\H}{\mathcal{H}}
\newcommand{\half}{\frac12}
\newcommand{\hoo}{\H_{00}}
\let\phi\varphi
\usepackage{bbold}
\def\charfcn#1{\mathbb{1}_{#1}}

\newcommand{\usc}{up\-per-semi\-con\-tin\-uous}
\newcommand{\stars}{*_{s}}
\newcommand{\starr}{*_{r}}
\newcommand{\cc}{\mathscr{C}}
\def\sacc(#1;#2){\mathscr{G}(#1;#2)}
\newcommand{\ccra}{\sacc(G;r^{*}\A)}
\newcommand{\gux}{G\backslash X}
\newcommand{\lhlc}{locally Hausdorff, locally compact}
\newcommand{\ttwo}{Hausdorff}
\newcommand{\xmg}{X/G}
\newcommand{\cchrb}{\sacc(H;r^{*}\B)}
\newcommand{\ccxe}{\sacc(X;\E)}
\newcommand{\myth}{\Upsilon{}}
\newcommand{\xssx}{X*_{s}X}

\theoremstyle{definition}
\newtheorem{claim}[thm]{Claim}

\newcommand{\bb}{\mathcal{B}^{b}}
\renewcommand{\L}{\mathcal{L}}
\newcommand\M{\mathcal{M}}

\newcommand\atensor{\odot}
\newcommand{\gcbatho}{\cc(G)\atensor\H_{0}}
\newcommand{\ccgatho}{\cc(G)\atensor\H_{0}}

\newcommand\N{\mathcal{N}}
\newcommand\bboc{\mathcal{B}^{1}_{c}}
\newcommand\bbocg{\mathscr{B}^{1}(G)}
\newcommand\bbocgatho{\bbocg\atensor\H_{0}}
\newcommand\D{\mathcal{D}}
\newcommand{\ilt}{inductive limit topology}
\renewcommand\AA{\mathscr{A}_{0}}
\newcommand\spxi{\zeta}

\newcommand\hoop{\hoo'}

\newcommand\strutornot{\relax}
\def\ipu<#1,#2>{\langle #1 \ipcomma #2 \rangle_{\strutornot u}}
\def\ipv<#1,#2>{\langle #1 \ipcomma #2 \rangle_{\strutornot v}}
\def\bipu<#1,#2>{\bigl\langle #1 \ipcomma #2 \bigr\rangle_{\strutornot u}}
\def\bipv<#1,#2>{\bigl\langle #1 \ipcomma #2 \bigr\rangle_{\strutornot
    v}}
\def\bipb #1<#2,#3>{\bigl\langle #2 \ipcomma #3 \bigr\rangle_{\strutornot #1}}

\newcommand\ccgathoop{\cc(G)\atensor\hoop}
\newcommand\hatpi{\check u}

\newcommand\trivial[1]{\mathscr{T}_{#1}}
\newcommand\lt{\operatorname{lt}}
\newcommand\group[1]{\mathfrak{#1}}
\newcommand\gG{\group{G}}
\newcommand\gK{\group{K}}
\newcommand\gH{\group{H}}
\newcommand\groupaction[1]{\check{#1}}
\newcommand\galpha{\groupaction{\alpha}}
\newcommand\gsigma{\groupaction{\sigma}}
\newcommand\gtau{\groupaction{\tau}}
\newcommand\gbeta{\groupaction{\beta}}

\SelectTips{cm}{}

\expandafter\ifx\csname BibSpec\endcsname\relax\else

\BibSpec{collection.article}{%
    +{}  {\PrintAuthors}                {author}
    +{,} { \textit}                     {title}
    +{.} { }                            {part}
    +{:} { \textit}                     {subtitle}
    +{,} { \PrintContributions}         {contribution}
    +{,} { \PrintConference}            {conference}
    +{}  {\PrintBook}                   {book}
    +{,} { }                            {booktitle}
    +{,} { }                            {series}
    +{,} { \voltext}                    {volume}
    +{,} { }                            {publisher}
    +{,} { }                            {organization}
    +{,} { }                            {address}
    +{,} { \PrintDateB}                 {date}
    +{,} { pp.~}                        {pages}
    +{,} { }                            {status}
    +{,} { \PrintDOI}                   {doi}
    +{,} { available at \eprint}        {eprint}
    +{}  { \parenthesize}               {language}
    +{}  { \PrintTranslation}           {translation}
    +{;} { \PrintReprint}               {reprint}
    +{.} { }                            {note}
    +{.} {}                             {transition}
    +{}  {\SentenceSpace \PrintReviews} {review}
}
\BibSpec{article}{%
    +{}  {\PrintAuthors}                {author}
    +{,} { \textit}                     {title}
    +{.} { }                            {part}
    +{:} { \textit}                     {subtitle}
    +{,} { \PrintContributions}         {contribution}
    +{.} { \PrintPartials}              {partial}
    +{,} { }                            {journal}
    +{}  { \textbf}                     {volume}
    +{}  { \PrintDatePV}                {date}
    +{,} { \eprintpages}                {pages}
    +{,} { }                            {status}
    +{,} { \PrintDOI}                   {doi}
    +{,} { available at \eprint}        {eprint}
    +{}  { \parenthesize}               {language}
    +{}  { \PrintTranslation}           {translation}
    +{;} { \PrintReprint}               {reprint}
    +{.} { }                            {note}
    +{.} {}                             {transition}
    +{}  {\SentenceSpace \PrintReviews} {review}
}
\fi

\renewcommand{\MR}[1]{\relax}
\allowdisplaybreaks[2]

\begin{document}
\date{20 June 2007}
\title[Equivalence Theorem]{Renault's Equivalence Theorem for Groupoid
  Crossed Products}

\author{Paul S. Muhly}
\address{Department of Mathematics\\
University of Iowa\\
Iowa City, IA 52242}
\email{pmuhly@math.uiowa.edu}

\author{Dana P. Williams}
\address{Department of Mathematics\\
Dartmouth College\\
Hanover, NH 03755-3551}
\email{dana.williams@dartmouth.edu}

\begin{abstract}
  We provide an exposition and proof of Renault's equivalence theorem
  for crossed products by  \lhlc{}
  groupoids.  Our approach stresses the bundle approach, concrete
  imprimitivity bimodules and is a preamble to a detailed treatment of
  the Brauer semigroup for a \lhlc{} groupoid.
\end{abstract}

\maketitle

\tableofcontents

\section{Introduction}
\label{introduction}

Our objective in this paper is to present an exposition of the theory
of groupoid actions on so-called upper-semicontinuous-$\cs$-bundles
and to present the rudiments of their associated crossed product
$\cs$-algebras. In particular, we shall extend the equivalence theorem
from \cite{mrw:jot87} and \cite{ren:jot87}*{Corollaire~5.4} to cover
locally compact, but not necessarily Hausdorff, groupoids acting on
such bundles. Our inspiration for this project derives from
investigations we are pursuing into the structure of the \emph{Brauer
  semigroup}, $S(G)$, of a locally compact groupoid $G$, which is
defined to be a collection of Morita equivalence classes of actions of
the groupoid on upper-semicontinuous-$\cs$-bundles.
The semigroup $S(G)$ arises in numerous guises in the literature and
one of our goals is to systematize their theory.  For this purpose, we
find it useful to work in the context of groupoids that are not
necessarily Hausdorff.  It is well known that complications arise when
one passes from Hausdorff groupoids to non-Hausdorff groupoids and
\emph{some} of them are dealt with in the literature.  Likewise,
conventional wisdom holds that there is no significant difference
between upper-semicontinuous-$\cs$-bundles and ordinary $\cs$-bundles;
one needs only to be careful.  However, there are subtle points in
both areas and it is fair to say that they have not been addressed or
collated in a fashion that is suitable for our purposes or for other
purposes where such structures arise.  Consequently, we believe that
it is useful and timely to write down complete details in one place
that will serve the needs of both theory and applications.

The non-Hausdorff locally compact spaces that enter the theory are not
arbitrary.  They are what is known as locally Hausdorff.  This means
that each point has a Hausdorff neighborhood. Nevertheless, such a
space need not have any non-trivial continuous functions.  As Connes
observed in \cites{con:lnm79, con:pspm80}, one has to replace
continuous functions by linear combinations of functions that are
continuous with compact support when restricted to certain locally
Hausdorff and locally compact sets, but are not continuous globally.
While at first glance, this looks like the right replacement of
continuous compactly supported functions in the Hausdorff setting, it
turns out that these functions are a bit touchy to work with, and
there are some surprises with which one must deal.  We begin our
discussion, therefore, in Section \ref{sec:locally-hausd-spac} by
reviewing the theory.  In addition to recapping some of the work in
the literature, we want to add a few comments of our own that will be
helpful in the sequel.  There are a number of ``standard'' results in
the Hausdorff case which are considerably more subtle in the \lhlc{}
case.  In Section~\ref{sec:cox-algebras} we turn to \cox-algebras.
The key observation here is that every \cox-algebra is actually the
section algebra of an \usc-\cs-bundle.  Since our eventual goal is the
equivalence theorem (Theorem~\ref{thm-renault}), we have to push the
envelope slightly and look at \usc-Banach bundles over \lhlc{} spaces.

In Section~\ref{sec:group-cross-prod}, we give the definition of, and
examine the basic properties of, groupoid crossed products.  Here we
are allowing (second countable) \lhlc{} groupoids acting  on
$C_{0}(\go)$-algebras.    In Section~\ref{sec:equiv-group-dynam} we
state the main object of this effort: Renault's equivalence theorem.  

Our version of the proof of the equivalence theorem requires some
subtle machinations with approximate identities and
Section~\ref{sec:appr-ident} is devoted to the details.  The other
essential ingredients of the proof require that we talk about
covariant representations of groupoid dynamical systems and prove a
disintegration theorem analogous to that for ordinary groupoid
representations.  This we do in Section~\ref{sec:covar-repr}.  With
all this machinery in hand, the proof of the equivalence theorem is
relatively straightforward and the remaining details are given in
Section~\ref{sec:proof}.  

In Section~\ref{sec:another-example} and
Section~\ref{sec:main-example} we look at two very important
applications of the equivalence theorem inspired by the constructions
and results in \cite{kmrw:ajm98}.

Since the really deep part of the proof of the equivalence theorem is
Renault's disintegration theorem (Theorem~\ref{thm-ren-4.2}), and
since that result --- particularly the details for \lhlc{} groupoids
--- is hard to sort out of the literature, we have included a complete
proof in Appendix~\ref{sec:muhlys-proof-theorem}.  Since that proof
requires some gymnastics with the analogues of Radon measures on
\lhlc{} spaces, we have also included a brief treatment of the results
we need in Appendix~\ref{sec:radon-measures}.

\subsection*{Assumptions}

Because Renault's disintegration
result is mired in
direct integral theory, it is necessary to restrict to second
countable groupoids and separable \cs-algebras for our main results.
We have opted to make those assumptions throughout --- at least
wherever possible.  In addition, we have adopted the common
conventions that all homomorphisms between \cs-algebras are
presumed to be $*$-preserving, and that representations of \cs-algebras are
assumed to be nondegenerate.  

\section{Locally Hausdorff Spaces, Groupoids 
and Principal $G$-spaces}
\label{sec:locally-hausd-spac}

In applications to noncommutative geometry --- in particular, to the
study of foliations --- in applications to group representation theory,
and in applications to the study of various dynamical systems, the
groupoids that arise often fail to be Hausdorff.  They are, however,
\emph{locally Hausdorff}, which means that each point has a
neighborhood that is Hausdorff in the relative topology.  Most of the
non-Hausdorff, but locally Hausdorff spaces $X$ we shall meet will,
however, also be locally compact. That is, each point in $X$ will have
a Hausdorff, compact neighborhood.\footnote{We do not follow Bourbaki
  \cite{bou:generalc}, where a space is compact if and only if it
  satisfies the every-open-cover-admits-a-finite-subcover-condition
  \emph{and} is Hausdorff} In such a space compact sets need not be
closed, but, at least, points are closed.

Non-Hausdorff, but locally Hausdorff spaces often admit a paucity of
continuous compactly supported functions. Indeed, as shown in the
discussion following \cite{khoska:jram02}*{Example~1.2}, there may be
no non-zero functions in $C_{c}(X)$. Instead, the accepted practice is
to use the following replacement for $C_{c}(X)$ introduced by Connes
in \citelist{\cite{con:pspm80}\cite{con:lnm79}}.  If $U$ is a
Hausdorff open subset of $X$, then we can view functions in $C_{c}(U)$
as functions on $X$ by defining them to be zero off $U$.  Unlike the
Hausdorff case, however, these extended functions may no longer be
continuous, or compactly supported on $X$.\footnote{Recall that the
  support of a function is the \emph{closure} of the set on which the
  function is nonzero.  Even though functions in $C_{c}(U)$ vanish off
  a compact set, the closure in $X$ of the set where they don't vanish
  may not be compact.}  Connes's replacement for $C_{c}(X)$ is the
subspace, $\cc(X)$, of the complex vector space of functions on $X$
spanned by the elements of $C_{c}(U)$ for all open Hausdorff subsets
$U$ of $X$.  Of course, if $X$ is Hausdorff, then $\cc(X)=C_{c}(X)$.
The notation $C_{c}(X)$ is often used in place of $\cc(X)$.  However,
since elements of $\cc(X)$ need be neither continuous nor compactly
supported, the $C_{c}$ notation seems ill-fitting.  Nevertheless, if
$f\in\cc(X)$, then there is a compact set $K_{f}$ such that $f(x)=0$
if $x\notin K_{f}$. As is standard, we will say that a net
$\set{f_{i}}\subset \cc(X)$ converges to $ f\in\cc(X)$ in the
\emph{inductive limit topology} on $\cc(X)$ if there is a compact set
$K$, independent of $i$, such that $f_{i}\to f$ uniformly and each
$f_{i}(x)=0$ if $x\notin K$.

While it is useful
for many purposes, the introduction of $\cc(X)$ is no panacea:
$\cc(X)$ is not closed under pointwise products, in general, and
neither is it closed under the process of ``taking the modulus'' of a
function.  That is, if $f\in\cc(X)$ it need not be the case that
$|f|\in\cc(X)$ \cite{pat:groupoids99}*{p.~32}.  A straightforward
example illustrating the problems with functions in $\cc(X)$ is the
following.

\begin{example}
  \label{ex-not-closed-abs}
  As in \cite{khoska:jram02}*{Example~1.2}, we form a groupoid $G$ as
  the topological quotient of $\Z\times[0,1]$ where for all $t\not=0$
  we identify $(n,t)\sim (m,t)$ for all $n,m\in\Z$.  (Thus as a set,
  $G$ is the disjoint union of $\Z$ and $(0,1]$).  If $f\in C[0,1]$,
  then we let $f^{n}$ be the function in $C\bigl(\set
  n\times[0,1]\bigr)\subset \cc(G)$ given by
  \begin{equation*}
    f^{n}(m,t):=
    \begin{cases}
      f(t)&\text{if $t\not=0$,}\\
      f(0)&\text{if $m=n$ and $t=0$ and} \\
      0&\text{otherwise.}
    \end{cases}
  \end{equation*}
  Then in view of \cite{khoska:jram02}*{Lemma~1.3}, every $F\in\cc(G)$
  is of the form
  \begin{equation*}
    F=\sum_{i=1}^{k} f_{i}^{n_{i}}
  \end{equation*}
  for functions $f_{1},\dots,f_{k}\in C[0,1]$ and integers $n_{i}$.
  In particular, if $F\in\cc(G)$ then we must have
  \begin{equation}
    \label{eq:61}
    \sum_{n} F(n,0)=\lim_{t\to0^{+}}F(0,t).
  \end{equation}
  Let $g(t)=1$ for all $t\in [0,1]$, and let $F\in\cc(G)$ be defined
  by $F=g^{1}-g^{2}$.  Then
  \begin{equation*}
    F(n,t)=
    \begin{cases}
      \phantom{-} 1&\text{if $t=0$ and $n=1$,}\\
      -1&\text{if $t=0$ and $n=2$ and}\\
      \phantom{-}0&\text{otherwise.}
    \end{cases}
  \end{equation*}
  Not only is $F$ an example of a function in $\cc(G)$ which is not
  continuous on $G$, but $|F|=\max(F,-F)=F^{2}$ fails to satisfy
  \eqref{eq:61}.  Therefore $|F|\notin \cc(G)$ even though $F$ is.
  This also shows that $\cc(G)$ is not closed under pointwise products
  nor is it a lattice: if $F,F'\in\cc(G)$, it does
  not follow that either $\max(F,F')\in\cc(G)$ or $\min(F,F')\in
  \cc(G)$.
\end{example}

We shall always assume that the \lhlc{} spaces $X$ with which we deal
are second countable, i.e., we shall assume there is a countable basis
of open sets. Since points are closed, the Borel structure on $X$
generated by the open sets is countably separated.  Indeed, it is
standard.  The reason is that every second countable, compact
Hausdorff space is Polish \cite{wil:crossed}*{Lemma~6.5}.  Thus $X$
admits a countable cover by standard Borel spaces.  It follows that
$X$ can be expressed as a disjoint union of a sequence of standard
Borel spaces, and so is standard.

The functions in $\cc(X)$ are all Borel.  By a \emph{measure} on $X$
we mean an ordinary, positive measure $\mu$ defined on the Borel
subsets of $X$ such that the restriction of $\mu$ to each Hausdorff
open subset $U$ of $X$ is a Radon measure on $U$.  That is, the
measures we consider restrict to regular Borel measures on each
Hausdorff open set and, in particular, they assign finite measure to
each compact subset of a Hausdorff open set.  (Recall that for second
countable locally compact Hausdorff spaces, Radon measures are simply
regular Borel measures.)  If $\mu$ is such a measure, then every
function in $\cc(X)$ is integrable.  (For more on Radon measures on
\lhlc{} spaces, see Appendix~\ref{sec:radon-meas-locally}.)

Throughout, $G$ will denote a locally Hausdorff, locally compact
groupoid.  Specifically we assume that $G$ is a groupoid endowed with
a topology such that
\begin{itemize}
\item [G1:] the groupoid operations are continuous,
\item [G2:] the unit space $\go$ is Hausdorff,
\item [G3:] each point in $G$ has a compact Hausdorff neighborhood,
  and
\item [G4:] the range (and hence the source) map is open.
\end{itemize}

A number of the facts about non-Hausdorff groupoids that we shall use
may be found in \cite{khoska:jram02}.  Another helpful source is the
paper by Tu \cite{tu:doc04}. Note that as remarked in
\cite{khoska:jram02}*{\S1B}, for each $u\in\go$,
$G^{u}:=\set{\gamma\in G:r(\gamma)=u}$ must be Hausdorff.  To see
this, recall that $\set u$ is closed in $G$, and observe that
\begin{equation*}
  G\stars G =\set{(\gamma,\eta)\in G\times G:s(\gamma)=s(\eta)}
\end{equation*}
is closed in $G\times G$.  Since $(\gamma,\eta)\mapsto
\gamma\eta^{-1}$ is continuous from $G\stars G$ to $G$, the diagonal
\begin{align*}
  \Delta(G^{u})&:= \set{(\gamma,\gamma)\in G^{u}\times G^{u}} \\
  &= \set{(\gamma,\eta)\in G\stars G:\gamma\eta^{-1}=u} \cap
  G^{u}\times G^{u}
\end{align*}
is closed in $G^{u}\times G^{u}$.  Hence $G^{u}$ is Hausdorff, as
claimed.  Of course, if $G$ is Hausdorff, then $\go$ is closed since
$\go = \set{\gamma \in G: {\gamma}^2 = \gamma}$ and convergent nets
have unique limits.  Conversely, if $G$ is not Hausdorff, then to see
that $\go$ fails to be closed, let $\gamma_{i}$ be a net in $G$
converging to both $\gamma$ and $\eta$ (with $\eta\not=\gamma$).
Since $\go$ is Hausdorff by [G2], we must have $s(\gamma)=s(\eta)$.
Then $\gamma_{i}^{-1}\gamma_{i}\to \gamma^{-1}\eta$ (as well as to
$\gamma^{-1}\gamma$).  Therefore $s(\gamma_{i})$ must converge to
$\gamma^{-1}\eta\notin\go$. Therefore $G$ is Hausdorff \emph{if and
  only if} $\go$ is closed in $G$.

\begin{remark}
  \label{rem-non-t2-isotropy}
  Suppose that $G$ is a non-Hausdorff, \lhlc{} groupoid.  Then there
  are distinct elements $\gamma$ and $\eta$ in $G$ and a net
  $\set{\gamma_{i}}$ converging to both $\gamma$ and $\eta$.  Since
  $\go$ is \ttwo, $s(\gamma_{i})\to u=s(\gamma)=s(\eta)$, and
  $r(\gamma_{i})\to v=r(\gamma)=r(\eta)$.  In particular,
  $\gamma^{-1}\eta$ is a non-trivial element of the isotropy group
  $G_{u}^{u}$.  In particular, a principal \lhlc{} groupoid must be
  Hausdorff.
\end{remark}

Since each $G^{u}$ is a locally compact Hausdorff space, $G^{u}$ has
lots of nice Radon measures.  Just as for Hausdorff
locally compact groupoids, a \emph{Haar system} on $G$ is a family of
measures on $G$, $\set{\lambda^{u}}_{u\in\go}$ on $G$, such that:
\begin{enumerate}
\item For each $u\in\go$, $\lambda^{u}$ is supported on $G^{u}$ and
  the restriction of $\lambda^{u}$ to $G^{u}$ is a regular Borel
  measure.
\item For all $\eta\in G$ and $f\in \cc(G)$,
  \begin{equation*}
    \int_{G} f(\eta\gamma)\,d\lambda^{s(\eta)}(\gamma) = \int_{G}
    f(\gamma) \,d\lambda^{r(\eta)}(\gamma)
  \end{equation*}
\item For each $f\in\cc(G)$,
  \begin{equation*}
    u\mapsto \int_{G}f(\gamma)\,d\lambda^{u}(\gamma)
  \end{equation*}
  is continuous and compactly supported on $\go$.
\end{enumerate}
We note in passing that Renault \cites{ren:jot87, ren:groupoid} and
Paterson \cite{pat:groupoids99}*{Definition~2.2.2} assume that the
measures in a Haar system $\set{\lambda^{u}}_{u\in\go}$ have full
support; i.e., they assume that $\supp(\lambda^{u}) = G^{u}$, whereas
Khoshkam and Skandalis don't (see \cite{khoska:jram02} and
\cite{khoska:jot04}.) It is easy to see that the union of the supports
of the $\lambda^{u}$ is an invariant set for the left action of $G$ on
$G$ (in a sense to be discussed in a moment).  If this set is all of
$G$, then we say that the Haar system is \emph{full}.  All of our
groupoids will be assumed to have full Haar systems and we shall not
add the adjective ``full'' to any Haar system we discuss.  Note that
if a groupoid satisfies G1, G2 and G3 and has a Haar system, then it
must also satisfy G4 \cite{pat:groupoids99}*{Proposition~2.2.1}.

If $X$ is a $G$-space,\footnote{Actions of groupoids on topological
  spaces are discussed in several places in the literature.  For
  example, see \cite{kmrw:ajm98}*{p.~912}.} then let
$G*X=\set{(\gamma,x):s(\gamma)=r(x)}$ and define $\Theta:G*X\to
X\times X$ by $\Theta(\gamma,x):=(\gamma\cdot x, x)$.  We say that $X$
is a \emph{proper $G$-space} if $\Theta$ is a proper
map.\footnote{Recall that a map $f:A\to B$ is proper if
  $f\times\id_{C}:A\times C\to B\times C$ is a a closed map for every
  topological space $C$ \cite{bou:generalc}*{I.10.1, Definition 1}.
  For the case of group actions, see \cite{bou:generalc}*{III.4}}
\begin{lemma}
  \label{lem-proper}
  Suppose a locally Hausdorff, locally compact groupoid $G$ acts on a
  locally Hausdorff, locally compact space $X$.  Then $X$ is a proper
  $G$-space if and only if $\Theta^{-1}(W)$ is compact in $G*X$ for
  all compact sets $W$ in $X\times X$.
\end{lemma}
\begin{proof}
  If $\Theta$ is a proper map, then $\Theta^{-1}(W)$ is compact
  whenever $W$ is by \cite{bou:generalc}*{I.10.2, Proposition~6}.

  Conversely, assume that $\Theta^{-1}(W)$ is compact whenever $W$ is.
  In view of \cite{bou:generalc}*{I.10.2, Theorem~1(b)}, it will
  suffice to see that $\Theta$ is a closed map.  Let $F\subset X*G$ be
  a closed subset, and let $E:=\Theta(F)$.  Suppose that
  $\set{(\gamma_{i},x_{i})} \subset F$ and that
  $\Theta(\gamma_{i},x_{i})=(\gamma_{i}\cdot x_{i},x_{i})\to (y,x)$.
  Let $W$ be a compact \emph{Hausdorff} neighborhood of $(y,x)$.
  Since $F$ is closed, $\Theta^{-1}(W)\cap F$ is compact and
  eventually contains $(\gamma_{i},x_{i})$.  Hence we can pass to a
  subnet, relabel, and assume that $(\gamma_{i},x_{i})\to (\gamma,z)$
  in $F\cap \Theta^{-1}(W)$.  Then $(\gamma_{i}\cdot x_{i},x_{i})\to
  (\gamma\cdot z, z)$ in $W$.  Since $W$ is Hausdorff, $z=x$ and
  $\gamma\cdot x=y$.  Therefore $(y,x)=(\gamma\cdot x,x)$ is in $E$.
  Hence $E$ is closed.  This completes the proof.
\end{proof}

\begin{remark}
  If $X$ is Hausdorff, the proof is considerably easier.  In fact, it
  suffices to assume only that $\Theta^{-1}(W)$
  pre-compact.\footnote{In the Hausdorff case, ``pre-compact'' and
    ``relatively compact'' refer to set whose closure is compact.  In
    potentially non-Hausdorff situations, such as here, we use
    ``pre-compact'' for a set which is \emph{contained} in a compact
    set.  In particular, a pre-compact set need not have compact
    closure.  (For an example, consider
    \cite{khoska:jram02}*{Example~1.2}.)}
\end{remark}

\begin{definition}
  \label{def-principal}
  A $G$-space $X$ is called \emph{free} if the equation $\gamma\cdot
  x=x$ implies that $\gamma=r(x)$.  A free and proper $G$-space is
  called a \emph{principal} $G$-space.
\end{definition}

If $X$ is a $G$ space, then we denote the orbit space by $\gux$.  The
orbit map $q:X\to\gux$ is continuous and open
\cite{muhwil:plms395}*{Lemma~2.1}.  Our next observation is that, just
as in the Hausdorff case, the orbit space for a proper $G$-space has
regularity properties comparable to those of the total space.

\begin{lemma}
  \label{lem-proper-orbit-space}
  Suppose that $X$ is a locally Hausdorff, locally compact
  \emph{proper} $G$-space.  Then $\gux$ is a locally Hausdorff,
  locally compact space.  In particular, if $C$ is a compact subset of
  $X$ with a compact Hausdorff neighborhood $K$, then $q(C)$ is
  Hausdorff in $\gux$.
\end{lemma}
\begin{proof}
  It suffices to prove the last assertion.  Suppose that $\set{x_{i}}$
  is a net in $C$ such that $G\cdot x_{i}$ converges to $G\cdot y$ and
  $G\cdot z$ for $y$ and $z$ in $C$.  It will suffice to see that
  $G\cdot y = G\cdot z$.  After passing to a subnet, and relabeling,
  we can assume that $x_{i}\to x$ in $C$ and that there are
  $\gamma_{i}\in G$ such that $\gamma_{i}\cdot x_{i}\to y$.  We may
  assume that $x_{i},\gamma_{i}\cdot x\in K$.  Since
  $\Theta^{-1}(K\times K)$ is compact and since
  $\set{(\gamma_{i},x_{i})}\subset \Theta^{-1}(K\times K)$, we can
  pass to a subnet, relabel, and assume that $(\gamma_{i},x_{i})\to
  (\gamma,w)$ in $\Theta^{-1}(K\times K)$.  Since $K$ is Hausdorff, we
  must have $w=x$.  Thus $\gamma_{i}\cdot x_{i}\to \gamma\cdot x$.
  Since $y\in C\subset K$, we must have $\gamma\cdot x=y$.  But then
  $G\cdot x=G\cdot y$.  Similarly, $G\cdot x=G\cdot z$.  Thus $G\cdot
  y=G\cdot z$, and we're done.
\end{proof}

\begin{example}
  \label{ex-g-on-g}
  If $G$ is a locally Hausdorff, locally compact groupoid, then the
  left action of $G$ on itself is free and proper.  In fact, in this
  case, $G*G=G^{(2)}$ and $\Theta$ is homeomorphism of $G^{(2)}$ onto
  $G\stars G=\set{(\gamma,\eta):s(\gamma)=s(\eta)}$ with inverse
  $\Phi(\beta,\alpha)=(\beta\alpha^{-1},\alpha)$.  Since $\Phi$ is
  continuous, $\Phi(W)=\Theta^{-1}(W)$ is compact whenever $W$ is.
\end{example}

\begin{remark}
  \label{rem-not2onnont2}
  If $G$ is a non-Hausdorff, \lhlc{} groupoid, then as the above
  example shows, $G$ acts (freely and) properly on itself.  Since this
  is a fundamental example --- perhaps even \emph{the} fundamental
  example --- we will have to tolerate actions on non-Hausdorff
  spaces.  It should be observed, however, that a \emph{Hausdorff}
  groupoid $G$ can't act properly on a non-Hausdorff space $X$.  If
  $G$ is Hausdorff, then $\go$ is closed and $\go*X$ is closed in
  $G*X$.  However $\Theta(\go*X)$ is the diagonal in $X\times X$,
  which if closed if and only if $X$ is \ttwo.
\end{remark}

\begin{remark}\label{rem-schultz}
  If $X$ is a proper $G$-space, and if $K$ and $L$ are compact subsets
  of $X$, then
  \begin{equation*}
    P(K,L):=\set{\gamma\in G:K\cap \gamma\cdot L\not=\emptyset}
  \end{equation*}
  is compact --- consider the projection onto the first factor of the
  compact set $\Theta^{-1}(K\times L)$.  If $X$ is Hausdorff, the
  converse is true; see, for example,
  \cite{anaren:amenable00}*{Proposition~2.1.9}.  However, the converse
  fails in general.  In fact, if $X$ is any non-Hausdorff, locally
  Hausdorff, locally compact space, then $X$ is, of course, a
  $G$-space for the trivial group(oid) $G=\set e$.  But in this case
  $\Theta(G*X)=\Delta(X):= \set{(x,x)\in X\times X:x\in X}$.  But
  $\Delta(X)$ is closed if and only if $X$ is Hausdorff.  Therefore,
  if $X$ is not Hausdorff, $\Theta$ is not a closed map, and therefore
  is not a proper map.\footnote{Notice that $\Theta^{-1}(K\times
    L)=\set e\times K\cap L$, and $K\cap L$ need not be compact even
    if both $K$ and $L$ are.}  Of course, in this example, $P(K,L)$ is
  trivially compact for any $K$ and $L$.  In \cite{ren:jot87}, it is
  stated that $X$ is a proper $G$-space whenever $P(K,L)$ is
  relatively compact for all $K$ and $L$ compact in $X$.  As this
  discussion shows, this is not true in the non-Hausdorff case.  If
  ``relatively compact'' in interpreted to mean contained in a compact
  set (as it always is here), then it can be shown that $P(K,L)$ is
  relatively compact for all $K$ and $L$ compact in $X$ if and only if
  $\Theta^{-1}(W)$ is relatively compact for all compact $W$
  \cite{shu:private88}.
\end{remark}

As Remark~\vref{rem-schultz} illustrates, there can be subtleties
involved when working with \lhlc{} $G$-spaces.  We record here some
technical results, most of which are routine in the Hausdorff case,
which will be of use later.

Recall that a subset $U\subset G$ is called \emph{conditionally
  compact} if $VU$ and $UV$ are pre-compact whenever $V$ is
pre-compact in $G$.  We say that $U$ is \emph{diagonally compact} if
$UV$ and $VU$ are compact whenever $V$ is compact.  If $U$ is a
diagonally compact neighborhood of $\go$, then its interior is a
conditionally compact neighborhood.  We will need to see that $G$ has
a fundamental system of diagonally compact neighborhoods of $\go$.
The result is based on a minor variation, of
\cite{ren:groupoid}*{Proof of Proposition~2.1.9} and
\cite{muhwil:ms90}*{Lemma~2.7} that takes into account the possibility
that $G$ is not Hausdorff.

\begin{lemma}
  \label{lem-diagonally-cpt}
Suppose that $G$ is a \lhlc{} groupoid.
  If $\go$ is paracompact, then $G$ has a fundamental system of
  diagonally compact neighborhoods of $\go$.
\end{lemma}
\begin{remark}
  If $G$ is second countable, then so is $\go$.  Hence $\go$ is always
  paracompact under our standing assumptions.
\end{remark}

\begin{proof}
  Let $V$ be any neighborhood of $\go$ in $G$.  Since $\go$ is
  paracompact, the shrinking lemma (cf.,
  \cite{rw:morita}*{Lemma~4.32}) implies that there is a locally
  finite cover $\set{K_{i}}$ of $\go$ such that each $K_{i}$ is a
  compact subset of $\go$ and such that the interiors of the $K_{i}$
  cover $\go$.  In view of the local finiteness, any compact subset of
  $\go$ meets only finitely many $K_{i}$.

  Let $U_{i}'$ be a compact neighborhood of $K_{i}$ in $G$ with
  $U_{i}'\subset V$.  Let $U_{i}:= U_{i}'\cap s^{-1}(K_{i})\cap
  r^{-1}(K_{i})$.  Since $s^{-1}(K_{i})$ and $r^{-1}(K_{i})$ are
  closed, $U_{i}$ is a compact set whose interior contains the
  interior of $K_{i}$, and
  \begin{equation*}
    K_{i}\subset U_{i}\subset V\cap s^{-1}(K_{i})\cap r^{-1}(K_{i}).
  \end{equation*}
  Therefore
  \begin{equation*}
    U:=\bigcup U_{i}
  \end{equation*}
  is a neighborhood of $\go$.  If $K$ is any compact subset of $\go$,
  then
  \begin{equation*}
    U\cap s^{-1}(K)=\bigcup_{K\cap K_{i}\not=\emptyset}U_{i}\cap
    s^{-1}(K). 
  \end{equation*}
  Since $s^{-1}(K)$ is closed and the union is finite, $U\cap
  s^{-1}(K)$ is compact.  Similarly, $r^{-1}(K)\cap U$ is compact as
  well.  Since $U\cdot K= (U\cap s^{-1}(K))\cdot K$, the former is
  compact as is $K\cdot U$.  Thus, $U$ is a diagonally compact
  neighborhood of $\go$ contained in $V$.
\end{proof}

\begin{remark}
  \label{rem-nbhd-base}
  We have already observed that if $G$ is not Hausdorff, then $\go$ is
  not closed in $G$.  Since points in $G$ are closed, it nevertheless
  follows that $\go$ is the intersection of all neighborhoods $V$ of
  $\go$ in $G$.  In particular, Lemma~\vref{lem-diagonally-cpt}
  implies that $\go$ is the intersection of all conditionally compact,
  or diagonally compact, neighborhoods of $\go$, provided $\go$ is
  paracompact.
\end{remark}

\begin{lemma}
  \label{lem-t2-nbhd}
  Suppose that $G$ is a \lhlc{} groupoid and that $K\subset \go$ is
  compact.  Then there is a neighborhood $W$ of $\go$ in $G$ such that
  $WK=W\cap r^{-1}(K)$ is \ttwo.
\end{lemma}
\begin{proof}
  Let $u\in K$ and let $V_{u}$ be a \ttwo{} neighborhood of $u$ in
  $G$.  Let $C_{u}\subset\go$ be a closed neighborhood of $u$ in $\go$
  such that $C_{u}\subset V_{u}$.  Let $W_{u}:= r^{-1}(\go\setminus
  C_{u})\cup V_{u}$.  Then $W_{u}$ is a neighborhood of $\go$ and
  $W_{u}C_{u} \subset V_{u}$.  Let $u_{1},\dots,u_{n}$ be such that
  $K\subset\bigcup_{i} C_{u_{i}}$, and let $W:=\bigcap W_{u_{i}}$.
  
  Suppose that $\gamma$ and $\eta$ are elements of $W\cdot K$ which
  can't be separated.  Then there is a $u\in K$ such that
  $r(\gamma)=u=r(\eta)$ (Remark~\vref{rem-non-t2-isotropy}).  Say
  $u\in C_{u_{i}}$.  Then $\gamma,\eta\in W_{u_{i}}$, and consequently
  both are in $V_{u_{i}}$.  Since the latter is \ttwo, $\gamma=\eta$.
  Thus $WK$ is \ttwo.
\end{proof}

\begin{lemma}
  \label{lem-open-nbhd}
  Suppose that $G$ is a \lhlc{} groupoid and that $X$ is a \lhlc{}
  $G$-space.  If $V$ is open in $X$ and if $K\subset V$ is compact,
  then there is a neighborhood $W$ of $\go$ in $G$ such that $W\cdot
  K\subset V$.
\end{lemma}
\begin{proof}
  For each $x\in K$ there is a neighborhood $U_{x}$ of $r(x)$ in $G$
  such that $U_{x}\cdot K\subset V$.  Let $x_{1},\dots,x_{n}$ be such
  that $\bigcup r(U_{x_{i}})\supset r(K)$.  Let $W:= \bigcup
  U_{x_{i}}\cup r^{-1}(\go\setminus r(K))$.  Then $W$ is a
  neighborhood of $\go$ and $W\cdot K\subset \bigl(\bigcup U_{x_{i}}
  \bigr) \cdot K\subset V$.
\end{proof}

The next lemma is a good example of a result that is routine in the
Hausdorff case, but takes a bit of extra care in general.

\begin{lemma}
  \label{lem-key-proper-nbhd}
  Suppose that $G$ is a \lhlc{} groupoid and that $X$ is a \lhlc{}
  \emph{free and proper}
  (right) $G$-space.  If $W$ is a neighborhood of $\go$ in $G$, then
  each $x\in X$ has a neighborhood $V$ such that the inclusion
  $(x,x\cdot \gamma)\in V\times V$ implies that $\gamma\in W$.
\end{lemma}
\begin{proof}
Fix $x\in X$.
  Let $C$ be a compact \ttwo{} neighborhood of $x$ in $X$.  If the
  lemma were false for $x$, then for each neighborhood $V$ of $x$ such that
  $V\subset C$, there would be a $\gamma_{V}\notin W$ and a $x_{V}\in
  V$ such that $( x_{V},x_{V}\cdot \gamma_{V})\in V\times V$. This
  would yield a net $\{(x_{V},\gamma_{V})\}_{\{V\subset C\}}$. Since
  \begin{equation*}
    A=\set{(x,\gamma)\in X\times G:\text{$x\in C$ and $x\cdot \gamma\in C$}}
  \end{equation*}
  is compact, we could pass to a subnet, relabel, and assume that
  $(x_{V},\gamma_{V})\to (y,\gamma)$ in $A$.  Since $C$ is \ttwo{} and
  since $x_{V}\to x$ while $x_{V}\cdot\gamma_{V}\to x$, we would have
  $x=y$ and $ x\cdot \gamma = x$.  Therefore, we would find that
  $\gamma=s(x)\in W$.  On the other hand, since $W$ is open and since
  $\gamma_{V}\notin W$ for all $V$ we would find that $\gamma\notin
  W$.  This would be a contradiction, and completes the proof.
\end{proof}

The next proposition is the non-\ttwo{} version of Lemmas 2.9~and 2.13
from~\cite{mrw:jot87}.

\begin{prop}
  \label{prop-mrw2.13}
  Suppose that $G$ is a \lhlc{} groupoid with Haar system
  $\set{\lambda^{u}}_{u\in\go}$.  Let $X$ be a \lhlc{} free and proper
  (right) $G$-space, let $q:X\to\xmg$ be the quotient map, and let
  $V\subset X$ be a \ttwo{} open set such that $q(V)$ is \ttwo.
  \begin{enumerate}
  \item If $\psi\in C_{c}(V)$, then
    \begin{equation*}
      \lambda(\psi)\bigl(q(x)\bigr) =\int_{G} \psi(x\cdot
      \gamma)\,d\lambda^{s(x)} _{G}(\gamma)
    \end{equation*}
    defines an element $\lambda(\psi)\in C_{c}\bigl(q(V)\bigr)$.
  \item If $d\in C_{c}\bigl(q(V)\bigr)$, then there is a $\psi\in
    C_{c}(V)$ such that $\lambda(\psi)=d$.
  \end{enumerate}
\end{prop}
\begin{cor}
  \label{cor-l-surjective}
  The map $\lambda$ defined in part~\partref1 of
  Proposition~\vref{prop-mrw2.13} extends naturally to a surjective
  linear map $\lambda:\cc(X)\to\cc(\xmg)$ which is continuous in the
  inductive limit topology.
\end{cor}
\begin{proof}[Proof of Corollary~\vref{cor-l-surjective}]
  Let $V$ be a \ttwo{} open subset of $X$, and let $\psi\in C_{c}(V)$.
  We need to see that $\lambda(\psi)\in\cc(\xmg)$.  Let $W$ be a open
  neighborhood of $\supp_{V}\psi$ with a compact neighborhood
  contained in $V$.\footnote{Here we use the notation $\supp_{V}$ to
    describe the support of a function on $V$ relative to $V$ as
    opposed to all of $X$. Recall that the support of a continuous
    function is the closure of the set where the function is not zero,
    and since $X$ is not necessarily Hausdorff, the closure of a set
    relative to a subset such as $V$ need not be the same as the
    closure of the subset in $X$.}  Then
  Lemma~\vref{lem-proper-orbit-space} implies that $q(W)$ is Hausdorff,
  and Proposition~\vref{prop-mrw2.13} implies that $\lambda(\psi)\in
  C_{c}\bigl(q(W)\bigr)$.  It follows that $\lambda$ extends to a
  well-defined linear surjection.  The statement about the inductive
  limit topology is clear.
\end{proof}

\begin{remark}
  \label{rem-q-system}
  In the language of \cite{ren:jot87}, the first part of the
  proposition says that the Haar system on $G$ induces a $q$-system on
  $X$ --- see \cite{ren:jot87}*{p.~69}.
\end{remark}

\begin{proof}[Proof of Proposition~\ref{prop-mrw2.13}]
  Let $D=\supp_{V}\psi$.  Since $V$ is locally compact Hausdorff,
  there is an open set $W$ and a compact set $C$ such that
  \begin{equation*}
    D\subset W \subset C \subset V.
  \end{equation*}
  Let $\Theta:X\times G\to X\times X$ be given by
  $\Theta(x,\gamma)=(x,x\cdot \gamma)$.  Since the $G$-action is
  proper,
  \begin{equation*}
    A:=\Theta^{-1}(C\times C)=\set{(x,\gamma)\in X\times G:\text{$x
        \in C$ and $x\cdot \gamma\in C$}}
  \end{equation*}
  is compact.  Moreover, if $\set{(x_{i},\gamma_{i})}$ is a net in $A$
  converging to both $(x,\gamma)$ and $(y,\eta)$ in $A$, then since
  $C$ is \ttwo, we must have $x=y$.  Then $\set{x_{i}\cdot
    \gamma_{i}}$ converges to both $x\cdot \gamma$ and $x\cdot \eta$
  in the Hausdorff set $C$.  Thus $x\cdot \gamma=x\cdot \eta$, and
  since the action is free, we must have $\gamma=\eta$.  In sum, $A$
  is \ttwo.

  Let $F:C\times G\to \C$ be defined by $F(x,\gamma)=\psi(x\cdot
  \gamma)$.  Notice that $F$ vanishes off $A$.  Let
  $K:=\operatorname{pr}_{2}(A)$ be the projection onto the second
  factor; thus, $K$ is compact in $G$.  Unfortunately, we see no
  reason that $K$ must be Hausdorff.  Nevertheless, we can cover $K$
  by \ttwo{} open sets $V_{1},\dots,V_{n}$.  Let $A_{j}:=A\cap(C\times
  V_{j})$, let $\set{f_{j}}$ be a partition of unity in $C(A)$
  subordinate to $\set{A_{j}}$ and let
  $F_{j}(x,\gamma):=f_{j}(x,\gamma)F(x,\gamma)$.  Then $F_{j}\in
  C_{c}(A_{i})$.
  \begin{claim}
    If we extend $F_{j}$ by setting to be $0$ off $A$, we can view
    $F_{j}$ as an element of $C_{c}(C\times V_{j})$.
  \end{claim}
  \begin{proof}[Proof of Claim]
    Suppose that $\set{(x_{i},\gamma_{i})}$ is a net in $C\times
    V_{i}$ converging to $(x,\gamma)$ in $C\times V_{i}$.  Let
    \begin{equation*}
      B:= (C\times G)\cap \Theta^{-1}(X\times
      W)=\set{(x,\gamma):\text{$x\in C$ and $x\cdot \gamma\in W$}}.
    \end{equation*}
    Then $B$ is open in $C\times G$ and $B\subset A$.  If
    $(x,\gamma)\in B$, then $(x_{i},\gamma_{i})$ is eventually in $B$
    and $F_{j}(x_{i},\gamma_{i})\to F_{j}(x,\gamma)$ (since $F_{j}$ is
    continuous on $A$).

    On the other hand, if $(x,\gamma)\notin B$, then
    $F_{j}(x,\gamma)=0$.  If $\set{F(x_{i},\gamma_{i})}$ does not
    converge to $0$, then we can pass to a subnet, relabel, and assume
    that there is a $\delta>0$ such that
    \begin{equation*}
      |F_{j}(x_{i},\gamma_{i})|\ge\delta\quad\text{for all $i$.}
    \end{equation*}
    This means that $f_{j}(x_{i},\gamma_{i})\not=0$ for all $i$.
    Since $f_{j}$ has compact support in $A_{j}$, we can pass to a
    subnet, relabel, and assume that $(x_{i},\gamma_{i})\to (y,\eta)$
    in $A_{j}$.  Since $C$ is \ttwo, $y=x$.  Since $V_{j}$ is \ttwo,
    $\eta=\gamma$.  Therefore $(x_{i},\gamma_{i})\to (x,\gamma)$ in
    $A$.  Since $F_{j}$ is continuous on $A$,
    $F_{j}(x,\gamma)\ge\delta$.  Since $\delta>0$, this is a
    contradiction.  This completes the proof of the claim.
  \end{proof}

  Since $C\times V_{j}$ is \ttwo, we may approximate $F_{j}$ in
  $C_{c}(C\times V_{j})$ by sums of functions of the form
  $(x,\gamma)\mapsto g(x)h(\gamma)$, as in
  \cite{mrw:jot87}*{Lemma~2.9} for example.  Hence
  \begin{equation*}
    x\mapsto \int_{G} F_{j}(x,\gamma)\,d\lambda_{G}^{s(x)}(\gamma)
  \end{equation*}
  is continuous.

  Suppose that $\set{x_{i}}$ is a net in $V$ such that $q(x_{i})\to
  q(x)$ (with $x\in V$).  If $q(x)\notin q(D)$, then since $q(D)$ is
  compact and hence closed in the \ttwo{} set $q(V)$, we eventually
  have $q(x_{i})\notin q(D)$.  Thus we eventually have
  $\lambda(\psi)\bigl(q(x_{i})\bigr)=0$, and $\lambda(\psi)$ is
  continuous at $q(x)$.  On the other hand, if $q(x)\in q(W)$, then we
  may as well assume that $x_{i}\to x$ in $C$.  But on $C$,
  \begin{equation*}
    x\mapsto\lambda(\psi)\bigl(q(x)\bigr) =
    \int_{G}F(x,\gamma)\,d\lambda_{G}^{s(x)} (\gamma) = \sum_{j}
    \int_{G}F_{j} (x,\gamma) \,d\lambda_{G}^{s(x)} (\gamma)
  \end{equation*}
  is continuous.  This completes the proof of part~\partref1.

  For part~\partref2, assume that $d\in C_{c}\bigl(q(V)\bigr)$.  Then
  $\supp_{q(V)} d$ is of the form $q(K)$ for a compact set $K\subset
  V$.  Let $g\in C_{c}(V)$ be strictly positive on $K$.  Then
  $\lambda(g)$ is strictly positive on $\supp_{q(V)}d$, and
  $\lambda\bigl(g\lambda(g) d\bigr)=d$.

\end{proof}

\begin{lemma}
  \label{lem-fund-groupoid-lemma}
  Suppose that $H$ and $G$ are \lhlc{} groupoids and that $X$ is a
  $(H,G)$-equivalence.  Let $\xssx=\set{(x,y)\in X\times
    X:s(x)=s(y)}$.  Then $\xssx$ is a principal $G$-space for the
  diagonal $G$-action.  If $\tau(x,y)$ is the unique element in $H$
  such that $\tau(x,y)\cdot y=x$, then $\tau:\xssx\to H$ is continuous
  and factors though the orbit map. Moreover, $\tau$ induces a
  homeomorphism of $\xssx/G$ with $H$.
\end{lemma}
\begin{proof}
  Clearly, $\xssx$ is a principal $G$-space and $\tau$ is a
  well-defined map on $\xssx$ onto $H$.  Suppose that
  $\set{(x_{i},y_{i})}$ converges to $(x,y)$.  Passing to a subnet,
  and relabeling, it will suffice to show that
  $\set{\tau(x_{i},y_{i})}$ has a subnet converging to $\tau(x,y)$.
  Let $L$ and $K$ be \ttwo{} compact neighborhoods of $x$ and $y$,
  respectively.  Since we eventually have
  $\set{(\tau(x_{i},y_{i}),y_{i})}$ in $\Theta^{-1}(K\times L)$, we
  can pass to a subnet, relabel, and assume that
  $\bigl(\tau(x_{i},y_{i}),y_{i}\bigr)\to (\eta,z)$ in
  $\Theta^{-1}(K\times L)$.  In particular, since $L$ is \ttwo, we
  must have $z=y$.  Since $\eta\cdot y\in K$, $x_{i}\to \eta\cdot y$
  and since $K$ is \ttwo, we must have $x=\eta\cdot y$.  Thus
  $\eta=\tau(x,y)$.  This shows that $\tau$ is continuous.

  Clearly $\tau$ is $G$-equivariant.  If $\tau(x,y)=\tau(z,w)$, then
  $s_{X}(x) = r_{H}\bigl(\tau(x,y)\bigr) =s_{X}(z)$.  Since $X$ is an
  equivalence, $z=x\cdot \gamma $ for some $\gamma\in G$.  Similarly,
  $r_{X}(y) =s_{H}\bigl(\tau(x,y)\bigr) =r_{X}(w)$, and
  $y=w\cdot\gamma'$ for some $\gamma'\in G$.  Therefore $\tau$ induces
  a bijection of $\xssx$ onto $H$.  To see that $\tau$ is open, and
  therefore a homeomorphism as claimed, suppose that
  $\tau(x_{i},y_{i})\to \tau(x,y)$.  After passing to a subnet and
  relabeling, it will suffice to see that $\set{(x_{i},y_{i})}$ has a
  subnet converging to $(x,y)$.  Let $L$ and $K$ be \ttwo{} compact
  neighborhoods of $x$ and $y$, respectively.  Since
  $\Theta^{-1}(L\times K)$ is compact, we can pass to a subnet,
  relabel, and assume that $\bigl(\tau(x_{i},y_{i}),y_{i}\bigr)\to
  (\eta,z)$ in $\Theta^{-1}(L\times K)$.  Since $K$ is \ttwo, $z=y$.
  On the other hand, we must have $x_{i}=\tau(x_{i},y_{i})\cdot
  y_{i}\to \tau(x,y)\cdot y=x$.  This completes the proof.
\end{proof}

\section{\cox-algebras}
\label{sec:cox-algebras}

A \emph{\cox-algebra} is a \cs-algebra $A$ together with a
nondegenerate homomorphism $\iota_{A}$ of $\cox$ into the center of
the multiplier algebra $M(A)$ of $A$.  The map $\iota_{A}$ is normally
suppressed and we write $f\cdot a$ in place of $\iota_{A}(f)a$.  There
is an expanding literature on \cox-algebras which describe their basic
properties; a partial list is
\cites{nil:iumj96,kas:im88,bla:bsmf96,echwil:jot01,wil:crossed}.  An essential
feature of \cox-algebras is that they can be realized as sections of a
bundle over $X$.  Specifically, if $C_{0,x}(X)$ is the ideal of
functions vanishing at $x\in X$, then
$I_{x}:=\overline{C_{0,x}(X)\cdot A}$ is an ideal in $A$, and
$A(x):=A/I_{x}$ is called the fibre of $A$ over $x$.  The image of
$a\in A$ in $A(x)$ is denoted by $a(x)$.

We are interested in fibred \cs-algebras as a groupoid $G$ must act on
the sections of a bundle that is fibred over the unit space (or over
some $G$-space).  In \cite{ren:jot87} and in \cite{kmrw:ajm98}, it was
assumed that the algebra $A$ was the section algebra of a \cs-bundle
as defined, for example, by Fell in \cite{fd:representations1}.
However recent work has made it clear that the notion of a \cs-bundle,
or for that matter a Banach bundle, as defined in this way is
unnecessarily restrictive, and that it is sufficient to assume only
that $A$ is a $C_{0}(\go)$-algebra
\cites{gal:94,gal:kt99,khoska:jot04,khoska:jram02}.  However, our
approach here, as in \cite{kmrw:ajm98} (and in \cite{ren:jot87}),
makes substantial use of the total space of the underlying bundle.
Although it predates the term ``\cox-algebra'', the existence of a
bundle whose section algebra is a given $\cox$-algebra goes back to
\cites{hof:74,hofkei:lnm79,hof:lnim77}, and to \cite{dg:banach}.  We
give some of the basic definitions and properties here for the sake of
completeness.

This definition is a minor variation on
\cite{dg:banach}*{Definition~1.1}.
\begin{definition}
  \label{def-usc-bundle}
  An \emph{\usc-Banach bundle} over a topological space $X$ is a
  topological space $\A$ together with a continuous, open surjection
  $p=p_{\A}:\A\to X$ and complex Banach space structures on each fibre
  $\A_{x}:= p^{-1}(\set x)$ satisfying the following axioms.
  \begin{itemize}
  \item [B1:] The map $a\mapsto\|a\|$ is upper semicontinuous from
    $\A$ to $\R^{+}$.  (That is, for all $\epsilon>0$,
    $\set{a\in\A:\|a\|\ge\epsilon}$ is closed.)
  \item [B2:] If $\A*\A:=\set{(a,b)\in\A\times\A:p(a)=p(b)}$, then
    $(a,b)\mapsto a+b$ is continuous from $\A*\A$ to $\A$.
  \item [B3:] For each $\lambda\in\C$, $a\mapsto \lambda a$ is
    continuous from $\A$ to $\A$.
  \item [B4:] If $\set{a_{i}}$ is a net in $\A$ such that $p(a_{i})\to
    x$ and such that $\|a_{i}\|\to 0$, then $a_{i}\to0_{x}$ (where
    $0_{x}$ is the zero element in $\A_{x}$).
  \end{itemize}
\end{definition}

Since $\set{a\in\A:\|a\|<\epsilon}$ is open for all $\epsilon>0$, it
follows that whenever $a_{i}\to0_{x}$ in $\A$, then $\|a_{i}\|\to 0$.
Therefore the proof of
\cite{fd:representations1}*{Proposition~II.13.10} implies that
\begin{itemize}
\item [B3$'$:] The map $(\lambda,a)\to\lambda a$ is continuous from
  $\C\times\A$ to $\A$.
\end{itemize}
\begin{definition}
  \label{def-cs-bundle}
  An \emph{\usc-\cs-bundle} is an \usc-Banach bundle $p_{\A}:\A\to X$
  such that each fibre is a \cs-algebra and such that
  \begin{itemize}
  \item [B5:] The map $(a,b)\mapsto ab$ is continuous from $\A*\A$ to
    $\A$.
  \item [B6:] The map $a\mapsto a^{*}$ is continuous from $\A$ to
    $\A$.
  \end{itemize}

\end{definition}

If axiom B1 is replaced by
\begin{itemize}
\item [B1$'$:] The map $a\mapsto \|a\|$ is continuous,
\end{itemize}
then $p:\A\to X$ is called a \emph{Banach bundle} (or a
\emph{\cs-bundle}).  Banach bundles are studied in considerable detail
in \S\S13--14 of Chapter~II of \cite{fd:representations1}.  As
mentioned above, the weaker notion of an \usc-Banach bundle is
sufficient for our purposes.  In fact, in view of the connection with
\cox-algebras described below, it is our opinion that \usc-Banach
bundles, and in particular \usc-\cs-bundles, provide a more natural
context in which to work.

If $p:\A\to X$ is an \usc-Banach bundle, then a continuous function
$f:X\to \A$ such that $p\circ f=\id_{X}$ is called a \emph{section}.
The set of sections is denoted by $\sa_{}(X;\A)$.  We say that
$f\in\sa_{}(X;\A)$ \emph{vanishes at infinity} if the the closed set
$\set{x\in X:|f(x)|\ge\epsilon}$ is compact for all $\epsilon>0$.  The
set of sections which vanish at infinity is denoted by
$\sa_{0}(X;\A)$, and the latter is easily seen to be a Banach space
with respect to the supremum norm: $\|f\|=\sup_{x\in X}\|f(x)\|$ (cf.
\cite{dg:banach}*{p.~10}); in fact, $\sa_{0}(X;\A)$ is a Banach
\cox-module for the natural \cox-action on sections.\footnote{We also
  use $\sa_{c}(X;\A)$ for the vector space of sections with compact
  support (i.e., $\set{x\in X:f(x)\not=0_{x}}$ has compact closure).}
In particular, the uniform limit of sections is a section.  Moreover,
if $p:\A\to X$ is an \usc-\cs-bundle, then the set of sections is
clearly a $*$-algebra with respect to the usual pointwise operations,
and $\sa_{0}(X;\A)$ becomes a \cox-algebra with the obvious
$\cox$-action.  However, for arbitrary $X$, there is no reason to
expect that there are any non-zero sections --- let alone non-zero
sections vanishing at infinity or which are compactly supported.  An
\usc-Banach bundle is said to have \emph{enough sections} if given
$x\in X$ and $a\in\A_{x}$ there is a section $f$ such that $f(x)=a$.
If $X$ is a Hausdorff locally compact space and if $p:\A\to X$ is a
Banach bundle, then a result of Douady and Soglio-H\'erault implies
there are enough sections \cite{fd:representations1}*{Appendix~C}.
Hofmann has noted that the same is true for \usc-Banach bundles over
Hausdorff locally compact spaces \cite{hof:lnim77} (although the
details remain unpublished \cite{hof:74}). In the situation we're
interested in --- namely seeing that a \cox-algebra is indeed the
section algebra of an \usc-\cs-bundle --- it will be clear that there
are enough sections.
\begin{prop}[Hofmann, Dupr\'e-Gillete]
  \label{prop-cox-usc}
  If $p:\A\to X$ is an \usc-\cs-bundle over a locally compact
  Hausdorff space $X$ (with enough sections), then $A:=\sa_{0}(X;\A)$
  is a \cox-algebra with fibre $A(x)=\A_{x}$.  Conversely, if $A$ is a
  \cox-algebra then there is an \usc-\cs-bundle $p:\A\to X$ such that
  $A$ is (isomorphic to) $\sa_{0}(X;\A)$.
\end{prop}
\begin{proof}
This is proved in \cite{wil:crossed}*{Theorem~C.26}.
\end{proof}

The next observation is useful and has a straightforward proof which
we omit.  (A similar result is proved in
\cite{wil:crossed}*{Proposition~C.24}.) 
\begin{lemma}
  \label{lem-dense-sections}
  Suppose that $p:\A\to X$ is an \usc-Banach bundle over a locally
  compact Hausdorff space $X$, and that $B$ is a subspace of
  $A=\sa_{0}(X;\A)$ which is closed under multiplication by functions
  in $\cox$ and such that $ \set{f(x):f\in B} $ is dense in $A(x)$ for
  all $x\in X$.  Then $B$ is dense in~$A$.
\end{lemma}
As an application, suppose that $p:\A\to X$ is an \usc-\cs-bundle over
a locally compact Hausdorff space $X$.  Let $A=\sa_{0}(X;\A)$ be the
corresponding \cox-algebra.  If $\tau:Y\to X$ is continuous, then the
pull-back $\tau^{*}\A$ is an \usc-\cs-bundle over $Y$.  If $Y$ is
\ttwo, then as in \cite{raewil:tams85}, we can also form the the
balanced tensor product $\tau^{*}(A):=C_{0}(Y)\tensor_{\cox}A$ which
is the quotient of $C_{0}(Y)\tensor A$ by the balancing ideal
$I_{\tau}$ generated by
\begin{equation*}
  \set{\phi (f\circ\tau)\tensor a-\phi\tensor 
    f\cdot a:\text{$\phi\in C_{0}(Y)$,
      $f\in\cox$ and $a\in A$}
  }.
\end{equation*}
If $\phi\in C_{0}(Y)$ and $a\in A$, then $\psi(\phi \tensor a)(y):=
\phi(y)a\bigl(\tau(y)\bigr)$ defines a homomorphism of
$C_{0}(Y)\tensor A$ into $\sa_{0}(Y;\tau^{*}\A)$ which factors through
$\tau^{*}(A)$, and has dense range in view of
Lemma~\vref{lem-dense-sections}.  As in the proof of
\cite{raewil:tams85}*{Proposition~1.3}, we can also see that this map
is injective and therefore an isomorphism.  Since pull-backs of
various sorts play a significant role in the theory, we will use this
observation without comment in the sequel.
\begin{remark}
  \label{rem-sections-of-pullbacks}
  Suppose that $p:\A\to X$ is an \usc-\cs-bundle over a locally
  compact Hausdorff space $X$.  If $\tau:Y\to X$ is continuous, then
  $f\in \sa_{c}(Y;\tau^{*}\A)$ if and only if there is a continuous,
  compactly supported function $\check f:Y\to \A$ such that
  $p\bigl(\check f(y)\bigr)=\tau(y)$ and such that
  $f(y)=\bigl(y,\check f(y)\bigr)$.  As is customary, we will not
  distinguish between $f$ and $\check f$.
\end{remark}

Suppose that $p:\A\to X$ and $q:\B\to X$ are \usc-\cs-bundles.  As
usual, let $A=\sa_{0}(X;\A)$ and $B=\sa_{0}(X;\B)$.  Any continuous
bundle map
\begin{equation}\label{eq:49}
  \xymatrix{\A\ar[dr]_{p}\ar[rr]^{\Phi}&&\B\ar[dl]^{q}\\&X}
\end{equation}
is determined by a family of maps $\Phi(x):A(x)\to B(x)$.  If each
$\Phi(x)$ is a homomorphism (of \cs-algebras), then we call $\Phi$ a
\emph{\cs-bundle map}.  A \cs-bundle map $\Phi$ induces a
\cox-homomorphism $\phi:A\to B$ given by
$\phi(f)(x)=\Phi\bigl(f(x)\bigr)$.

Conversely, if $\phi:A\to B$ is a \cox-homomorphism, then we get
homomorphisms $\phi_{x}:A(x)\to B(x)$ given by
$\phi_{x}\bigl(a(x)\bigr) =\phi(a)(x)$.  Then $\Phi(x):=\phi_{x}$
determines a bundle map $\Phi:\A\to \B$ as in \eqref{eq:49}.  It is
not hard to see that $\Phi$ must be continuous: Suppose that $a_{i}\to
a$ in $\A$.  Let $f\in A$ be such that $f\bigl(p(a)\bigr)=a$.  Then
$\phi(f)\bigl(p(a)\bigr) =\Phi(a)$ and
\begin{equation*}
  \|\Phi(a_{i})-\phi(f)\bigl(p(a_{i})\bigr)\|\le
  \|a_{i}-f\bigl(p(a_{i})\bigr) \|\to 0.
\end{equation*}
Therefore $\Phi(a_{i})\to \Phi(a)$ by the next lemma (which shows that
the topology on the total space is determined by the sections).
\begin{lemma}
  \label{lem-top-sections}
  Suppose that $p:\A\to X$ is an \usc-Banach-bundle.  Suppose that
  $\set{a_{i}}$ is a net in $\A$, that $a\in \A$ and that
  $f\in\sa_{0}(X;\A)$ is such that $f\bigl(p(a)\bigr)=a$.  If
  $p(a_{i})\to p(a)$ and if $\|a_{i}-f\bigl(p(a_{i})\bigr) \|\to 0$,
  then $a_{i}\to a$ in $\A$.
\end{lemma}
\begin{proof}
  We have $a_{i}-f\bigl(p(a_{i})\bigr)\to 0_{p(a)}$ by axiom~B4.
  Hence
  \begin{equation*}
    a_{i}=(a_{i}-f\bigl(p(a_{i})\bigr) + f\bigl(p(a_{i})\bigr)\to
    0_{p(a)} + a =a.\qed
  \end{equation*}
  \renewcommand{\qed}{}
\end{proof}
\begin{remark}
  \label{rem-unique-bundle}
  If $\sa_{0}(X;\A)$ and $\sa_{0}(X;\B)$ are isomorphic \cox-algebras,
  then $\A$ and $\B$ are isomorphic as \usc-\cs-bundles.  Hence in
  view of Proposition~\ref{prop-cox-usc}, every \cox-algebra is the section
  algebra of a unique \usc-\cs-bundle (up to isomorphism).
\end{remark}

\begin{remark}
  \label{rem-iso}
  If $\A$ and $\B$ are \usc-\cs-bundles over $X$ and if $\Phi:\A\to
  \B$ is a \cs-bundle map such that each $\Phi(x)$ is an isomorphism,
  then $\Phi$ is bicontinuous and therefore a \cs-bundle isomorphism.
\end{remark}
\begin{proof}
  We only need to see that if $\Phi(a_{i})\to \Phi(a)$, then $a_{i}\to
  a$.  After passing to a subnet and relabeling, it suffices to see
  that $\set{a_{i}}$ has a subnet converging to $a$.  But
  $p(a_{i})=q\bigl(\Phi(a_{i})\bigr)$ must converge to $p(a)$.  Since
  $p$ is open, we can pass to a subnet, relabel, and assume that there
  is a net $b_{i}\to a$ with $p(b_{i})=p(a_{i})$.  But we must then
  have $\Phi(b_{i})\to \Phi(a)$.  Then
  $\|b_{i}-a_{i}\|=\|\Phi(b_{i}-a_{i})\|\to 0$.  Therefore
  $b_{i}-a_{i}\to 0_{p(a)}$, and
  \begin{equation*}
    a_{i}=(a_{i}-b_{i})+b_{i}\to 0_{p(z)}+a=a.\qed
  \end{equation*}
  \renewcommand{\qed}{}
\end{proof}

If $p:\E\to X$ is an \usc-Banach bundle over a locally Hausdorff,
locally compact space $X$, then as in the scalar case, there may not
be any non-zero sections in $\sa_{c}(X;\E)$.  Instead, we proceed as
in \cite{ren:jot87} and let $\sacc(X;\E)$ be the complex vector space
of functions from $X$ to $\E$ spanned by sections in
$\sa_{c}(U;\E\restr U)$, were $U$ is any open Hausdorff subset of $X$
and $\E\restr U:=p^{-1}(U)$ is viewed as an \usc-Banach bundle over
the locally compact Hausdorff space $U$.\footnote{In the sequel, we
  will abuse notation a bit and simply write $\sa_{c}(U;\E)$ rather
  than the more cumbersome $\sa_{c}(U;\E\restr U)$.}  We say $\E$ has
enough sections if given $e\in \E$, there is a $f\in\sacc(X;\E)$ such
that $f\bigl(p(e)\bigr) = e$.  By Hofmann's result \cite{hof:lnim77},
$\E$ always has enough sections.
\begin{remark}
  \label{rem-ilt}
  Suppose that $p:\E\to X$ is an \usc-Banach bundle over a locally
  Hausdorff, locally compact space $X$.  Then we say that a net
  $\set{z_{i}}_{i\in I}$ converges to $z$ in the inductive limit
  topology on $\sacc(X;\E)$ if $z_{i}\to z$ uniformly and if there is
  a compact set $C$ in $X$ such that all the $z_{i}$ and $z$ vanish
  off $C$.  (We are not claiming that there is a topology on
  $\sacc(X;\E)$ in which these are the \emph{only} convergent nets.)
\end{remark}

\begin{lemma}
  \label{lem-prop2.14-bundles}
  Suppose that $Z$ is a \lhlc{} principal (right) $G$-bundle, that
  $p:\B\to Y$ is an \usc-Banach bundle and that $\sigma:Z/G\to Y$ is a
  continuous, open map.  Let $q:Z\to Z/G$ be the orbit map.  If
  $f\in\sacc(X;(\sigma\circ q)^{*}\B)$, then the equation
  \begin{equation*}
    \lambda(f)\bigl(x\cdot G\bigr) = \int_{G} f(x\cdot
    \gamma)\,d\lambda^{s(x)}(\gamma) 
  \end{equation*}
  defines an element $\lambda(f)\in\sacc(Z/G;\sigma^{*}\B)$.
\end{lemma}
\begin{proof}
  We can assume that $f\in\sa_{c}(V;(\sigma\circ q)^{*}\B)$ with $V$ a
  \ttwo{} open set in $Z$.  Using an approximation argument, it
  suffices to consider $f$ of the form
  \begin{equation*}
    f(x)=g(x)a\bigl(\sigma\bigl(q(x)\bigr)\bigr)
  \end{equation*}
  where $g\in C_{c}(V)$ and $a\in \sa_{c}(\sigma(q(V));\B)$.  The
  result follows from Corollary~\vref{cor-l-surjective}.
\end{proof}
\begin{remark}
  \label{rem-hypotheses}
  The hypotheses in Lemma \ref{lem-prop2.14-bundles} may seem a bit
  stilted at first glance.  However, they are precisely what are
  needed to handle induced bundle representations of groupoids.  We
  will use this result is the situation where $X$ is a principal
  $G$-space, $Z=\xssx$, $H$ is the associated imprimitivity groupoid
  $Y=\ho$, $\sigma([x,y])=r_{X}(x)$ and $q=r_{H}$.
\end{remark}

\section{Groupoid Crossed Products}
\label{sec:group-cross-prod}

In this section we want to review what it means for a locally
Hausdorff, locally compact groupoid $G$ to act on a
$C_{0}(\go)$-algebra $\A$ by isomorphisms.  Such actions will be
called \emph{groupoid dynamical systems} and will be denoted $\aga$.
We also discuss the associated crossed product $\acg$.  Fortunately,
there are several nice treatments in the literature upon which one can
draw \citelist{\cite{gal:94}*{\S7} \cite{pat:groupoids99}
  \cite{gal:kt99}\cite{gal:cm01}*{\S2}
  \cite{khoska:jram02}*{\S1}\cite{khoska:jot04}*{\S2.4 \& \S3}
  \cite{ren:jot87}}. However, as in \cite{gal:94}*{\S7}, we intend to
emphasize the underlying bundle structure.  Otherwise, our treatment
follows the excellent exposition in
\citelist{\cite{khoska:jram02}\cite{khoska:jot04}}.  We remark that
Renault uses the more restrictive definition of \cs-bundle in
\cite{ren:jot87}.  In our formulation, the groupoid analogue of a
strongly continuous group of automorphisms of a \cs-algebra arises
from a certain type of action of the groupoid on the total space of an
\usc-\cs-bundle.
\begin{definition}
  \label{def-dynamical-sys}
  Suppose that $G$ is a locally Hausdorff, locally compact groupoid
  and that $A$ is a $C_{0}(\go)$-algebra such that $A=\sa_{0}(\go;\A)$
  for an \usc-\cs-bundle $\A$ over $\go$.  An action $\alpha$ of $G$
  on $A$ by $*$-isomorphisms is a family
  $\set{\alpha_{\gamma}}_{\gamma\in G}$ such that
  \begin{enumerate}
  \item for each $\gamma$, $\alpha_{\gamma}:A\bigl(s(\gamma)\bigr)\to
    A\bigl(r(\gamma)\bigr)$ is an isomorphism,
  \item $\alpha_{\eta\gamma}=\alpha_{\eta}\circ\alpha_{\gamma}$ for
    all $(\eta,\gamma)\in G^{(2)}$ and
  \item $\gamma\cdot a:=\alpha_{\gamma}(a)$ defines a continuous
    action of $G$ on $\A$.
  \end{enumerate}
  The triple $\aga$ is called a (groupoid) dynamical system.
\end{definition}

Our next lemma implies that our definition coincides with that in
\cite{khoska:jot04} (where the underlying bundle structure is not
required), but first we insert a remark to help with the notation.

\begin{remark}
  \label{restriction-of-pullbacks}
  Suppose that $\aga$ is a dynamical system, and let $A$ be the
  $C_{0}(\go)$-algebra, $\sa_{0}(\go;\A)$.  We may pull back $\A$ to
  $G$ with $s$ and $r$ to get two \usc-\cs-bundles $s^*\A$ and
  $r^*\A$ on $G$ (See Remark \ref{rem-sections-of-pullbacks}).  If
  $U$ is a subset of $G$, we may restrict these bundles to $U$,
  getting bundles on $U$ which we denote by ${s\restr U}^*\A$ and
  ${r\restr U}^*\A$.  If $U$ is open and Hausdorff in $G$, then we
  may form the $C_{0}(U)$-algebras, $\sa_{0}(U;s^*\A)$ and
  $\sa_{0}(U;r^*\A)$, which we denote by ${s\restr U}^*(A)$ and
  ${r\restr U}^*(A)$, respectively.  Then, by the discussion in the
  two paragraphs preceding Lemma \ref{lem-top-sections}, we see that
  there bijective correspondence between bundle isomorphisms between
  ${s\restr U}^*\A$ and ${r\restr U}^*\A$ and
  $C_0(U)$-isomorphisms between ${s\restr U}^*(A)$ and ${r\restr
    U}^*(A)$.  (See also the comments prior to
  Remark~\vref{rem-sections-of-pullbacks}.) 
\end{remark}

\begin{lemma}
  \label{lem-same-def}
  Suppose that $\aga$ is a dynamical system and let $A$ be the
  $C_{0}(\go)$-algebra, $\sa_{0}(\go;\A)$.  If $U\subset G$ is open
  and Hausdorff, then
  \begin{equation}\label{eq:52p}
    \alpha_{U}(f)(\gamma):= \alpha_{\gamma}\bigl(f(\gamma)\bigr)
  \end{equation}
  defines a $C_{0}(U)$-isomorphism of ${s\restr U}^{*}(A)$ onto
  ${r\restr U}^{*}(A)$.  If $V\subset U$ is open, then viewing
  ${s\restr V}^{*}(A)$ as an ideal in ${r\restr U}^{*}(A)$,
  $\alpha_{V}$ is the restriction of $\alpha_{U}$.

  Conversely, if $A$ is a $C_{0}(\go)$-algebra and if for each open,
  Hausdorff subset $U\subset G$, there is a $C_{0}(U)$-isomorphism
  $\alpha_{U}:{s\restr U}^{*}(A)\to{r\restr U}^{*}A$ such that
  $\alpha_{V}$ is the restriction of $\alpha_{U}$ whenever $V\subset
  U$, then there are well-defined isomorphisms
  $\alpha_{\gamma}:A\bigl(s(\gamma)\bigr)\to A\bigl(r(\gamma)\bigr)$
  satisfying \eqref{eq:52p}.  If in addition,
  $\alpha_{\eta\gamma}=\alpha_{\eta}\circ \alpha_{\gamma}$ for all
  $(\eta,\gamma)\in G^{(2)}$, then $\aga$ is a dynamical system.
\end{lemma}
\begin{proof}
  If $\aga$ is a dynamical system, then the statements about the
  $\alpha_{U}$ are easily verified.

  On the other hand, if the $\alpha_{U}$ are as given in the second
  part of the lemma, then the
  $\alpha_{\gamma}:=(\alpha_{U})_{\gamma}:A\bigl(s(\gamma)\bigr) \to
  A\bigl(r(\gamma)\bigr)$ are uniquely determined due to the
  compatibility condition on the $\alpha_{U}$'s.  It only remains to
  check that $\gamma\cdot a:=\alpha_{\gamma}(a)$ defines a continuous
  action of $G$ on $\A$.  Suppose that $(\gamma_{i},a_{i})\to
  (\gamma,a)$ in $\set{(\gamma,a):s(\gamma)=p(a)}$.  We need to prove
  that $\gamma_{i}\cdot a_{i}\to \gamma\cdot a$ in $\A$.  We can
  assume that there is a Hausdorff neighborhood $U$ of $\gamma$
  containing all $\gamma_{i}$.  Let $g\in{s\restr U}^{*}(A)$ be such
  that $g(\gamma)=a$.  We have
  \begin{equation*}
    \alpha_{U}(g)(\gamma_{i})\to \alpha_{U}(g)(\gamma):=\gamma\cdot a. 
  \end{equation*}
  Also
  \begin{equation*}
    \|\alpha_{U}(g)(\gamma_{i})-\gamma_{i}\cdot a_{i}\| = \|
    \alpha_{\gamma_{i}}\bigl( g(\gamma_{i})-a_{i}\bigr) \| =
    \|g(\gamma_{i})-a_{i}\| \to 0.
  \end{equation*}
  Therefore $\alpha_{U}(g)(\gamma_{i})-\gamma_{i}\cdot a_{i}\to
  0_{p(a)}$, and
  \begin{equation*}
    \gamma_{i}\cdot a_{i} =\alpha_{U}(g)(\gamma_{i})+
    \bigl(\gamma_{i}\cdot a_{i} - \alpha_{U}(g)(\gamma_{i})\bigr) \to
    \gamma\cdot a + 0_{p(a)}=\gamma\cdot a.\qed
  \end{equation*}
  \renewcommand{\qed}{}
\end{proof}

In the Hausdorff case, the crossed product is a completion of the
compactly supported sections $\sa_{c}(G;r^{*}\A)$ of the pull-back of
$\A$ via the range map $r:G\to\go$.  In the non-Hausdorff case, we
must find a substitute for $\sa_{c}(G;r^{*}\A)$.  As in
Section~\ref{sec:cox-algebras}, we let $\ccra$ be
the subspace of functions from $G$ to $\A$ spanned by elements in
$\sa_{c}(U;{r}^{*}\A)$ for all open, Hausdorff sets $U\subset G$.
(Elements in $\sa_{c}(U;{r}^{*}\A)$ are viewed as functions on $G$ as
in the definition of $\cc(G)$.)
\begin{prop}
  \label{prop-star-alg}
  If $G$ is a locally Hausdorff, locally compact groupoid with Haar
  system $\set{\lambda^{u}}$, then $\ccra$ becomes a $*$-algebra with
  respect to the operations
  \begin{equation*}
    f*g(\gamma)=\int_{G}f(\eta)\alpha_{\eta}\bigl(g(\eta^{-1}\gamma)\bigr) 
    \,d \lambda^{r(\gamma)} (\eta) \quad\text{and}\quad
    f^{*}(\gamma)=\alpha_{\gamma}\bigl(f(\gamma^{-1})^{*}\bigr).
  \end{equation*}
\end{prop}
The proof of the proposition is fairly routine, the only real issue
being to see that the formula for $f*g$ defines an element of $\ccra$.
For this, we need a preliminary observation.
\begin{lemma}
  \label{lem-cty-int}
  Suppose that $U$ and $W$ are Hausdorff open subsets of $G$.  Let
  $U\starr W=\set{(\eta,\gamma):r(\eta)=r(\gamma)}$, and let
  $r^{*}\A=\set{(\eta,\gamma,a):r(\eta)=r(\gamma)=p(a)}$ be the
  pull-back.  If $F\in\sa_{c}(U\starr W;r^{*}\A)$, then
  \begin{equation*}
    f(\gamma)=\int_{G} F(\eta,\gamma)\,d\lambda^{r(\gamma)}(\eta)
  \end{equation*}
  defines a section in $\sa_{c}(W;{r}^{*}\A)$.
\end{lemma}
\begin{proof}
  If $F_{i}\to F$ in the inductive limit topology on $\sa_{c}(U\starr
  W;r^{*}\A)$, then it is straightforward to check that $f_{i}\to f$
  in the inductive limit topology on $\sa_{c}(W;{r}^{*}\A)$.  Thus it
  will suffice to consider $F$ of the form $(\eta,\gamma)\mapsto
  h(\eta,\gamma)a\bigl(r(\gamma)\bigr)$ for $h\in C_{c}(U\starr W)$
  and $a\in\sa_{c}(\go;\A)$.  Since $U\starr W$ is closed in $U\times
  W$, we can assume that $h$ is the restriction of $H\in C_{c}(U\times
  W)$.  Since sums of the form $(\eta,\gamma)\mapsto
  h_{1}(\eta)h_{2}(\gamma)$ are dense in $C_{c}(U\starr W)$, we may as
  well assume that $H(\eta,\gamma)=h_{1}(\eta)h_{2}(\gamma)$ for
  $h_{1}\in C_{c}(U)$ and $h_{2}\in C_{c}(W)$.  But then
  \begin{equation*}
    f(\gamma)=h_{2}(\gamma)a\bigl(r(\gamma)\bigr)
    \int_{G}h_{1}(\eta)\,d\lambda^{r(\gamma)}(\eta), 
  \end{equation*}
  which is clearly in $\sa_{c}(W;{r}^{*}\A)$.
\end{proof}
\begin{proof}[Proof of Proposition~\vref{prop-star-alg}]
  To see that convolution is well-defined, it suffices to see that
  $f*g\in\ccra$ when $f\in \sa_{c}(U;{r}^{*}\A)$ and $g\in
  \sa_{c}(V;{r}^{*}\A)$ for Hausdorff open sets $U$ and $V$.  We
  follow the argument of \cite{khoska:jram02}*{p.~52--3}.  
Since $U$ is Hausdorff, and therefore regular, there is an open set
$U_{0}$ and a compact set $K_{f}$ such that
\begin{equation*}
  \supp f\subset U_{0}\subset K_{f} \subset U.
\end{equation*}
Given $\gamma\in \supp g$, we have $K_{f}\gamma\gamma^{-1}\subset U$.
Even if $K_{f}\gamma\gamma^{-1}$ is empty, using the local compactness
of $G$ and the continuity of multiplication, we can find an open set
$W_{\gamma}$ and a compact set $K_{\gamma}$ such that $\gamma\in
W_{\gamma}\subset K_{\gamma}$ with
\begin{equation*}
  K_{f}K_{\gamma}K_{\gamma}^{-1}\subset U.
\end{equation*}
\begin{claim}
  $U_{0}W_{\gamma}$ is Hausdorff.
\end{claim}
\begin{proof}[Proof of the Claim]
  Suppose that $\eta_{i}\gamma_{i}$ converges to both $\alpha$ and
  $\beta$ in $U_{0}W_{\gamma}$.  Since $W_{\gamma}\subset K_{\gamma}$,
  we can pass to a subnet, relabel, and assume that
  $\gamma_{i}\to\gamma\in K_{\gamma}$.  Then $\eta_{i}$ converges to
  both $\alpha\gamma^{-1}$ and $\beta\gamma^{-1}$.  This the latter
  are both in $U_{0}W_{\gamma}K_{\gamma}^{-1}\subset
  K_{f}K_{\gamma}K_{\gamma}^{-1}\subset U$, and since $U$ is
  Hausdorff, we must have $\alpha\gamma^{-1}=\beta\gamma^{-1}$.  But
  then $\alpha=\beta$.  This proves the claim.
\end{proof}

Since $\supp g$ is compact, we can find open sets $U_{1},\dots, U_{n}$
and $W_{1},\dots,W_{n}$ such that $\supp g\subset \bigcup W_{i}$,
$\supp f\subset U_{i}$ and $U_{i}W_{i}$ is Hausdorff.  If we let
$U':=\bigcap U_{i}$ and use a partition of unity to write $g=\sum g_{i}$
with each $\supp g_{i}\subset V_{i}$, then we can view $f\in
C_{c}(U')$ and replace $g$ by $g_{i}$.  Thus we may assume that $f\in
\sa_{c}(U;{r}^{*}\A)$ and $g\in \sa_{c}(V;{r}^{*}\A)$ for Hausdorff
open sets $U$ and $V$ with $UV$ Hausdorff as well.  Next observe
that $(\eta,\gamma)\mapsto (\eta,\eta\gamma)$ is a homeomorphism of
$B:=\set{(\eta,\gamma)\in U\times V:s(\eta)=r(\gamma)}$ onto an
\emph{open} subset $B'$ of $U*_{r}UV$.  On $B'$, we can define a
continuous function with compact support by
  \begin{equation*}
    (\eta_{1},\gamma_{1})\mapsto f(\eta_{1})\alpha_{\eta_{1}}
    \bigl(g(\eta_{1}^{-1}\gamma_{1})\bigr) .
  \end{equation*}
  Extending this function to be zero off the open subset $B'$, we get
  a section (as in Remark~\vref{rem-sections-of-pullbacks})
  $\phi\in\sa_{c}(U*_{r} UV;r^{*}\A)$ such that
  \begin{equation*}
    \phi(\eta_{1},\gamma_{1})=f(\eta_{1})g(\gamma_{1})\quad\text{for
      all $(\eta_{1},\gamma_{1})\in B$.}
  \end{equation*}
  It follows from Lemma~\vref{lem-cty-int}, that $f*g\in
  \sa_{c}(UV;r^{*}\A) \subset \sacc(G;r^{*}\A)$.

  The remaining assertions required to prove the proposition are
  routine to verify.
\end{proof}

If $f\in\ccra$, then $\gamma\mapsto \|f(\gamma)\|$ is
upper-semicontinuous on open Hausdorff subsets, and is therefore
integrable on $G$ with respect to any Radon measure.  Thus we can
define the \emph{$I$-norm} by
\begin{equation*}
  \|f\|_{I}=\max\set{\sup_{u\in\go}
    \int_{G}\|f(\gamma)\|\,d\lambda^{u}(\gamma), \sup_{u\in\go}
    \int_{G}\|f(\gamma)\|\,d\lambda_{u}(\gamma)}.  
\end{equation*}
The \emph{crossed product} $\acg$ is defined to be the enveloping
\cs-algebra of $\ccra$.  Specifically, we define the (universal)
\cs-norm by
\begin{equation*}
  \|f\|:=\sup\{\,\|L(f)\|:\text{$L$ is a $\|\cdot\|_{I}$-decreasing
    $*$-representation of $\ccra$}\,\}.
\end{equation*}
Then $\acg$ is the completion of $\ccra$ with respect to
$\|\cdot\|$. (The
notation $A\rtimes_{\alpha}G$ would also be appropriate; it is used in
\cite{khoska:jot04} for example.  Since our approach is a bundle one,
using a notation that includes the bundle seems appropriate.)

\begin{example}
  \label{ex-new-scalar-case}
  As in the case of ordinary \cs-dynamical systems (see \cite{wil:crossed}*{Example~2.33}), the a groupoid
  \cs-algebra $\cs(G,\lambda)$ is a degenerate case of the groupoid
  crossed product. Let $G$ be a \lhlc{} groupoid with Haar system
  $\set{\lambda^{u}}_{u\in\go}$.  Let $\A=\trivial{\go}$ be the
  trivial bundle $\go\times\C$.  Then $G$ acts by isomorphisms on
  $\trivial{\go}$ by left translation:
  \begin{equation*}
    \lt^{G}_{\gamma}\bigl(s(\gamma),z\bigr)=\bigl(r(\gamma),z\bigr). 
  \end{equation*}
Then it is routine to check that $\bigl(\trivial{\go},G,\lt^{G}\bigr)$
is a dynamical system with $\trivial{\go}\rtimes_{\lt^{G}}G$
isomorphic to the groupoid \cs-algebra $\cs(G,\lambda)$.
\end{example}

Since a group is a groupoid whose unit space is a single point, we can
view ordinary dynamical systems and their crossed products as trivial
examples of groupoid dynamical systems and crossed products.   A more
interesting class examples arise as follows.

\begin{example}
  \label{ex-new-group-dynamical}
  Suppose that $A$ is a \cox-algebra and that $X$ is a locally compact
  (Hausdorff) 
  $\gG$-space for a locally compact group $\gG$.  Suppose that
  \begin{equation*}
    \galpha:\gG\to\Aut A
  \end{equation*}
is a strongly continuous automorphism group such that
\begin{equation}\label{eq:70}
  \galpha_{s}(\phi\cdot a)=\lt_{s}(\phi)\cdot \galpha_{s}(a),
\end{equation}
where $\lt_{s}(\phi)(x)=\phi(s^{-1}\cdot x)$.  For example, if $\Prim
A$ is Hausdorff, then we can let $X:=\Prim A$.  Then $A$ is a
\cox-algebra (via the Dauns-Hofmann Theorem), and $X$ is naturally
a $\gG$-space such that \eqref{eq:70} holds (see
\cite{wil:crossed}*{Lemma~7.1}). 

Let $G=\gG\times X$ be the transformation groupoid, and note that
$A=\sa_{0}(X;\A)$ for an \usc-bundle $\A$
\cite{wil:crossed}*{Theorem~C.26}.  Then we get a groupoid dynamical
system $(\A,G,\alpha)$ where
\begin{equation}
  \label{eq:71}
  \alpha_{(s,x)}\bigl(a(s^{-1}\cdot x)\bigr)=\galpha_{s}(a)(x).
\end{equation}
Then it is a matter of checking that $\A\rtimes_{\alpha}G$ is
isomorphic to the ordinary crossed product $A\rtimes_{\galpha}\gG$ via
the map that sends $f\in C_{c}(G,A)$ to $\tilde f\in
\sa_{c}(G;r^{*}\A)$ given by
\begin{equation}\label{eq:73}
  \tilde f(s,x)=\Delta(s)^{\frac12}f(s)(x),
\end{equation}
where $\Delta$ is the modular function on $\gG$.\footnote{The modular
  function is introduced simply because it is traditional to use the
  modular function as part of the definition of the involution on
  $C_{c}(\gG,A)\subset A\rtimes_{\galpha}\gG$ and we need $f\mapsto
  \tilde f$ to be $*$-preserving.  The indicated map on
  dense subalgebras is isometric because representations which are
  continuous in the inductive limit topology are in fact bounded with
  respect to the universal norms.  For the same reason, it is
  possible, although not traditional, to define the involution on
  ordinary crossed products without the modular function.  In the
  latter case, we could dispense with the modular function in
  \eqref{eq:73}.} 
For example,
\begin{align*}
  \Delta(s)^{\frac12}f*g(s)(x)&=\int_{\gG}
  \Delta(r)^{\frac12}
  f(r)\galpha_{r}\bigl(\Delta(r^{-1}s)^{\frac12}g(r^{-1}s)\bigr) \,ds\,(x) \\ 
  &= \int_{\gG} \Delta(r)^{\frac12} f(r)(x)
  \galpha_{r}\bigl(\Delta(r^{-1}s)^{\frac12}g(s^{-1}r)\bigr) (x) \, ds \\ 
  &= \int_{G} \tilde f(r,x) \alpha_{(s,r)} \bigl(\tilde g(s^{-1}r,
  s^{-1}\cdot x)\bigr) \, ds \\
  &=\tilde f*\tilde g(s,x).
\end{align*}
\end{example}
\section{Renault's Equivalence Theorem}
\label{sec:equiv-group-dynam}

In this section, we want to extend Renault's definition
\cite{ren:jot87}*{Definition~5.3} of an equivalence between two
dynamical systems $(H,\B,\beta)$ and $\aga$ to the setting of
\usc-\cs-bundles, and to give a precise statement of the his
Equivalence Theorem in this context.  In doing so, we also give an
explicit description of the pre-\ib{} between $\sacc(H;r^{*}\B)$ and
$\sacc(G;r^{*}\A)$.

\begin{definition}
  \label{def-equi-dyn-sys}
  An \emph{equivalence} between dynamical systems $(\B,H,\beta)$ and
  $\aga$ is an \usc-Banach bundle $p_{\E}:\E\to X$ over a
  $(H,G)$-equivalence $X$ together with $B\bigl(r(x)\bigr) \sme
  A\bigl(s(x)\bigr)$-\ib{} structures on each fibre $\E_{x}$ and
  commuting (continuous) $H$- and $G$-actions on the left and right,
  respectively, of $\E$ such that the following additional properties
  are satisfied.
  \begin{enumerate}
  \item (Continuity) The maps induced by the \ib{} inner products from
    $\E*\E\to \B$ and $\E*\E\to \A$ are continuous as are the maps
    $\B*\E\to \E$ and $\E*\A\to\E$ induced by the \ib{}
    actions.\label{item:2}
  \item (Equivariance) The groupoid actions are equivariant with
    respect to the bundle map $p_{\E}:\E\to X$; that is,  $p_{\E}(\eta\cdot
    e)=\eta\cdot p_{\E}(e)$ and  $p_{\E}(e\cdot \gamma)=p_{\E}(e)\cdot
    \gamma$.
  \item (Compatibility) The groupoid actions are compatible with the
    \ib{} structure:
    \begin{align*}
      \lip\B<\eta\cdot e,\eta\cdot f>&=\beta_{\eta}\bigl(\lip\B<e,f>
      \bigr) & \eta\cdot (b\cdot e)&=\beta_{\eta}(b)\cdot (\eta\cdot e)\\
      \rip\A<e\cdot\gamma,f\cdot\gamma>&=\alpha^{-1}_{\gamma}\bigl(
      \rip\A<e,f> \bigr) &(e\cdot a)\cdot \gamma&=(e\cdot \gamma)\cdot
      \alpha_{\gamma}^{-1}(a).
    \end{align*}
  \item (Invariance) The $H$-action commutes with the $\A$-action on
    $\E$ and the $G$-action commutes with the $\B$-action.  That is,
    $\eta\cdot (e\cdot a)=(\eta\cdot e)\cdot a$ and $(b\cdot e)\cdot
    \gamma =b\cdot (e\cdot \gamma)$.\label{item:1}
  \end{enumerate}
\end{definition}

\begin{lemma}
  \label{lem-ip-invar}
  As a consequence of invariance we have
  \begin{equation*}
    \lip\B<e\cdot \gamma,f\cdot \gamma>=\lip\B<e,f>\quad\text{and}
    \quad \rip\A<\eta\cdot e,\eta\cdot f>=\rip\A<e,f>
  \end{equation*}
  for all $e,f\in\E$, $\eta\in H$ and $\gamma\in G$.
\end{lemma}
\begin{proof}
  If $g\in\E$, then using invariance, we have
  \begin{align*}
    \lip\B<e,f>\cdot(g\cdot\gamma)&= \bigl(\lip\B<e,f>\cdot
    g\bigr)\cdot\gamma \\
    &= \bigl(e\cdot \rip\A<f,g>\bigr)\cdot\gamma\\
    &=(e\cdot \gamma)\cdot\rip\A<f\cdot\gamma,g\cdot \gamma> \\
    &=\lip\B<e\cdot \gamma,f\cdot\gamma>\cdot (g\cdot\gamma).
  \end{align*}
  The first equation follows and the second follows by symmetry.
\end{proof}

\begin{remark}
  Since our inner products are full, the converse of
  Lemma~\vref{lem-ip-invar} holds as well.  That is, if the inner
  products are invariant under the ``other'' groupoid action, then
  invariance holds.
\end{remark}

\begin{example}
  \label{ex-basic-equivalence}
  An important and instructive example of
  Definition~\vref{def-equi-dyn-sys} is to see that $(\A,G,\alpha)$ is
  equivalent to itself via $p:r^{*}\A\to G$.  Recall that
  \begin{equation*}
    r^{*}\A:=\set{(\gamma,a)\in
      G\times\A:r(\gamma)=p_{\A}(a)}.
  \end{equation*}
  We equip the fibre over $\gamma$ with a $A\bigl(r(\gamma)\bigr) \sme
  A\bigl(s(\gamma)\bigr)$-\ib{} structure as follows:
  \begin{align*}
    \blip A(r(\gamma))<(\gamma,a),(\gamma,b))> &:=
    ab^{*} &a\cdot(\gamma,b) &:=(\gamma,ab) \\
    \brip A(s(\gamma))<(\gamma,a),(\gamma,b))> &:=
    \alpha_{\gamma}^{-1}(a^{*}b) & (\gamma,b)\cdot a &=
    (\gamma,b\alpha_{\gamma}(a)).
  \end{align*}
  We let $G$ act on the right and left of $r^{*}\A$ as follows:
  \begin{align*}
    \beta\cdot(\gamma,a)&:= (\beta\gamma,ab)& (\gamma,a)\cdot\beta &=
    (\gamma\beta, a).
  \end{align*}
  At this point, it is a simple matter to verify that axioms
  \partref1--\partref4 of Definition~\vref{def-equi-dyn-sys} are
  satisfied.
\end{example}

\begin{thm}[\cite{ren:jot87}*{Corollaire~5.4}]
  \label{thm-renault}
Suppose that $G$ and $H$ are second countable \lhlc{} groupoids with
Haar systems $\set{\lambda_{G}^{u}}_{u\in\go}$ and
$\set{\lambda_{H}^{v}}_{v\in\ho}$, respectively.
  If $p_{\E}:\E\to X$ is an equivalence between $(\B,H,\beta)$ and
  $\aga$, then $\X_{0}=\sacc(X;\E)$ becomes a
  $\cp(\B,H,\beta)\sme\acg$-pre-\ib{} with respect to the following
  operations:
  \begin{align}
    \bchip<z,w>(\eta)&:=\int_{G}\blip\B<z(\eta\cdot x \cdot \gamma) ,
    \eta\cdot w(x\cdot \gamma)>\,d\lambda_{G}^{s(x)}(\gamma), \label{eq:1}\\
    f\cdot z (x) &:= \int_{H} f(\eta)\cdot \bigl(\eta\cdot
    z(\eta^{-1}\cdot x)\bigr)\,d\lambda_{H}^{r(x)}(\eta),\label{eq:3}\\
    z\cdot g(x)&:= \int_{G} \bigl(z(x\cdot \gamma)\cdot
    \gamma^{-1}\bigr) \cdot \alpha_{\gamma}\bigl(g(\gamma^{-1})\bigr)
    \,d\lambda_{G}^{s(x)}(\gamma) \text{ and}\label{eq:4}\\
    \acgip<w,v>(\gamma)&:=\int_{H} \brip\A<w(\eta^{-1}\cdot y \cdot
    \gamma^{-1}) , v(\eta^{-1}\cdot
    y)\cdot\gamma^{-1}>\,d\lambda_{H}^{r(y)}(\eta). \label{eq:5}
  \end{align}
\end{thm}
\begin{remark}
  \label{rem-ind-x}
  Since $X$ is a $(H,G)$-equivalence, the equation $r(x)=r(y)$ implies
  that $y=x\cdot \gamma'$ for some $\gamma'\in G$.  Thus in
  \eqref{eq:1} we are free to choose any $x\in
  r_{X}^{-1}\bigl(s_{H}(\eta)\bigr)$.  On the other hand, we can replace
  $x$ in \eqref{eq:1} by $y:= \eta\cdot x$ and obtain
  \begin{equation}
    \label{eq:2}
    \bchip<z,w>(\eta):=\int_{G}\blip\B<z(y \cdot \gamma) ,
    \eta\cdot w(\eta^{-1}\cdot y\cdot
    \gamma)>\,d\lambda_{G}^{s(y)}(\gamma), \tag{$\ref{eq:1}'$}
  \end{equation}
where any $y\in r_{X}^{-1}\bigl(r_{G}(\eta)\bigr)$ will do.
  Similarly, in \eqref{eq:5}, we are free to choose any $y\in
  s_{X}^{-1}\bigl(s_{G}(\gamma)\bigr)$.
\end{remark}

\begin{remark}
  \label{rem-well-defined}
  Checking that \eqref{eq:1}--\eqref{eq:5} take values in the
  appropriate spaces of functions is a bit fussy in the non-\ttwo{}
  case.  We can suppose that $z\in\sa_{c}(U;\E)$ and
  $w\in\sa_{c}(V;\E)$, where $U$ and $V$ are \ttwo{} open subsets of
  $X$.  Then $U*_{s}V$ is a \ttwo{} open subset of $\xssx$.  Let
  $q:\xssx\to H$ be the ``orbit'' map (cf.\
  Lemma~\vref{lem-fund-groupoid-lemma}).  We get an element $f\in
  \sa_{c}(U*_{s}V;(r_{H}\circ q)^{*}\B)$ defined by
  \begin{equation*}
    f(x,y):=\blip\B<z(x),\tau(x,y)\cdot w(y)>,
  \end{equation*}
  where $\tau(x,y)$ is defined as in
  Lemma~\ref{lem-fund-groupoid-lemma}.  Then the obscure hypotheses of
  Lemma~\vref{lem-prop2.14-bundles} have been cooked up so that we can
  conclude that
  \begin{equation*}
    \bchip<z,w>(\sigma)=\lambda(f)\bigl(q(\sigma\cdot x,x)\bigr)
  \end{equation*}
  defines an element in $\sa_{c}(q(U*_{s}V);r_{H}^{*}\B)$ as
  required.

  To see that \eqref{eq:3} defines an element of $\sacc(X;\E)$, we
  proceed exactly as in the proof for the convolution in
  Proposition~\vref{prop-star-alg}.  We assume
  $f\in\sa_{c}(V;r^{*}\B)$ and $z\in \sa_{c}(U;\E)$.  Using partitions
  of unity, we can assume that the open set $V\cdot U$ is \ttwo.  The
  map $(\eta,x)\mapsto (\eta,\eta\cdot x)$ is a homeomorphism of
  $B=\set{(\eta,x)\in V\times U:s(\eta)=r(x)}$ onto an open subset
  $B'$ of $V*_{r}V\cdot U=\set{(\sigma,y)\in V\times V\cdot
    U:r(\sigma)=r(y)}$.  Hence the integrand in \eqref{eq:3} is a
  section $h\in\sa_{c}(V*_{r}V\cdot U;r^{*}\E)$.  An argument
  analogous to that in Lemma~\vref{lem-cty-int} shows that $f\cdot
  z\in \sacc(X;\E)$.
\end{remark}

\begin{remark}
  \label{rem-basic-ex}
  It is worth noting that, in Example~\vref{ex-basic-equivalence}, the
  inner-products and actions set out in \eqref{eq:1}--\eqref{eq:5} are
  the natural ones:
  \begin{align*}
    \bchip<z,w>&=z*w^{*}&f\cdot z&=f*z \\
    \acgip<w,v>&=w^{*}*v& z\cdot g&=z*g.
  \end{align*}
  For example, we start with \eqref{eq:1}:
  \begin{align*}
    \bchip<z,w>(\eta)&=\int_{G} \blip A(r(\gamma))<z(\eta\cdot x\cdot
    \gamma), \eta\cdot w(x\cdot \gamma)>\,d\lambda_{G}^{s(x)}(\gamma) \\
    &=\int_{G}z(\eta\gamma)\eta\cdot w(\gamma)^{*}
    \,d\lambda_{G}^{s(\eta)}(\gamma) \\
    &=\int_{G} z(\gamma)\alpha_{\eta}\bigl( w
    (\eta^{-1}\gamma)^{*}\bigr)
    \,d\lambda_{G}^{r(\eta)}(\gamma) \\
    &= \int_{G} z(\gamma)
    \alpha_{\gamma}\bigl(w^{*}(\gamma^{-1}\eta)\bigr)
    \,d\lambda_{G}^{r(\eta)}(\gamma) \\
    &=z*w^{*}(\eta)
  \end{align*}
  Similarly, we can start with \eqref{eq:3}:
  \begin{align*}
    f\cdot z(x)&= \int_{G}f(\eta)\cdot\bigl(\eta\cdot z(\eta^{-1}\cdot
    x)\bigr) \,d\lambda_{G}^{r(x)}(\eta) \\
    &=\int_{G}f(\eta)\alpha_{\eta}\bigl(z(\eta^{-1}x)\bigr)
    \,d\lambda_{G}^{r(x)} (\eta) \\
    &= f*z(x).
  \end{align*}
  We can do the same with \eqref{eq:4} and \eqref{eq:5}, or appeal to
  symmetry (as described below).
\end{remark}

It will be helpful to see that equivalence of dynamical systems is
completely symmetric.  Let $\E$ be an equivalence between
$(\B,H,\beta)$ and $\aga$.  Let $\Es$ be the underlying topological
space of $\E$, let $\dual:\E\to\Es$ be the identity map and define
$\ps:\Es\to\Xop$ by
$\ps\bigl(\dual(e)\bigr)=\bigl(p(e)\bigr)^{\text{\normalfont op}}$.
Then as a Banach space, the fibre $\Es_{\op{x}}=\E_{x}$ and we can
give $\Es_{\op{x}}$ the dual $A\bigl(r(\op x)\bigr)\sme B\bigl(s(\op
x)\bigr)$-\ib{} structure of the dual module $(\E_{x})^{*}$.  Then
$\ps:\Es\to\Xop$ is a $\aga\sme(\B,H,\beta)$ equivalence.
Furthermore, if we define $\Phi:\sacc(X;\E)\to\sacc(\Xop;\Es)$ by
$\Phi(f)(\op x):=\dual\bigl(f(x)\bigr)$, then we can easily compute
that
\begin{align*}
  \btlip \cp(\A,G,\alpha)<\Phi(z),\Phi(w)>(\gamma) &= \int_{H}
  \blip\A<\Phi(z)(\gamma\cdot \op x \cdot \eta), \gamma\cdot
  \Phi(w)(\op x\cdot \eta)> \,d\lambda_{H}^{s(\op x)}(\eta) \\
  &=\int_{H} \blip\A<\dual\bigl(z(\eta^{-1}\cdot x\cdot
  \gamma^{-1})\bigr), \gamma\cdot \dual\bigl(w(\eta^{-1}\cdot
  x)\bigr)> \,d\lambda_{H}^{r(x)}\\
  &=\int_{G}\blip\A<\dual\bigl(z(\eta^{-1}\cdot x\cdot
  \gamma^{-1}))\bigr) , \dual\bigl(w(\eta^{-1}\cdot x)\cdot
  \gamma^{-1}\bigr) > \,d\lambda_{H}^{r(x)}
  \\
  &= \brip\A<z(\eta^{-1}\cdot x\cdot \gamma^{-1})) , w(\eta^{-1}\cdot
  x)\cdot
  \gamma^{-1} >  \,d\lambda_{H}^{r(x)} \\
  &= \acgip<z,w>.
\end{align*}
Equally exciting computations give us the following lemma.
\begin{lemma}
  \label{lem-sym}
  With $\Phi:\sacc(X;\E)\to\sacc(\Xop;\Es)$ defined as above, we have
  \begin{align*}
    \btlip \cp(\A,G,\alpha)<\Phi(z),\Phi(w)>&= \acgip<z,w> & \btrip
    \cp(\B,H,\beta )<\Phi(w),\Phi(v)> &= \bchip<w,v> \\
    g\cdot \Phi(z)&= \Phi(z\cdot g^{*}) & \Phi(z)\cdot f &= \Phi(f^{*}
    \cdot z).
  \end{align*}
\end{lemma}

This lemma can be very useful.  For example, once we show that
$\btlip\cp(\B,H,\beta)<z,z>$ is positive for all $z$, it follows by
symmetry that $\btlip \cp(\A,G,\alpha)<\Phi(z),\Phi(z)>$ is positive
for all $z$.
But by Lemma~\vref{lem-sym}, we must have $\acgip<z,z>$ positive.

Now for example, we show that the left-inner product respects
left-module action:
\begin{align*}
  \bchip<f\cdot z,w>(\eta)&= \int_{G} \blip\B<f\cdot z(\eta\cdot
  x\cdot \gamma), \eta\cdot w (x\cdot \gamma)>
  \,d\lambda_{G}^{s(x)}(\gamma) \\
  &= \int_{H}\int_{G} \blip\B<f(\sigma)\cdot\bigl(\sigma\cdot z\cdot
  (\sigma^{-1} \eta\cdot x\cdot \gamma)\bigr),\eta\cdot w(x\cdot
  \gamma)> \\
  &\hskip
  40ex\lambda_{G}^{s(x)}(\gamma)\,d\lambda_{H}^{r(\eta)}(\sigma)
  \\
  &= \int_{H} f(\sigma) \beta_{\sigma} \Bigl( \int_{G} \blip\B
  <z(\sigma^{-1}\eta\cdot x \cdot \gamma), \sigma^{-1}\eta\cdot
  w(x\cdot \gamma) >\,d\lambda_{G}^{s(x)}(\gamma) \Bigr)
  \\
  &\hskip55ex d\lambda_{H}^{r(\eta)} (\sigma) \\
  &= \int_{H}f(\sigma)
  \beta_{\sigma}\Bigl(\bchip<z,w>(\sigma^{-1}\gamma)\Bigr)
  \,d\lambda_{H}^{r(\eta)}(\sigma) \\
  &= f*\bchip<z,w>(\eta).
\end{align*}

By symmetry and by applying Lemma~\vref{lem-sym}, it also follows that
\begin{equation*}
  \acgip<w,v\cdot g>=\acgip<w,v>*g.
\end{equation*}
Similar computations show that \eqref{eq:3} defines a left-action, and
it automatically follows that \eqref{eq:4} is a right action by
symmetry.

Next we check that
\begin{equation}
  \label{eq:6}
  \bchip<z,w>\cdot v = z\cdot \acgip<w,v>.
\end{equation}
(For the sake of honesty, not to mention motivation, we should admit
that we started with \eqref{eq:1} and \eqref{eq:3}, and then used
\eqref{eq:6} to compute what \eqref{eq:4} and \eqref{eq:5} should be.)
Anyway, to check \eqref{eq:6} we compute
\begin{align*}
  \bchip<z,w>&\cdot v(x) = \int_{H} \bchip<z,w>\cdot \bigl(\eta\cdot
  v(\eta^{-1} \cdot x)\bigr) \,d\lambda_{H}^{r(x)}(\eta) \\
  &=\int_{H}\int_{G} \blip\B<z(\eta\cdot y\cdot\gamma), \eta\cdot
  w(y\cdot \gamma)> \cdot \bigl(\eta\cdot v(\eta^{-1}\cdot x)\bigr)
  \\
  &\hskip40ex d\lambda_{G}^{s(y)}(\gamma)\,d\lambda_{H}^{s(x)}(\eta) \\
  \intertext{which, after replacing $y$ by $\eta^{-1}\cdot x$ and
    taking advantage of invariance
    (Definition~\ref{def-equi-dyn-sys}(\ref{item:1})), is}
  &=\int_{H}\int_{G} \blip\B<z(x\cdot\gamma),\eta\cdot
  w(\eta^{-1}\cdot x \cdot \gamma)>\cdot \bigl(\eta\cdot
  v(\eta^{-1}\cdot x)\cdot\gamma\bigr)
  \cdot \gamma^{-1} \\
  &\hskip47ex d\lambda_{G}^{s(y)}(\gamma)\,d\lambda_{H}^{s(x)}(\eta)
  \\
  \intertext{which, since $\E_{x\cdot \gamma}$ is an \ib, is} &=
  \int_{H}\int_{G} \bigl(z(x\cdot \gamma)\cdot \brip\A<\eta\cdot
  w(\eta^{-1} \cdot x \cdot \gamma),\eta\cdot v(\eta^{-1}\cdot x)\cdot
  \gamma>\bigr) \cdot \gamma^{-1} \\
  &\hskip47ex d\lambda_{G}^{s(y)}(\gamma)\,d\lambda_{H}^{s(x)}(\eta)
  \\
  \intertext{which, in view of Lemma~\ref{lem-ip-invar}, is} &=
  \int_{H}\int_{G} \bigl(z(x\cdot \gamma)\cdot \gamma^{-1}\bigr) \cdot
  \alpha_{\gamma}\bigl(\brip\A<w(\eta^{-1} \cdot x \cdot
  \gamma), v(\eta^{-1}\cdot x)\cdot \gamma>\bigr) \\
  &\hskip47ex d\lambda_{G}^{s(y)}(\gamma)\,d\lambda_{H}^{s(x)}(\eta) \\
  &=\int_{G} \bigl(z(x\cdot \gamma)\cdot \gamma^{-1}\bigr) \cdot
  \alpha_{\gamma}\bigl( \acgip<w,v>(\gamma^{-1})\bigr)
  \,d\lambda_{G}(\gamma) \\
  &= z\cdot \acgip<w,v>(x).
\end{align*}

\begin{example}[The Scalar Case]
  \label{ex-new-scalar-thm}
  Suppose that $G$ and $H$ are second countable \lhlc{} groupoids with
  Haar systems $\set{\lambda^{u}}_{u\in \go}$ and
  $\set{\beta^{v}}_{v\in\ho}$, respectively.  Then if $X$ is a
  $(H,G)$-equivalence, we can make $\trivial{X}=X\times\C$ into a
  $(\trivial{\ho},H,\lt^{H})\sme
  (\trivial{\go},G,\lt^{G})$-equivalence in the obvious way.  Then
  Theorem~\ref{thm-renault} implies that $\cs(H,\beta)\cong
  \trivial{\ho} \rtimes_{\lt^{H}}H$ and $\cs(G,\lambda)\cong
  \trivial{\go} \rtimes_{\lt^{G}}G$ are Morita equivalent.  Therefore
  we recover the main theorem from \cite{mrw:jot87}.
\end{example}

\begin{example}[Morita Equivalence over $T$]
  \label{ex-new-ib}
  Let $p_{\A}:\A\to T$ and $p_{\B}:\B\to T$ be \usc-\cs-bundles over a
locally compact Hausdorff space $T$.  As usual, let $A=\sa_{0}(T;\A)$
and $B=\sa_{0}(T;\B)$ be the associated $C_{0}(T)$-algebras.  We can
view the topological space $T$ as a groupoid --- the so-called
\emph{co-trivial groupoid} --- and then we get dynamical systems
$(\A,T,\id)$ and $(\B,T,\id)$.  If $q:\XX\to T$ is a $(\A,T,\id)\sme
(\B,T,\id)$-equivalence, then in the case $p_{\A}$ and $p_{\B}$ are
\cs-bundles, $q$ is what we called an $\A\sme\B$-\ib{} bimodule in
\cite{kmrw:ajm98}*{Definition~2.17}.\footnote{In
  \cite{kmrw:ajm98}*{Definition~2.17}, the hypothesis that the inner
  products should be continuous on $\XX*\XX$ was inadvertently  omitted.}
As in \cite{kmrw:ajm98}*{Proposition~2.18}, $\X:=\sa_{0}(T;\XX)$ is a
$A\sme_{T} B$-\ib.  Just as in the Banach bundle case, the converse
holds: if $\X$ is a $A\sme_{T} B$-\ib, then there is an \usc-Banach
bundle $q:\XX\to T$ such that $X\cong \sa_{0}(T;\XX)$.  In the Banach
bundle case, this follows from
\cite{fd:representations1}*{Theorem~II.13.18~and Corollary~II.14.7}.
The proof in the \usc-Banach bundle case is similar (and invokes
\cite{dg:banach}*{Proposition~1.3}).
\end{example}

\begin{example}[Raeburn's Symmetric Imprimitivity Theorem]
  \label{ex-new-symmetric}
  Perhaps \emph{the} fundamental Morita equivalence result for
  ordinary crossed products is the Symmetric Imprimitivity Theorem due
  to Raeburn \cite{rae:ma88}.  We want to see here that, at least in
  the separable case, the result follows from
  Theorem~\ref{thm-renault}.  We follow the notation and treatment
  from \cite{wil:crossed}*{Theorem~4.1}.  The set-up is as follows.
  We have commuting free and proper actions of
  locally compact groups $\gK$ and $\gH$ on the left and right,
  respectively, of a locally compact space $P$ together with commuting
  actions $\galpha$ and $\gbeta$ on a \cs-algebra $D$.  In order to
  apply the equivalence theorem, we assume that $\gK$, $\gH$ and $P$
  are second countable and that $D$ is separable.

Then, as in \cite{wil:crossed}*{\S3.6}, we can form the induced
algebras $B:=\Ind_{\gH}^{P}(D,\gbeta)$ and
$A:=\Ind_{\gK}^{P}(D,\galpha)$, and the diagonal actions
\begin{equation*}
  \gsigma:\gK\to\Aut\Ind_{\gH}^{P}(D,\gbeta)\quad\text{and} \quad
  \gtau:\gH\to \Aut\Ind_{\gK}^{P}(D,\galpha)
\end{equation*}
defined in \cite{wil:crossed}*{Lemma~3.54}.  The Symmetric
Imprimitivity Theorem implies that
\begin{equation}
  \label{eq:72}
  \Ind_{\gH}^{P}(D,\gbeta)\rtimes_{\gsigma}\gK\quad\text{is Morita
    equivalent to}\quad \Ind_{\gK}^{P}(D,\galpha)\rtimes_{\gtau}\gH.
\end{equation}

Since $B=\Ind_{\gH}^{P}(D,\gbeta)$ is clearly a
$C_{0}(P/\gH)$-algebra, $B=\sa_{0}(P/\gH;\B)$ for an \usc-\cs-bundle
$\B$.  The fibre $B(p\cdot \gH)$ over $p\cdot \gH\in P/gH$ can be identified
with $\Ind_{\gH}^{p\cdot \gH}(D,\beta)$.  Of course, for any $q\in
p\cdot \gH$, the map $f\mapsto f(q)$ identifies $B(p\cdot \gH)$ with $A$.
However, this identification is not natural, and we prefer to view
$B(p\cdot \gH)$ as functions on $p\cdot \gH$.  It will be convenient to denote
elements of $\B$ as pairs $(p\cdot \gH,f)$ where $f$ is an appropriate
function on $p\cdot \gH$.  As in
Example~\vref{ex-new-group-dynamical}, we can realize
$\Ind_{\gH}^{P}(D,\gbeta)\rtimes_{\gsigma}\gK$ as a groupoid crossed
product $\B\rtimes_{\sigma}H$, where $H$ is the transformation
groupoid $H:=\gK\times P/\gH$ and $\sigma_{(t,p\cdot \gH)}$ is defined as
follows.  Given $f\in \Ind_{\gH}^{P}(\beta)$, we can view
$f\restr{t^{-1}\cdot p\cdot \gH}$ as an element of $B(t^{-1}\cdot p\cdot \gH)$,
and we get an element of $B(p\cdot \gH)$ by
\begin{equation*}
  \sigma_{(t,p\cdot \gH)}(f)(q)=\galpha_{t}\bigl(f(t^{-1}\cdot q) \bigr). 
\end{equation*}
In a similar way, we can realize
$\Ind_{\gK}^{P}(D,\galpha)\rtimes_{\gtau}\gH$ as a groupoid crossed
product $\A\rtimes_{\tau}G$ where $G$ is the transformation groupoid
$G:=\gK\backslash P\times \gH$ and $\tau_{(h,K\cdot p)}$ is given by
\begin{equation*}
  \tau_{(h,K\cdot p)}(f)(q)=\gbeta_{h}\bigl(f(q\cdot h)\bigr).
\end{equation*}

We want to derive \eqref{eq:72} from the equivalence theorem by
showing that the trivial bundle $\E=P\times A$ is a $(\B,H,\tau)\sme
(\A,G,\sigma)$ equivalence.  We have to equip $\E_{p}=\set{(p,a):a\in
  D}$ with a $B(p\cdot \gH)\sme A(\gK\cdot p)$-\ib{} structure and
specify the $H$ and $G$ actions on $\E$.  Standard computations show
that we get an \ib{} structure using the following inner-products and
actions: 
\begin{align*}
  \blip B(\gK\cdot p)<(p,a),(p,b)>(p\cdot h)&= \gbeta_{h}^{-1}(ab^{*})
  &(p\cdot \gH,f)\cdot (p,a) 
&=(p,f(p)a)
 \\
\brip A(\gK\cdot p)<(p,a),(p,b)>(t\cdot p)
&= \galpha_{t}(a^{*}b)
&(p,a)\cdot (\gK\cdot
p,f)&=(p,af(p)). 
\end{align*}
The $H$ and $G$ actions are given by
\begin{align*}
  (t,p\cdot \gH)\cdot (t^{-1}\cdot p,a)&= (p,\galpha_{t}(a))&
  (p,a)\cdot (h,\gK\cdot p)&= (p\cdot h,\gbeta_{h}^{-1}(a)).
\end{align*}
Since $P$ is a $(H,G)$-equivalence, it is now simply
a matter of checking axioms
\partref1, \partref2, \partref3 and
\partref4 of Definition~\vref{def-equi-dyn-sys}. 

Checking part~\partref1 (Continuity) at first seems awkward because
the bundles $\A$ and $\B$ are only specified indirectly.  However we
can do what we need using sections.  For example, we have the
following observation.
\begin{lemma}
  \label{lem-sections}
  The map $\E*\E\to\B$ is continuous if and only if $p\mapsto
  \blip\B<f(p),g(p)>$ is in $\sa_{c}(P;r_{P}^{*}\B)$ for all
  $f,g\in\sa_{c}(P;\E)$. 
\end{lemma}
\begin{proof}
  The $(\Longrightarrow)$ direction is immediate.  For the other
  direction, assume that $a_{i}\to a$ and $b_{i}\to b$ in $\E$ with
  $p(a_{i}) =x_{i}=p(b_{i})$ converging to $p(a)=x=p(b)$.  Then we
  need to see that $\lip\B<a_{i},b_{i}>\to \lip\B<a,b>$ in $\B$.  For
  this, it suffices to find a section $F\in\sa_{c}(P;r_{P}^{*}\B)$
  such that $F(x)=\blip\B<a,b>$ and such that
  $\|F(x_{i})-\lip\B<a_{i},b_{i}>\|\to 0$ (see
  Lemma~\vref{lem-top-sections}). Thus, we can take
  $F(p):=\blip\B<f(p),g(p)>$, where $f(x)=a$ and $g(x)=b$. 
\end{proof}
However, even with Lemma~\vref{lem-sections} in hand, there is still a
bit of work to do.  Let $P\starr P:=\set{(p,q)\in P\times
  P:p\cdot\gH=q\cdot \gH}$.  Then $P\starr P$ is locally compact and the
properness of the action implies that 
there is a continuous map $\theta:P\starr P\to \gH$ such that $q\cdot
\theta(q,p)=p$.  We know from \cite{raewil:tams85}, for example, that
$r_{P}^{*}(B)= \sa_{0}(P;r_{P}^{*}\B)$ is (isomorphic to) the balanced
tensor product
\begin{equation*}
  C_{0}(P)\tensor_{C_{0}(P/\gH)}B.
\end{equation*}
Therefore sections in $ \sa_{0}(P;r_{P}^{*}\B)$ are given by
continuous bounded $D$-valued functions on $P\starr P$ such that
$p\mapsto \|F(p,\cdot)\|$ vanishes at infinity on $P$ and such that
\begin{equation*}
  F(p,q\cdot h)=\beta_{h}^{-1}\bigl(F(p,q)\bigr).
\end{equation*}
Therefore if $f,g\in C_{c}(P,D)$, then we get a section $F$ in $
\sa_{0}(P;r_{P}^{*}\B)$ by defining 
\begin{equation*}
  F(p,q)=\blip\B<f(p),g(p)>(q) = \blip\B<f(p),g(p)>\bigl(p\cdot
  \theta(p,q)\bigr) =\beta_{\theta(q,p)}\bigl(f(p)g(p)^{*}\bigr). 
\end{equation*}
Then the continuity of the map from $\E*\E\to\B$ can be derived easily
from Lemma~\vref{lem-sections}.  The rest of part~\partref1 follows
similarly.

Part~\partref2 (Equivariance) is built in, and both part~\partref3
(Compatibility) and part~\partref4 (Invariance) follow from
straightforward computations.  Thus $\E$ is the desired equivalence.
\end{example}

We will return to the proof of the equivalence theorem in
\S\ref{sec:proof}.  In the meantime, we need to build up a bit of
technology.  In particular, we need some special approximate
identities, and we need to know that representations of crossed
products are the integrated form of covariant representations in a
manner that parallels that for ordinary dynamical systems.

\section{Approximate Identities}
\label{sec:appr-ident}

In this section, we assume throughout that $\E$ implements an
equivalence between the groupoid dynamical systems $(H,\B,\beta)$ and
$\aga$ as laid out in Definition~\vref{def-equi-dyn-sys}.  Notice that
since $H$ and $G$ are possibly non-\ttwo{} \lhlc{} groupoids, we have
to allow that our $(H,G)$-equivalence $X$ may not be Hausdorff as
well.

\begin{lemma}
  \label{lem-b-approx-id}
  Let $B=\sa_{0}(\ho;\B)$ act on $\sacc(X;\E)$ in the natural way: $
  b\cdot z(x):=b\bigl(r(x)\bigr) \cdot z(x)$.  If $\set{b_{l}}$ is an
  approximate identity for $B$, then for all $z\in\sacc(X;\E)$,
  $b_{l}\cdot z$ converges to $z$ in the inductive limit topology.
\end{lemma}
\begin{proof}
  Fix $\epsilon>0$ and a Hausdorff open set $U\subset X$.  Let $z\in
  \sa_{c}(U;\E)$.  It will suffice to see that there is an $l_{0}$
  such that $l\ge l_{0}$ implies that
  \begin{equation*}
    \|b_{l}\bigl(r(x)\bigr)\cdot z(p)-z(p)\|<\epsilon \quad\text{for
      all $p$.}
  \end{equation*}
  Let $C$ be a compact subset of $U$ such that $z$ vanishes off $U$.
  Since $\E_{x}$ is a left Hilbert $B\bigl(r(x)\bigr)$-module,
  $b_{l}\bigl(r(x)\bigr)\cdot z(x)$ converges to $z(x)$ for each
  $x$.\footnote{The relative topology on $\E_{x}$ is the Banach space
    topology \cite{dg:banach}*{p.~10}.}  Since $e\mapsto \|e\|$ is
  upper semicontinuous, there is a cover of $C$ by open sets
  $V_{1},\dots V_{n}$ such that $V_{i}\subset U$ and such that there
  is a $a_{i}\in \sa_{0}(\ho;\B)$ such that
  \begin{equation*}
    \|a_{i}\bigl(r(x)\bigr) \cdot z(x)-z(x)\|<\delta\quad\text{for all
      $x\in V_{i}$,}
  \end{equation*}
  where $\delta=\min(\epsilon/3,\epsilon/(3\|z\|_{\infty}+1))$.  Let
  $\phi_{i}\in C_{c}^{+}(U)$ be such that $\supp\phi_{i}\subset V_{i}$
  and such that $\sum\phi_{i}(x)=1$ if $x\in C$.  Define $a\in
  \sa_{c}(U;{r}^{*}\B)$ by
  \begin{equation*}
    a(x)=\sum_{i} a_{i}\bigl(r(x)\bigr)\phi_{i}(x).
  \end{equation*}
  Then for all $x\in X$,
  \begin{equation}\label{eq:54}
    \|a(x)\cdot z(x)-z(x)\|<\delta.
  \end{equation}
  We can find a $l_{0}$ such that $l\ge l_{0}$ implies
  \begin{equation*}
    \|b_{l}\bigl(r(x)\bigr) a_{i}\bigl(r(x)\bigr) -
    a_{i}\bigl(r(x)\bigr) \|<\frac\epsilon3\quad\text{for all $i$ and
      all $x$.}
  \end{equation*}
  Then
  \begin{equation}\label{eq:55}
    \|b_{l}\bigl(r(x)\bigr) \cdot a(x) -a(x)\| \le \sum_{i}
    \|b_{l}\bigl(r(x)\bigr) a_{i}\bigl(r(x)\bigr) -
    a_{i}\bigl(r(x)\bigr) \| \phi_{i}(x) 
    <\frac\epsilon3.
  \end{equation}
  If $l\ge l_{0}$, we have $\|b_{l}\bigl(r(x)\bigr)\cdot z(x)-z(x)\|$
  bounded by
  \begin{multline*}
    \|b_{l}\bigl(r(x)\bigr) \bigl(z(x)-a(x)\cdot z(x)\bigr) \| + \|
    \bigl(b_{l}\bigl(r(x)\bigr)a(x)-a(x)\bigr)\cdot z(x)\| + \|
    a(x)\cdot z(x)-z(x)\|
  \end{multline*}
  Since $\|b_{l}\bigl(r(x)\bigr)\|\le 1$ for and $x$, and in view of
  \eqref{eq:54}~and \eqref{eq:55}, the above is bounded by $\epsilon$.
  This completes the proof.
\end{proof}

\begin{lemma}
  \label{lem-hilbert-mod}
  Suppose that $U$ is a Hausdorff open subset of $X$.  Then
  $\sa_{c}(U;\E)$ becomes a left pre-Hilbert
  $\sa_{0}(U;{r}^{*}\B)$-module where the left action and pre-inner
  product are given by
  \begin{equation*}
    b\cdot z(x):=b(x)z(x) \quad\text{and}\quad
    \lip \sa_{0}(U;{r}^{*}\B)<z,w>(x):=\blip B(r(x))<z(x),w(x)>.
  \end{equation*}
\end{lemma}
\begin{proof}
  The only issues are the positivity of the inner product and the
  density of the span of the range of the inner product.  But since
  every irreducible representation of the $C_{0}(U)$-algebra
  $\sa_{0}(U;{r}^{*}\B)$ factors through a fibre, to show positivity
  it will suffice to see that for each $x$, $\blip
  B(r(x))<z(x),w(x)>\ge0$ in $B\bigl(r(x)\bigr)$.  However, this
  follows since $\E$ is an equivalence.  Furthermore, the ideal
  \begin{align*}
    I_{x}&=\operatorname{span}\set{\lip
      \sa_{0}(U;{r}^{*}\B)<z,w>(x):z,w\in\sa_{c}(U;\E\restr U)} \\
    &=\operatorname{span}\set{\blip
      B(r(x))<z(x),w(x)>:z,w\in\sa_{c}(U;\E\restr U)}
  \end{align*}
  is dense in the fibre $\sa_{0}(U;{r}^{*}\B)(x)$ over $x$.  Since the
  ideal $I$ spanned by the inner product is a $C_{0}(X)$-module, it
  follows that $I$ is dense in $\sa_{0}(U;{r}^{*}\B)$.
\end{proof}

We also need the following observation which was used in
\cite{wil:jfa81}*{p.~75} with an inadequate reference.
\begin{lemma}
  \label{lem-jfa-missing}
  Suppose that $\X$ is a full right Hilbert $A$-module.  Then sums of
  the form
  \begin{equation*}
    \sum_{i=1}^{n}\rip A<x_{i},x_{i}>
  \end{equation*}
  are dense in $A^{+}$.
\end{lemma}
\begin{remark}
  \label{rem-left-right}
  We'll actually need the left-sided version of the result.  But this
  follows immediately by taking the dual module.  (Note that a
sum is really required; think of the usual $\K(\H)$-valued inner
product on a Hilbert space $\H$.)
\end{remark}
\begin{proof}
  Fix $a\in A^{+}$.  Then $a=b^{*}b$ and since $\X$ is full, we can
  approximate $b$ by a sum
  \begin{equation*}
    \sum_{i=1}^{r}\rip A<x_{i},y_{i}>.  
  \end{equation*}
  Therefore we can approximate $a$ by
  \begin{align}
    \sum_{ij}\rip A<x_{j},y_{j}>^{*}\rip A<x_{i},y_{i}>
    &=\sum_{ij}\brip A<x_{i}\cdot{\rip A<x_{j},y_{j}>},y_{i}> \notag \\
    &=\sum_{ij} \brip A <
    \Theta_{x_{i},x_{j}}(y_{j}),y_{i}>.\label{eq:34}
  \end{align}
  But $M:=\bigl(\Theta_{x_{i},x_{j}}\bigr)$ is a positive matrix in
  $M_{r}\bigl(\K(\X)\bigr)$ (\cite{rw:morita}*{Lemma~2.65}).  Thus
  there is a matrix $(T_{ij})\in M_{r}\bigl(\K(\X)\bigr)$ such that
  \begin{equation*}
    \Theta_{x_{i},x_{j}}=\sum_{i=1}^{r} T_{ik}^{*}T_{jk}.
  \end{equation*}
  Then \eqref{eq:34} equals
  \begin{equation*}
    \sum_{ijk} \brip A<T_{jk}(y_{j}),T_{ik}(y_{i})> = \sum_{k} \rip A <
    z_{k} ,z_{k}>,
  \end{equation*}
  where
  \begin{equation*}
    z_{k}=\sum_{i} T_{ik}(y_{i}).
  \end{equation*}
  This completes the proof.
\end{proof}

\begin{cor}
  \label{cor-lemma1}
  Suppose that $b$ is a positive element in $B=\sa_{0}(\ho;\B)$, that
  $C$ is compact subset of a Hausdorff open subset $U$ of $X$.  If
  $\epsilon>0$, then there are $z_{1},\dots,z_{n}\in\sa_{c}(U;\E)$
  such that
  \begin{equation*}
    \| b\bigl(r(x)\bigr) - \sum_{i=1}^{n} \blip
    B(r(x))<z_{i}(x),z_{i}(x)> \| <\epsilon\quad\text{for all $x\in C$.}
  \end{equation*}
\end{cor}
\begin{proof}
  There is a $d\in \sa_{0}(U;{r}^{*}\B)$ such that
  $d(x)=b\bigl(r(x)\bigr)$ for all $x\in C$.  In view of
  Lemma~\vref{lem-hilbert-mod}, Lemma~\vref{lem-jfa-missing} implies
  that there are $z_{i}$ such that
  \begin{equation*}
    \| d(x)-\sum_{i} \blip B(r(x))<z_{i}(x),z_{i}(x)>\|<\epsilon
  \end{equation*}
  for all $x$.  This suffices.
\end{proof}

Since we plan to build an approximate identity, we need to recognize
one when we see one.

\begin{prop}
  \label{prop-approx-id}
  Let $\set{b_{l}}_{l\in L}$ be an approximate identity for
  $B=\sa_{0}(\ho;\B)$.  Suppose that for each $4$-tuple
  $(K,U,l,\epsilon)$ consisting of a compact subset $K\subset\ho$, a
  conditionally compact neighborhood $U$ of $\ho$ in $H$, $l\in L$ and
  $\epsilon>0$ there is a
  \begin{equation*}
    e=e_{(K,U,l,\epsilon)}\in\cchrb
  \end{equation*}
  such that
  \begin{enumerate}
  \item $\supp e\subset U$,
  \item
    \begin{math}
      \begin{displaystyle}
        \int_{H} \| e(\eta) \| \,d\lambda_{H}^{u}(\eta)\le 4
      \end{displaystyle}
    \end{math} for all $u\in K$ and
  \item
    \begin{math}
      \begin{displaystyle}
        \Bigl\|\int_{H}e(\eta)
        \,\lambda_{H}^{u}(\eta)-b_{l}(u)\Bigr\|<\epsilon\quad
        \text{for all $u\in K$.}
      \end{displaystyle}
    \end{math}
  \end{enumerate}
  Then $\set{e_{(K,U,l,\epsilon)}}$, directed by increasing $K$ and
  $l$, and decreasing $U$ and $\epsilon$, is an approximate identity
  in the inductive limit topology for both the left action of $\cchrb$
  on itself, and of $\cchrb$ on $\ccxe$.
\end{prop}
\begin{proof}
  In view of Example~\vref{ex-basic-equivalence}, it suffices to treat
  just the case of the action of $\cchrb$ on $\ccxe$.  Let $V$ be a
  \ttwo{} open subset of $X$ and let $z\in \sa_{c}(V;\E)$.  It will
  suffice to see that $e_{m}\cdot z\to z$ in the inductive limit
  topology.

  Let $K_{1}:=\supp_{V}z$.  Lemma~\vref{lem-open-nbhd} implies that
  there is a diagonally compact neighborhood $W_{1}$ of $\ho$ in $H$
  such that $K_{2}:= W_{1}\cdot K_{1}\subset V$.  Using
  Lemma~\vref{lem-t2-nbhd}, and shrinking $W_{1}$ a bit if necessary,
  we can also assume that $W_{1}r(K_{2})$ is \ttwo.

\begin{claim}\label{claim-one}
  There is a conditionally compact neighborhood $U_{1}$ of $\ho$ in
  $H$ such that $U_{1}\subset W_{1}$ and such that $\eta \in U_{1}$
  implies that
  \begin{equation}
    \label{eq:53}
    \|\eta\cdot z(\eta^{-1}\cdot x)-z(x)\|<\epsilon\quad\text{for all
      $x\in X$.}
  \end{equation}
\end{claim}
\begin{proof}[Proof of Claim]
  Notice that if the left-hand side of \eqref{eq:53} is non-zero, then
  we must have $x$ in the compact set $K_{2}$.  Therefore if the claim
  were false, then for each $U\subset W_{1}$ there would be a
  $\eta_{U}\in U$ and a $x_{U}\in K_{2}$ such that
  \begin{equation}
    \label{eq:56}
    \|\eta_{U}\cdot z(\eta_{U}^{-1}\cdot x_{U})-z(x_{U})\|\ge\epsilon.
  \end{equation}
  Since we must also have each $\eta_{U}$ in the compact set
  $W_{1}\cdot r(K_{2})$, and since each $x_{U}$ is in the compact set
  $K_{2}$, there are subnets $\set{\eta_{U_{a}}}$ and
  $\set{x_{U_{a}}}$ converging to $\eta\in W_{1}r(K_{2})$ and $x\in
  K_{2}$, respectively.  For any $U\subset W_{1}$, we eventually have
  $\eta_{U_{a}}$ in $Ur(K_{2})\subset W_{1}r(K_{2})$.  Since $
  W_{1}r(K_{2})$ is \ttwo, we must have $\eta\in Ur(K_{2})$ for all
  $U$.  Therefore $\eta\in r(K_{2})$ in view of
  Remark~\vref{rem-nbhd-base}.  Therefore $\set{\eta_{U_{a}}^{-1}\cdot
    x_{U_{a}}}$ converges to $x$ in $V$.  Thus $\eta_{U_{a}}\cdot
  z(\eta_{U_{a}}^{-1}\cdot x_{U_{a}})\to z(x)$ in $\E$.  Since
  $e\mapsto \|e\|$ is upper semicontinuous, this eventually
  contradicts \eqref{eq:56}.  This completes the proof of the claim.
\end{proof}

Lemma~\vref{lem-b-approx-id} implies that we can choose $l_{1}$ such
that $l\ge l_{1}$ implies
\begin{equation*}
  \| b_{l}\bigl(r(x)\bigr)z(x) -z(x)\|<\epsilon\quad\text{for all
    $x\in X$}.
\end{equation*}

If $e=e_{(K,U,l,\epsilon)}$ with $K\supset r(K_{2})$, $U\subset U_{1}$
and $l\ge l_{1}$, then $\|e\cdot z(x)-z(x)\|=0$ if $r(x)\notin K$ and
if $r(x)\in K$ we compute that
\begin{align*}
  \|e\cdot z(x)-z(x)\| &\le \Bigl\| \int_{H} e(\eta)\bigl(\eta\cdot z
  (\eta^{-1} \cdot x)-z(x))\bigr) \,d\lambda_{H}^{r(x)}(\eta) \Bigr\| \\
  &\qquad\qquad + \Bigl\| \Bigl( \int_{H} e(\eta)
  \,d\lambda_{H}^{r(x)}(\eta)
  -b_{l}\bigl(r(x)\bigr) \Bigr)z(x) \Bigr\| \\
  &\qquad\qquad \qquad\qquad\qquad\qquad + \bigl\|
  b_{l}\bigl(r(x)\bigr)
  z(x) -z(x)\bigr\| \\
  &\le 4\epsilon +\epsilon\|z\|_{\infty}+\epsilon.
\end{align*}
Since $\supp (e\cdot z)\subset K_{2}$, this suffices.
\end{proof}

Now we can state and prove the key result we require on approximate
identities.  It is a natural extension of
\cite{mrw:jot87}*{Proposition~2.10} to our setting.  In fact, we will
make considerable use of the constructions from \cite{mrw:jot87}.

\begin{prop}
  \label{prop-main-approx-id}
  There is a net $\set{e_{\lambda}}$ in $\cchrb$ consisting of
  elements of the form
  \begin{equation*}
    e_{\lambda} = \sum_{i=1}^{n_{\lambda}}
    \bchip<z_{i}^{\lambda},z_{i}^{\lambda}>
  \end{equation*}
  with each $z_{i}^{\lambda}\in\ccxe$, which is an approximate
  identity for the left action of $\cchrb$ on itself and on $\ccxe$.
\end{prop}

\begin{proof}
  We will apply Proposition~\vref{prop-approx-id}.  Let $\set{b_{l}}$
  be as in that proposition, and let $(K,U,l,\epsilon)$ be given.

  Let $O_{1},\dots,O_{n}$ be pre-compact \ttwo{} open sets in $X$ such
  that $\set{r(O_{i})}$ cover $K$.  Let $\set{h_{i}}\subset
  C_{c}^{+}(\ho)$ be such that $\supp h_{i}\subset r(O_{i})$ and such
  that
  \begin{equation*}
    \sum_{i=1}^{n}h_{i}(u)=1\quad\text{if $u\in K$, and}\quad
    \sum_{i=1}^{n}h_{i}(u)\le 1\quad\text{for all $u$.}
  \end{equation*}
  Let $C_{i}$ be a compact set in $O_{i}$ such that
  \begin{equation}
    \label{eq:57}
    r(C_{i})=K\cap \supp h_{i}.
  \end{equation}
  Notice that $\bigcup r(C_{i})=K$, and that there are compact
  neighborhoods $D_{i}$ of $C_{i}$ such that $D_{i}\subset O_{i}$.

  For each $i$, we will produce $e_{i}$, which is a sum of
  inner-products required in the proposition, with the additional
  properties that
  \begin{enumerate}
  \item $\supp e_{i}\subset U$,\label{page-abc}
  \item if $u\in K$, then
    \begin{equation*}
      \int_{H}\| e_{i}(\eta)\|\,d\lambda_{H}^{u}(\eta)\le 
      2\Bigl(h_{i}(u)+\frac1n\Bigr),
    \end{equation*}
    and
  \item if $u\in K$, then
    \begin{equation*}
      \Bigl\| \int_{H} e_{i}(\eta)\,d\lambda_{H}^{u}(\eta) -
      h_{i}(u)b_{l}(u)\Bigr\| <\frac\epsilon n.
    \end{equation*}
  \end{enumerate}
  Then if $e:=\sum e_{i}$, we certainly have $\supp e\subset U$.
  Furthermore, if $u\in K$, then
  \begin{align*}
    \int_{H}\|e(\eta)\|\,d\lambda_{H}^{u}(\eta) &\le \sum_{i=1}^{n}
    \int_{H} \|e_{i}(\eta)\| \,d\lambda_{H}^{u}(\eta) \\
    &\le \sum_{i=1}^{n}2\Bigl(h_{i}(u)+\frac1n\Bigr) \\
    &\le 4.
  \end{align*}
  Moreover, if $u\in K$, then
  \begin{equation*}
    \Bigl\|\int_{H}e(\eta)\,d\lambda_{H}^{u}-b_{l}(u)\Bigr\| \le
    \sum_{i=1}^{n} \Bigl\|\int_{H}e_{i}(\eta)\,d\lambda_{H}^{u}(\eta) -
    h_{i}(u) b_{l}(u)\Bigr\| <\epsilon.
  \end{equation*}
  Therefore it will suffice to produce $e_{i}$'s as described above.

  Fix $i$, and let
  $\delta=\min\bigl(\frac12,\frac2n,\frac\epsilon5\bigr)$.  Use
  Corollary~\vref{cor-lemma1} to find $z_{j}\in\sa_{c}(O_{i};\E)$ such
  that
  \begin{equation}
    \label{eq:41}
    \Bigl\| h_{i}\bigl(r(x)\bigr)b_{l}\bigl(r(x)\bigr) -\sum_{j=1}^{m} \blip
    B(r(x))<z_{j}(x),z_{j}(x)>\Bigr\| <\delta \quad\text{for all $x\in
      D_{i}$.} 
  \end{equation}
  To make some of the formulas in the sequel a little easier to
  digest, we introduce the notation
  \begin{equation*}
    \myth(\eta,y):= \sum_{j=1}^{m} \blip
    B(r(x))<z_{j}(y),\eta\cdot z_{j}(\eta^{-1}\cdot y)>.
  \end{equation*}
  Notice that the summation in \eqref{eq:41} is
  $\myth\bigl(r(x),x\bigr)$.
  \begin{claim}
    There is a conditionally compact neighborhood $W$ of $\ho$ in $H$
    such that $W\subset U$ and such that $\eta \in W$ implies that
    \begin{equation}
      \label{eq:58}
      \|\myth(\eta,y)-\myth\bigl(r(y),y\bigr)\|<\delta\quad\text{for all
        $y\in X$.}
    \end{equation}
  \end{claim}
  \begin{proof}[Proof of Claim]
    The proof follows the lines of the proof of the claim on
    page~\pageref{claim-one}.  We just sketch the details here.

    Let $K_{0}$ be a compact subset of $O_{i}$ such that for all $1\le
    j\le m$ we have $z_{j}(x)=0$ if $x\notin K_{0}$. Let $W_{1}$ be a
    diagonally compact neighborhood of $\ho$ in $H$ such that
    $W_{1}\cdot K_{0}\subset O_{i}$ and such that $W_{1}r(K_{0})$ is
    \ttwo.  If the claim where false, then for each $W\subset W_{1}$
    we could find an $x_{W}\in W_{1}\cdot K_{0}$ and an $\eta_{W} \in
    W\cap \bigl(W_{1}r(K_{0})\bigr)$ such that
    \begin{equation}
      \label{eq:59}
      \|\myth(\eta_{W},x_{W})-\myth\bigl(r(x_{W}),x_{W}\bigr)\|\ge\delta>0.
    \end{equation}
    We could then pass to subnet, relabel, and assume that $x_{W}\to
    x\in W_{1}\cdot K_{0}$ and that $\eta_{W}\to r(x)$.  Since the net
    would eventually fall in $O_{i}$, $\myth(\eta_{W},x_{W})\to
    \myth\bigl(r(x),x\bigr)$, which would eventually contradict
    \eqref{eq:59}.
  \end{proof}

  We repeat some of the constructions from
  \cite{mrw:jot87}*{Proposition~2.10} --- taking care to remain in the
  \ttwo{} realm.  Let $V_{1},\dots,V_{k}$ be pre-compact open sets
  contained in $D_{i}$ which cover $C_{i}$, and such that $(x,x\cdot
  \eta)\in V_{j}\times V_{j}$ implies that $\eta\in W$.

  Since $\set{r(V_{j})}$ covers $r(C_{i})$, there are $d_{j}\in
  C_{c}^{+}(\ho)$ such that $\supp d_{i}\subset r(V_{i})$,
  $\sum_{j}d_{j}(u)=1$ if $u\in r(C_{i})$, and $\sum_{j}d_{j}(u)\le 1
  $ for all $u$.  Since the $G$-action on $X$ is free and proper,
  Proposition~\vref{prop-mrw2.13} implies that there are $\psi_{j}\in
  C_{c}^{+}(V_{j})$ such that
  \begin{equation*}
    d_{j}\bigl(r(x)\bigr) = \int_{G} \psi_{j}(x\cdot
    \gamma)\,d\lambda_{G}^{s(x)}(\gamma) .
  \end{equation*}
  Since the $V_{j}$ are all contained in $D_{i}$, there is a constant
  $M$ such that
  \begin{equation*}
    M:= \sup_{x\in X} \sum_{j=1}^{k}\int_{G} \charfcn{V_{i}}(x\cdot
    \gamma)\, d\lambda_{G}^{s(x)}(\gamma).
  \end{equation*}
  (To see this, let $\xi_{j}\in C_{c}^{+}(O_{i})$ be such that
  $\xi_{j}(x)=1$ for all $x\in V_{j}$.
  Proposition~\vref{prop-mrw2.13} implies that $\lambda(\xi_{j})\in
  C_{c}(\xmg)$.  Then $M\le \sum_{j=1}^{k}
  \|\lambda(\xi_{j})\|_{\infty}$.)

  Using \cite{mrw:jot87}*{Lemma~2.14}, we can find $\phi_{j}\in
  C_{c}^{+}(O_{i})$ with $\supp \phi_{j}\subset V_{j}$ such that
  \begin{equation}
    \label{eq:43}
    \Bigl|\psi_{j}(x)-\phi_{j}(x)\int_{H} \phi_{j}(\eta^{-1}\cdot
    x)\,d\lambda_{H}^{r(x)} (\eta)\Bigr|<\frac\delta M.
  \end{equation}
  The point is that
  \begin{multline*}
    \Bigl|\int_{H}\int_{G} \sum_{j=1}^{k} \phi_{j}(x\cdot
    \gamma)\phi_{j}(\eta^{-1}\cdot
    x\cdot\gamma)\,d\lambda_{G}^{s(x)}(\gamma)
    \,d\lambda_{H}^{r(x)}(\eta) - \sum_{j=1}^{k}d_{j}\bigl(r(x)\bigr)
    \Bigr| \\
    = \Bigl|\int_{H}\int_{G} \sum_{j=1}^{k} \phi_{j}(x\cdot
    \gamma)\phi_{j}(\eta^{-1}\cdot
    x\cdot\gamma)\,d\lambda_{G}^{s(x)}(\gamma)
    \,d\lambda_{H}^{r(x)}(\eta) \\- \sum_{j=1}^{k}
    \int_{G}\psi_{j}(x\cdot \gamma) \,d\lambda_{G}^{s(x)}(\gamma)
    \Bigr| \\
    = \Bigl| \int_{G}\Bigl(\sum_{i=1}^{k} \phi_{j}(x\cdot
    \gamma)\int_{H} \phi_{j} (\eta^{-1}\cdot x\cdot
    \gamma)\,d\lambda_{H}^{r(x)}(\eta) -
    \psi_{j}(x\cdot \gamma)\Bigr)\,d\lambda_{G}^{s(x)}(\gamma)\Bigr| \\
    \le \frac\delta M\int_{G}\charfcn{V_{j}}(x\cdot \gamma)
    \,d\lambda_{G}^{s(x)} (\gamma) < \delta
  \end{multline*}
  To make the formulas easier to read, let
  \begin{equation*}
    F(\eta,y):=\sum_{j=1}^{k} \phi_{j}(y)\phi_{j}(\eta^{-1}\cdot y).
  \end{equation*}
  Notice that our choice of $V_{j}$'s implies that
  \begin{equation}
    \label{eq:47p}
    F(\eta,y)=0\quad\text{if $\eta\notin W$ or $y\notin D_{i}$.}
  \end{equation}
  Then the above calculation implies that if $r(x)\in r(C_{i})$, then
  \begin{equation}
    \label{eq:45}
    \Bigl|
    \int_{H}\int_{G}F(\eta,x\cdot\gamma)\,d\lambda_{G}^{s(x)}(\gamma)
    \,d\lambda_{H}^{r(x)}(\eta) -1\Bigr|<\delta,
  \end{equation}
  while $\delta<\frac12$ implies that we always have
  \begin{equation}
    \label{eq:46}
    0\le \int_{H}\int_{G}F(\eta,x\cdot\gamma)\,d\lambda_{G}^{s(x)}(\gamma)
    \,d\lambda_{H}^{r(x)}(\eta) < 1+\delta\le 2.
  \end{equation}

  Define
  \begin{equation*}
    w_{jp}(x):=\phi_{j}(x)z_{p}(x)\quad\text{and}\quad e_{i}(\eta) :=
    \sum_{jp} \bchip<w_{jp},w_{jp}>(\eta).
  \end{equation*}
  Using \eqref{eq:2}, we have
  \begin{equation*}
    e_{i}(\eta)=\int_{G}F(\eta,x\cdot \gamma)\myth(\eta,x\cdot\gamma)
    \,d\lambda_{G}^{s(x)} (\gamma).
  \end{equation*}

  If $\eta\in W$, then we chose $W$ such that
  \begin{equation}\label{eq:48}
    \|\myth(\eta,y)-\myth\bigl(r(y),y\bigr)\|<\delta\quad\text{for all $y$}.
  \end{equation}
  On the other hand, if $y\in D_{i}$, then we also have
  \begin{equation}\label{eq:51}
    \|\myth\bigl(r(y),y)-
    h_{i}\bigl(r(y)\bigr)b\bigl(r(y)\bigr)\bigr) \|<\delta.
  \end{equation}
  Since we always have $\|b_{l}(u)\|\le1$, it follows that
  \begin{equation}
    \label{eq:50}
    \|\myth(\eta,y)\|\le h_{i}\bigl(r(y)\bigr)+2\delta\le
    h_{i}\bigl(r(y)\bigr) +\frac1n\quad\text{provided
      $\eta\in W$ and $y\in D_{i}$.}
  \end{equation}

  Next we want to see that $e_{i}$ has the properties laid out on
  page~\pageref{page-abc}.  Since \eqref{eq:47p} implies that $\supp
  e_{i}\subset W$ and since we chose $W\subset U$, condition \partref1
  is clearly satisfied.

  On the other hand, if $u\in K$ and $r(x)=u$, then
  \begin{align*}
    \int_{H}\|e_{i}(\eta)\|\,d\lambda_{H}^{u} (\eta) & \le
    \int_{H}\int_{G} F(\eta,x\cdot\gamma)\|\myth(\eta,x\cdot \gamma)\|
    \,d\lambda_{G}^{s(x)}(\gamma) \,d\lambda_{H}^{u}(\eta) \\
    \intertext{which, since $F(\eta,x\cdot \gamma)=0$ unless $\eta\in
      W$ and $x\cdot\gamma\in D_{i}$ allows us to use \eqref{eq:50},
      is} &\le
    \Bigl(h_{i}(u)+\frac1n\Bigr)\int_{G}\int_{H}F(\eta,x\cdot\gamma)
    \,d\lambda_{G}^{s(x)}(\gamma) \,d\lambda_{H}^{u}(\eta) \\
    \intertext{which, by \eqref{eq:46}, is} &\le
    2\Bigl(h_{i}(u)+\frac1n\Bigr).
  \end{align*}
  Thus, \partref2 is verified.

  Similarly,
  \begin{align*}
    \Bigl\| \int_{H}&e_{i}(\eta)\,d\lambda_{H}^{u}(\eta){} -
    h_{i}(u)b_{l}(u)\Bigr\|
    \\
    &= \Bigl\| \int_{H}\int_{G} F(\eta,x\cdot\gamma)\myth(\eta,x\cdot
    \gamma) \,d\lambda_{G}^{s(x)}(\gamma) \,d\lambda_{H}^{u}(\eta)
    -h_{i}\bigl(r(x)\bigr)b_{l}\bigl(r(x)\bigr) \Bigr\| \\
    &\le \int_{H}\int_{G} F(\eta,x\cdot\gamma)
    \|\myth(\eta,x\cdot\gamma)-h_{i}\bigl(r(x)\bigr)b_{l}
    \bigl(r(x)\bigr) \| \,d\lambda_{G}^{s(x)}(\gamma)
    \,d\lambda_{H}^{u}(\eta)
    \\
    &\qquad\qquad +\int_{H}\int_{G}|F(\eta,x\cdot\gamma)-1|
    \,d\lambda_{G}^{s(x)}(\gamma) \,d\lambda_{H}^{u}(\eta)
    \|b_{l}\bigl(r(x)\bigr) \|.
  \end{align*}
  Keeping in mind that $F(\eta,x\cdot\gamma)$ vanishes off $W\times
  D_{i}$, the first of these integrals is bounded by $4\delta$ in view
  of \eqref{eq:46}, \eqref{eq:48} and \eqref{eq:51}.  Using
  \eqref{eq:45} and the fact that $\|b_{l}\|\le 1$, the second
  integral is bounded by $\delta$.  Our choice of $\delta$ implies
  that $5\delta<\epsilon$.  Therefore \partref3 is satisfied, and the
  proposition follows from Proposition~\vref{prop-approx-id}.
\end{proof}

\section{Covariant Representations}
\label{sec:covar-repr}

A critical ingredient in understanding groupoid crossed products (or
groupoid \cs-algebras for that matter) is Renault's Proposition~4.2 in
his 1987 \emph{Journal of Operator Theory} paper \cite{ren:jot87}
(cf., Theorem~\vref{thm-ren-4.2}).  To appreciate it fully, and to
make the necessary adjustments to generalize it to crossed products
(Theorem~\vref{thm-4.2-prime}), we review unitary representations of
groupoids.

Let $\mu$ be a Radon measure on $\go$.  We get a Radon measures $\nu$
and $\nu^{-1}$ on $G$ via the equations
\begin{align}
  \label{eq:10}
  \nu(f) &:=
  \int_{\go}\int_{G}f(\gamma)\,d\lambda^{u}(\gamma)\,d\mu(u)
  \quad\text{for $f\in \cc(G)$ and }\\
  \nu^{-1}(f) &:= \int_{\go}\int_{G}f(\gamma)\,d\lambda_{u}(\gamma)\,
  d\mu(u)\quad\text{for $f\in \cc(G)$.} \label{eq:18}
\end{align}

In the event $\nu$ and $\nu^{-1}$ are equivalent measures, we say that
$\mu$ is \emph{quasi-invariant}.  The modular function
$d\nu/d\nu^{-1}$ is denoted by $\Delta$.  Thus
\begin{equation}
  \label{eq:19}
  \int_{\go}\int_{G}f(\gamma)\Delta(\gamma)\,d\lambda_{u}(\gamma)\,d\mu(u)
  = 
  \int_{\go}\int_{G}f(\gamma)\,d\lambda^{u}(\gamma)\,d\mu(u).
\end{equation}

\begin{remark}
Of course $\Delta$ is only determined $\nu$-almost everywhere.
  However, $\Delta$ can always be chosen to be a homomorphism from $G$
  to the positive reals, $\R^{+}_{\times}$.  The details are spelled
  out in the proof of \cite{muh:cbms}*{Theorem~3.15}. The idea is
  this: Owing to \cite{hah:tams78}*{Corollary 3.14} and 
  \cite{ren:groupoid}*{Proposition I.3.3}, \emph{any} choice of the
  Radon-Nikodym derivative $\Delta$ is what is called an almost everywhere
  homomorphism of $G$ into $\R^{+}_{\times}$.  This means that the set
  of points $(\gamma_1,\gamma_2)\in G^{(2)}$ such that
  $\Delta(\gamma_1\gamma_2)\neq \Delta(\gamma_1)\Delta(\gamma_2)$ is a
  null set with respect to the measure 
\begin{equation*}
 \nu^{(2)}:= \int_{\go}\lambda_{u}\times
  \lambda^{u},d\mu(u).
\end{equation*} 
Since $G$ is $\sigma$-compact, \cite{ram:am71}*{Theorem~5.2} and
\cite{ram:jfa82}*{Theorem~3.2} together imply that any almost
everywhere homomorphism from $G$ to any analytic \emph{groupoid} is
equal to a homomorphism almost everywhere.
\end{remark}

As noted in \cite{muh:cbms}*{Remark~3.18}, quasi-invariant measures
are easy to come by.  Let $\mu_{0}$ be any probability measure on
$\go$ and let $\nu_{0}:=\mu_{0}\circ\lambda$ be as in \eqref{eq:10}.
Then $\nu_{0}$ is $\sigma$-finite and is equivalent to a probability
measure $\nu$ on $G$.  As show in \cite{ren:groupoid}*{pp.~24--25},
$\mu=s_{*}\nu$ (that is, $\mu(E)=\nu\bigl(s^{-1}(E)\bigr)$ is
quasi-invariant, and it is also equivalent to $\mu_{0}$ if $\mu_{0}$
was quasi-invariant to begin with.

Given a quasi-invariant measure, the next step on the way to building
unitary representations of groupoids is a \emph{Borel Hilbert Bundle
  over a space $X$}.  As explained in \cite{muh:cbms}, these are
nothing more or less than the total space of a direct integral of
Hilbert spaces \emph{a la} Dixmier.  (See also
\cite{ram:am76}*{p.~264+} and \cite{wil:crossed}*{Appendix~F}) We
start with a collection
\begin{equation*}
  \HH:=\set{\H(x)}_{x\in X}
\end{equation*}
of complex Hilbert spaces.  Then the total space is the disjoint union
\begin{equation*}
  X*\HH:=\set{(x,h):h\in\H(x)},
\end{equation*}
and we let $\pi:X*\HH\to X$ be the obvious map.
\begin{definition}
  \label{def-borel-h-bundle}
  Let $\HH=\set{\H(x)}_{x\in X}$ be a family of Hilbert spaces.  Then
  $(X*\HH,\pi)$ is an analytic (standard) Borel Hilbert Bundle if
  $X*\HH$ has an analytic (standard) Borel structure such that
  \begin{enumerate}
  \item $E$ is a Borel subset of $X$ if and only if $\pi^{-1}(E)$ is
    Borel in $X*\HH$,
  \item there is sequence $\set{f_{n}}$ of sections such that
    \begin{enumerate}
    \item the maps $\tilde f_{n}:X*\HH\to\C$ are each Borel where
      \begin{equation*}
        \tilde f_{n}(x,h):=\bip(f_{n}(x)| h),
      \end{equation*}
    \item for each $n$ and $m$,
      \begin{equation*}
        x\mapsto \bip(f_{n}(x)|f_{m}(x))
      \end{equation*}
      is Borel, and
    \item the functions $\set{\tilde f_{n}}$, together with $\pi$,
      separate points of $X*\HH$.
    \end{enumerate}
  \end{enumerate}
\end{definition}

\begin{remark}
  A section $f:X\to X*\HH$ is Borel if and only if $x\mapsto
  \bip(f(x)|f_{n}(x))$ is Borel for all $n$.  In particular, if
  $B(X*\HH)$ is the set of Borel sections and if $f\in B(X*\HH)$, then
  $x\mapsto \|f(x)\|$ is Borel.  If $\mu$ is a measure on $X$, then
  the quotient $L^{2}(X*\HH,\mu)$ of
  \begin{equation*}
    \mathcal{L}^{2}(X*\HH,\mu)=\set{f\in B(X*\HH):
      \text{$x\mapsto \|f(x)\|^{2}$
        is integrable}} ,
  \end{equation*}
  where functions agreeing $\mu$-almost everywhere are identified, is
  a Hilbert space with the obvious inner product.  Thus
  $L^{2}(X*\HH,\mu)$ is nothing more than the associated \emph{direct
    integral}
  \begin{equation*}
    \int^{\oplus}_{X}\H(x)\,d\mu(x).
  \end{equation*}
\end{remark}

\begin{definition}
  \label{def-iso-bundle}
  If $X*\HH$ is a Borel Hilbert Bundle, then its \emph{isomorphism
    groupoid} is the groupoid
  \begin{equation*}
    \operatorname{Iso}(X*\HH) := \set{(x,V,y):V\in U(\H(y),\H(x))}
  \end{equation*}
  with the weakest Borel structure such that
  \begin{equation*}
    (x,V,y)\mapsto\bip (Vf(y)|g(x))
  \end{equation*}
  is Borel for all $f,g\in B(X*\HH)$.
\end{definition}

As a Borel space, $\operatorname{Iso}(X*\HH)$ is analytic or standard
whenever $X$ has the same property.

With the preliminaries in hand, we have the machinery to make the
basic definition for the analogue of a unitary representation of a
group.  Note that we must fix a Haar system in order to make sense of 
quasi-invariant measures.
\begin{definition}
  \label{def-uni-repn}
  A unitary representation of a groupoid $G$ with Haar system
  $\set{\lambda^{u}}_{u\in\go}$ is a triple $(\mu,\go*\HH,L)$
  consisting of a quasi-invariant measure $\mu$ on $\go$, a Borel
  Hilbert bundle $\go*\HH$ over $\go$ and a Borel homomorphism $\hat
  L:G\to\operatorname{Iso}(\go*\HH)$ such that
  \begin{equation*}
    \hat L(\gamma)=\bigl(r(\gamma),L_{\gamma},s(\gamma)\bigr).
  \end{equation*}
\end{definition}

Recall that the $\|\cdot\|_{I}$-norm was defined at the end of
Section~\ref{sec:group-cross-prod}. 

\begin{prop}
  \label{prop-integrate}
  If $(\mu,\go*\HH,L)$ is a unitary representation of a \lhlc{}
  groupoid $G$, then we obtain a $\|\cdot\|_{I}$-norm bounded
  representation of
  $\cc(G)$ on
  \begin{equation*}
    \H:=\int^{\oplus}_{\go}\H(u)\,d\mu(x)=L^{2}(\go*\HH,\mu),
  \end{equation*}
  called the integrated form of $(\mu,\go*\HH,L)$, determined by
  \begin{equation}\label{eq:20}
    \bip(L(f)h|k)=\int_{G} f(\gamma)
    \bip(L_{\gamma}\bigl(h(s(\gamma))\bigr)|
    {k\bigl(r(\gamma)\bigr)})\Delta(\gamma)^{-\half}\, d\nu(\gamma).
  \end{equation}
\end{prop}
\begin{remark}
  Equation \eqref{eq:20} is convenient as it avoids dealing with
  vector-valued integration.  However, it is sometimes more convenient
  in computations to realize that \eqref{eq:20} is equivalent to
  \begin{equation}
    \label{eq:21}
    L(f)h(u) = \int_{G}f(\gamma)L_{\gamma}\bigl(h(s(\gamma)\bigr)
    \Delta(\gamma) ^{-\half}\,d\lambda^{u}(\gamma).
  \end{equation}
These sorts of vector-valued integrals are discussed in
\cite{wil:crossed}*{\S1.5}. 
  In any event, showing that $L$ is a homomorphism of
  $\cc(G)$ into $B(\H)$ is fairly straightforward and requires only
  that we recall that $\Delta$ is a homomorphism (at least almost
  everywhere).  The quasi-invariance, in the form of $\Delta$, is used
  to show that $L$ is $*$-preserving.  These assertions will follow
  from the more general results for covariant representations proved
  in 
  Proposition~\vref{prop-int-cov-repn}.
\end{remark}

We turn our attention now to the principal result in the theory:
\cite{ren:jot87}*{Proposition~4.2}.  A proof in the \ttwo{}
case is given in \cite{muh:cbms}.  This result provides very general
conditions under which a representation of a groupoid $\cs$-algebra is
the integrated form of a unitary representation of the groupoid.  In
fact, it covers representations of $\cc(G)$ acting on
pre-Hilbert spaces.  A complete proof will be given in
Appendix~\ref{sec:muhlys-proof-theorem}, but for the remainder of this
section, we will show how it may be extended to representations of
groupoid crossed products $\sacc(G;r^{*}\A)$ in the setting of
not-necessarily-Hausdorff locally compact groupoids acting on
\usc-\cs-bundles (see Theorem~\vref{thm-4.2-prime}).

\begin{thm}[Renault's Proposition 4.2]
  \label{thm-ren-4.2}
  Suppose that $\H_{0}$ is a dense subspace of a complex Hilbert space
  $\H$.  Let $L$ be a homomorphism from $\cc(G)$ into the algebra
  of linear maps on $\H_{0}$ such that
  \begin{enumerate}
  \item $\set{L(f)h:\text{$f\in \cc(G)$ and $h\in \H_{0}$}}$ is
    dense in $\H$,
  \item for each $h,k\in\H_{0}$,
    \begin{equation*}
      f\mapsto\bip(L(f)h|k)
    \end{equation*}
    is continuous in the inductive limit topology on $\cc(G)$ and
  \item for $f\in \cc(G)$ and $h,k\in\H_{0}$ we have
    \begin{equation*}
      \bip(L(f)h|k)=\bip(h|{L(f^{*})}k).
    \end{equation*}

  \end{enumerate}
  Then each $L(f)$ is bounded and extends to an operator $\bar L(f)$ on
  $\H$ of norm at most $\|f\|_{I}$.  Furthermore, $\bar L$ is a
  representation of $\cc(G)$ on $\H$ and there is a unitary
  representation $(\mu,\go*\HH,U)$ of $G$ such that $\H\cong
  L^{2}(\go*\HH,\mu)$ and $\bar L$ is (equivalent to) the integrated form
  of $(\mu,\go*\HH,U)$.
\end{thm}

Returning to the situation where we have a covariant system $\aga$,
let $(\mu,\go*\HH,U)$ be a unitary representation and let
\begin{equation*}
  \H=\int^{\oplus}_{\go}\H(u)\,d\mu(u)=L^{2}(\go*\HH,\mu)
\end{equation*}
be the associated Hilbert space.  Recall that $D\in B(\H)$ is called
\emph{diagonal} if there is a bounded Borel function $\phi\in
L^{\infty}(\mu)$ such that $D=L_{\phi}$, where by definition
\begin{equation*}
  L_{\phi}h(u)=\phi(u)h(u).
\end{equation*}
The set of diagonal operators $\mathcal{D}$ is an abelian von-Neumann
subalgebra of $B(\H)$.  The general theory of direct integrals  is based on the following basic
observations (see for example \cite{wil:crossed}*{Appendix~F}).  An
operator $T$ belongs to $\mathcal{D}'$ if and only 
if there are operators $T(u)\in B(\H(u))$ such that
\begin{equation*}
  Th(u)=T(u)\bigl(h(u)\bigr)
\end{equation*}
for $\mu$-almost every $u\in \go$ \cite{wil:crossed}*{Theorem~F.21}. 
Moreover, if $A:=\sa_{0}(\go;\A)$ and if $M:A\to B(\H)$ is a
representation such that $M(A)\subset \mathcal{D}'$, then there are
representations $M_{u}:A\to
B(\H(u))$ such that
\begin{equation}\label{eq:68}
  M(a)h(u)=M_{u}\bigl(a\bigr)\bigl(h(u)\bigr)\quad\text{for
    $\mu$-almost all $u$.}
\end{equation}
 Of course, the $M_{u}$
are only determined up to a $\mu$-null set, and it is customary to
write
\begin{equation*}
  M=\int^{\oplus}_{\go}M_{u}\,d\mu(u).
\end{equation*}
An important example for the current discussion occurs when
we are given a $C_{0}(\go)$-linear
representation $M:A\to B(L^{2}(\go*\HH,\mu))$: that is,
\begin{equation}\label{eq:74}
  M(\phi\cdot a)=L_{\phi}M(a).
\end{equation}
Then it is easy to see that $M(A)\subset \mathcal{D}'$.  \emph{In
  addition}, it is not hard to see that \eqref{eq:74} implies that
for each $u$, $\ker
M_{u}\supset I_{u}$, where $I_{u}$ is the ideal of sections in $A$
vanishing at $u$.  In particular, we can view $M_{u}$ as a
representation of the fibre $A(u)$.  Thus \eqref{eq:68} becomes
\begin{equation}
  \label{eq:69}
  M(a)h(u)=M_{u}\bigl(a(u)\bigr)\bigl(h(u)\bigr).
\end{equation}

The remainder of this section is devoted to modifying the discussion
contained in \cite{ren:jot87} to cover the setting of \usc-Banach
bundles.  Although this is straightforward, we sketch the details for
convenience.

\begin{definition}
  \label{def-cov-repn}
  A covariant representation $(M,\mu,\go*\HH,U)$ of $\aga$ consists of
  a unitary representation $(\mu,\go*\HH,U)$ and a $C_{0}(\go)$-linear
  representation $M:A\to B\bigl(L^{2}(\go*\HH,\mu)\bigr)$ decomposing
  as in \eqref{eq:69} such that
  there is a $\nu$-null set $N$ such that for all $\gamma\notin N$,
  \begin{equation}
    \label{eq:22}
    U_{\gamma}M_{s(\gamma)}(b)=M_{r(\gamma)}
    \bigl(\alpha_{\gamma}(b)\bigr)U_{\gamma}  \quad\text{for all $b\in
      A\bigl(s(\gamma)\bigr)$.} 
  \end{equation}
\end{definition}

\begin{remark}
  \label{rem-new}
  Suppose that $(M,\mu,\go*\HH,U)$ is a covariant representation of
  $\aga$ as above.  Then by definition, the set $\Sigma$ of $\gamma\in
  G$ such that \eqref{eq:22} holds in $\nu$-conull.  Since $U$ and
  $\alpha$ are bona fide homomorphisms, it is not hard to see that
  $\Sigma$ is closed under multiplication.  By a result of Ramsay's
  (\cite{ram:am71}*{Lemma~5.2} or \cite{muh:cbms}*{Lemma~4.9}), there
  is a $\mu$-conull set $V\subset \go$ such that $G\restr V\subset
  \Sigma$.
\end{remark}

\begin{prop}
  \label{prop-int-cov-repn}
  If $(M,\mu,\go*\HH,U)$ is a covariant representation of $\aga$, then
  there is a $\|\cdot\|_{I}$-norm decreasing $*$-representation $R$ of
  $\sacc(G;r^{*}\A)$ given by
  \begin{align}\label{eq:23}
    \bip(R(f)h|k)&=\int_{G}\bip(M_{r(\gamma)}
    \bigl(f(\gamma)\bigr)U_{\gamma}h\bigl(s(\gamma)\bigr)
    | {{k\bigl(r(\gamma)\bigr)}}) \Delta(\gamma)^{-\half}\,d\nu(\gamma)\\
    \intertext{or} R(f)h(u) &= \int_{G}
    M_{u}\bigl(f(\gamma)\bigr)U_{\gamma}h\bigl(s(\gamma)\bigr)
    \Delta(\gamma)^{-\half} \,d\lambda^{u}(\gamma).\label{eq:24}
  \end{align}
\end{prop}
\begin{proof}
  Clearly, \eqref{eq:23} and \eqref{eq:24} define the same operator.
  Using \eqref{eq:23}, the quasi-invariance of $\mu$ and the usual
  Cauchy-Schwartz inequality  in $L^{2}(\nu)$ we have
  \begin{align*}
    \bigl|\bip(R(f)h|k)\bigr|&\le
    \int_{G}\|f(\gamma)\|\|h\bigl(s(\gamma)\bigr) \|
    \|{k\bigl(r(\gamma)\bigr)} \|\Delta(\gamma)^{-\half}\,d\nu(\gamma)
    \\
    &\le \Bigl(\int_{G}\|f(\gamma)\|
    \|h\bigl(s(\gamma)\bigr)\|^{2}\Delta(\gamma)^{-1}
    \, d\nu(\gamma)\Bigr)^{\half} \\
    &\qquad\qquad\qquad\qquad \Bigl(\int_{G}\|f(\gamma)\| \| k
    \bigl(r(\gamma)\bigr)
    \|^{2} \,d\nu(\gamma)\Bigr)^{\half} \\
    &\le \Bigl(\|f\|_{I}\|h\|^{2}\Bigr)^{\half}
    \Bigl(\|f\|_{I}\|k\|^{2}\Bigr)^{\half}
    \\
    &=\|f\|_{I}\|h\|\|k\|.
  \end{align*}
  Therefore $R$ is bounded as claimed.
  
  To see that $R$ is multiplicative, we invoke Remark~\vref{rem-new}
  to find $\mu$-conull set $V\subset\go$ such that \eqref{eq:22} holds
  for all $\gamma\in G\restr V$.  Then if $u\in V$, we have
  \begin{align*}
    R(f&{}*g)(h)(u)= \int_{G}M_{u}\bigl(f*g(\gamma)\bigr)U_{\gamma}
    h\bigl(s(\gamma)\bigr)
    \Delta(\gamma)^{-\half}\,d\lambda^{u}(\gamma)
    \\
    &=\int_{G}\int_{G}
    M_{u}\bigl(f(\eta)\alpha_{\eta}\bigl(g(\eta^{-1}\gamma)\bigr)\bigr)
    U_{\gamma} h\bigl(s(\gamma)\bigr) \Delta(\gamma)^{-\half} \,
    d\lambda^{u}(\eta) \,d\lambda^{u}(\gamma) \\
    &=\int_{G} M_{u} \bigl(f(\eta)\bigr) \int_{G} M_{u}
    \bigl(\alpha_{\eta}\bigl(g(\eta^{-1} \gamma)\bigr)\bigr)
    U_{\gamma} h\bigl(s(\gamma)\bigr) \Delta(\gamma)^{-\half} \,
    d\lambda^{u}(\gamma)\, d\lambda^{u}(\eta)\\
    &= \int_{G} M_{u}\bigl(f(\eta)\bigr) \int_{G} M_{u}
    \bigl(\alpha_{\eta}\bigl( g(\gamma)\bigr)\bigr) U_{\eta\gamma}
    h\bigl(s(\gamma) \bigr) \Delta(\eta\gamma)^{-\half}
    \,d\lambda^{s(\eta)} (\gamma)\, d\lambda^{u}(\eta) \\
    \intertext{Now since $U_{\eta\gamma}=U_{\eta}U_{\gamma}$,
      $\Delta(\eta\gamma) =\Delta(\eta)\Delta(\gamma)$ and since
      \begin{equation*}
        U_{\eta}M_{s(\eta)}(a)=M_{u}\bigl(\alpha_{\eta}(a)\bigr)U_{\eta}
      \end{equation*}
      because $u\in V$, we have} 
    &=\int_{G}
    M_{u}\bigl(f(\eta)\bigr)U_{\eta}\Bigl( \int_{G} M_{u}
    \bigl(g(\gamma)\bigr) U_{\gamma}h\bigl(s(\gamma)\bigr)
    \Delta(\gamma)^{-\half} \, d\lambda^{s(\eta)}(\gamma) \Bigr)
    \Delta(\eta)^{-\half} \, d\lambda^{u}(\eta) \\
    &=\int_{G} M_{u}\bigl(f(\eta)\bigr)
    U_{\eta}R(g)h\bigl(s(\eta)\bigr)
    \Delta(\eta)^{-\half} \, d\lambda^{u}(\eta) \\
    &=R(f)R(g)h(u).
  \end{align*}

  We also have to see that $R$ is $*$-preserving.  This will require
  the quasi-invariance of $\mu$.
  \begin{align*}
    \bip(R(f^{*})h|k) &= \int_{G}
    \bip(M_{r(\gamma)}\bigl(f^{*}(\gamma)\bigr)U_{\gamma}h\bigl(s(\gamma)\bigr)
    | {k\bigl(r(\gamma)\bigr)}) \Delta(\gamma)^{-\half} \, d\nu(\gamma) \\
    &= \int_{G}
    \bip(M_{r(\gamma)}\alpha_{\gamma}\bigl(f(\gamma)\bigr)^{*}
    U_{\gamma} h\bigl(s(\gamma)\bigr)|
    {k\bigl(r(\gamma)\bigr)})\Delta(\gamma)^{-\half} \,d\nu(\gamma) \\
    \intertext{which, since $\Delta(\gamma)^{-\half}\,d\nu(\gamma)$ is
      invariant under inversion, is} &= \int_{G}
    \bip(M_{s(\gamma)}\bigl(\alpha_{\gamma}^{-1}\bigl(f(\gamma)\bigr)\bigr)^{*}U_{
      \gamma}^{*}h\bigl(r(\gamma)\bigr)|k\bigl(s(\gamma)\bigr))
    \Delta(\gamma)^{-\half} \, d\nu(\gamma) \\
    \intertext{Now
      \begin{equation*}
        U_{\gamma}
        M_{s(\gamma)}
        \bigl(\alpha_{\gamma}^{-1}\bigl(a\bigl(r(\gamma)\bigr)\bigr)\bigr) =
        M_{r(\gamma)} \bigl(a\bigl(r(\gamma)\bigr)\bigr) U_{\gamma}
      \end{equation*}
      for $\nu$-almost all $\gamma$.  Thus}
    &=\int_{G}\bip(h\bigl(r(\gamma)\bigr)|{
      M_{s(\gamma)}\bigl(f(\gamma)\bigr)U_{\gamma}k\bigl(s(\gamma)\bigr)})
    \Delta(\gamma)^{-\half} \,d\nu(\gamma) \\
    &=\bip(h|{R(f)k}) \qed
  \end{align*}
  \renewcommand{\qed}{}
\end{proof}

The previous result admits a strong converse in the spirit of
Renault's Theorem~\vref{thm-ren-4.2}.  The extra generality will be
used in the proof of the equivalence theorem
(Theorem~\vref{thm-renault}). 
\begin{thm}[\cite{ren:jot87}*{Lemme~4.6}]
  \label{thm-4.2-prime}
  Suppose that $\H_{0}$ is a dense subspace of a complex Hilbert space
  $\H$ and that $\pi$ is a homomorphism from $\sacc(G;r^{*}\A)$ to
  the algebra linear operators on $\H_{0}$ such that
  \begin{enumerate}
  \item $\operatorname{span}\set{\pi(f)h:\text{$f\in \sacc(G;r^{*}\A)$
        and $h\in\H_{0}$}}$ is dense in $\H$,
  \item for each $h,k\in\H_{0}$,
    \begin{equation*}
      f\mapsto \bip(\pi(f)h|k)
    \end{equation*}
    is continuous in the inductive limit topology.
  \item for each $f\in\sacc(G;r^{*}\A)$ and all $h,k\in\H_{0}$
    \begin{equation*}
      \bip(\pi(f)h|k)=\bip(h|{\pi(f^{*})}k).
    \end{equation*}
  \end{enumerate}
  Then each $\pi(f)$ is bounded and extends to a bounded operator
  $\Pi(f)$ on $\H$ such that $\Pi$ is a representation of
  $\sacc(G;r^{*}\A)$ satisfying $\|\Pi(f)\|\le\|f\|_{I}$.
  Furthermore, there is a covariant representation $(M,\mu,\go*\HH,L)$
  such that $\Pi$ is equivalent to the corresponding integrated form.
\end{thm}
\begin{proof}
  Let
  $\hoo=\operatorname{span}\set{\pi(f)h:\text{$f\in\sacc(G;r^{*}\A)$
      and $h\in\H_{0}$}}$.  The first order of business is to define
  actions of $\cc(G)$ and $A:=\sa_{0}(\go;\A)$ on $\hoo$.  If $\phi\in
  \cc(G)$, $a\in A$ and $f\in\sacc(G;r^{*}\A)$, then we define
  elements of $\sacc(G;r^{*}\A)$ as follows:
  \begin{gather}
    \label{eq:25}
    \phi\cdot f(\gamma):=\int_{G}
    \phi(\eta)\alpha_{\eta}\bigl(f(\eta^{-1}\gamma)\bigr)
    \,d\lambda^{r(\gamma)}(\eta) ,\\
    a\cdot f(\gamma):=a\bigl(r(\gamma)\bigr)f(\gamma)\quad\text{and}
    \label{eq:27}\\
    f\cdot a(\gamma):=f(\gamma)\alpha_{\gamma}\bigl(a\bigl(s(\gamma)
    \bigr)\bigr).\label{eq:28}
  \end{gather}
  Note that if $\phi_{i}\to\phi$ and $f_{i}\to f$ in the inductive
  limit topology then $\phi_{i}\cdot f_{i}\to \phi\cdot f$ in the
  inductive limit topology.
  
  Suppose that
  \begin{equation*}
    \sum_{i}\pi(f_{i})h_{i}=0
  \end{equation*}
  in $\hoo$.  As a special case of
  Proposition~\vref{prop-main-approx-id}, we know that there is an
  approximate identity $\set{e_{j}}$ in $\sacc(G;r^{*}\A)$ for the
  inductive limit topology. Thus we have
  \begin{align*}
    \sum_{i}\pi(\phi\cdot f_{i})h_{i} &= \lim_{j} \sum_{i}
    \pi\bigl(\phi\cdot (e_{j}*f_{i})\bigr)h_{i} \\
    &=\lim_{j}\pi(\phi\cdot e_{i})\Bigl(\sum_{i}\pi(f_{i})h_{i}\Bigr) \\
    &=0
  \end{align*}
  Therefore we can define a linear operator $L(\phi)$ on $\hoo$ by
  \begin{equation*}
    L(\phi)\pi(f)h:=\pi(\phi\cdot f)h.
  \end{equation*}
  It is fairly straightforward to check that $L$ satisfies (a), (b) \&
  (c) of
  Theorem~\vref{thm-ren-4.2}.  Thus Renault's Proposition 4.2
  (Theorem~\vref{thm-ren-4.2}) applies and there is a unitary
  representation $(\mu,\go*\HH,L)$ of $G$ such that $\H
  =L^{2}(\go*\HH,\mu)$ and such that the original map
  $L$ is the integrated form of $(\mu,\go*\HH,L)$.

  The action of $A=\sa_{0}(\go;\A)$ on $\sacc(G;r^{*}\A)$ given by
  \eqref{eq:27} easily extends to $\widetilde A$.  Since $\widetilde
  A$ is a unital \cs-algebra,
  \begin{equation}\label{eq:26}
    k:=(\|a\|^{2}1_{A}-a^{*}a)^{\half}
  \end{equation}
  is an element of $\tilde A$ for all $a\in A$.  Since it is easy to
  check that
  \begin{align*}
    \bip(\pi(a\cdot f)h|{\pi(g)h})&=\bip(h|{\pi\bigl((a\cdot
      f)^{*}*g\bigr)k}) \\
    &= \bip(\pi(f)h|{\pi(a^{*}\cdot g)k}),
  \end{align*}
  we can use \eqref{eq:26} to show that
  \begin{equation*}
    \Bigl\|\sum\pi(a\cdot f_{i})h_{i}\Bigr\|^{2}  = \|a\|^{2} \Bigl\|
    \sum \pi(f_{i})h_{i}\Bigr\|^{2} - \Bigl\| \sum \pi(k\cdot
    f_{i})(h_{i}) \Bigr\|^{2}.
  \end{equation*}
  It follows that
  \begin{equation*}
    M(a)\pi(f)h:=\pi(a\cdot f)h
  \end{equation*}
  defines a bounded operator on $\hoo$ which extends to a bounded
  operator $M(a)$ on $\H$ with $\|M(a)\|\le \|a\|$.  In particular,
  $M:A\to B(\H)$ is a $C_{0}(\go)$-linear 
representation of $A$ on $\H$.  Therefore $M$ decomposes as in
\eqref{eq:69}. 

  If $\phi\in \cc(G)$ and $a\in A$, then we define two
  \emph{different} elements of $\sacc(G;r^{*}\A)$ by
  \begin{equation*}
    a\tensor\phi(\gamma)=a\bigl(r(\gamma)\bigr)\phi(\gamma)\quad\text{and}
    \quad \phi\tensor a(\gamma) =
    \phi(\gamma)\alpha_{\gamma}\bigl(a\bigl(s(\gamma)\bigr)\bigr) .
  \end{equation*}
  If $g\in\sacc(G;r^{*}\A)$, then
  \begin{align*}
    (a\tensor\phi)*g(\gamma) &= \int_{G} a\bigl(r(\eta)\bigr)
    \phi(\eta) \alpha_{\eta} \bigl(g(\eta^{-1}\gamma)\bigr) \,d
    \lambda^{r(\gamma)}(\eta) \\
    &= a\bigl(r(\gamma)\bigr) \int_{G}\phi(\eta) \alpha_{\eta}
    \bigl(g(\eta^{-1} \gamma)\bigr) \,d\lambda^{r(\gamma)} (\eta) \\
    &= a\cdot \phi \cdot g(\gamma).
  \end{align*}
  Thus
  \begin{equation}\label{eq:29}
    \pi\bigl((a\tensor \phi)*g\bigr) = M(a)L(\phi)\pi(g).
  \end{equation}
  And a similar computation shows that
  \begin{equation}
    \label{eq:30}
    \pi(\phi\tensor a)=L(\phi)M(a).
  \end{equation}
  We conclude that for $h,k\in\hoo$,
  \begin{align}
    \bip(\pi(a\tensor\phi)h |k) &= \bip(M(a)L(\phi)h|k)\notag \\
    &= \int_{G} \phi(\gamma)
    \bip(M_{r(\gamma)}\bigl(a\bigl(r(\gamma)\bigr)\bigr) u_{\gamma}
    \bigl(h\bigl(s(\gamma)\bigr)\bigr) |
    {k\bigl(r(\gamma)\bigr)}) \Delta(\gamma)^{-\half}\,d\nu(\gamma)\notag \\
    &= \int_{G} \bip( M_{r(\gamma)}\bigl(a\tensor\phi(\gamma)\bigr)
    U_{\gamma} \bigl(h\bigl(s(\gamma)\bigr)\bigr) |
    {k\bigl(r(\gamma)\bigr)}) \Delta(\gamma)^{-\half}
    \,d\nu(\gamma).\label{eq:31}
  \end{align}
  Similarly,
  \begin{align}
    \bip(\pi(\phi\tensor a)h|k)&= \bip(L(\phi)M(a)H|k)\notag \\
    &= \int_{G} \phi(\gamma)
    \bip(U_{\gamma}M_{s(\gamma)}\bigl(a\bigl(s(\gamma)\bigr)\bigr)
    \bigl(h\bigl(s(\gamma)\bigr)\bigr) |
    {k\bigl(r(\gamma)\bigr)})\Delta(\gamma)^{-\half} \,d\nu(\gamma)
    \label{eq:32}
  \end{align}
  Since $\operatorname{span}\set{a\tensor\phi}$ is dense in
  $\sacc(G;r^{*}\A)$, \eqref{eq:31} must hold for all
  $f\in\sacc(G;r^{*}\A)$.  In particular, it must hold for
  $f=\phi\tensor a$, and \eqref{eq:32} must coincide with
  \begin{equation*}
    \int_{G}\phi(\gamma)\bip(M_{r(\gamma)}\bigl(\alpha_{\gamma}
    \bigl(a\bigl(s(\gamma)\bigr)\bigr)\bigr)
    U_{\gamma}\bigl(h\bigl(s(\gamma)\bigr)\bigr) |
    {k\bigl(r(\gamma)\bigr)})\Delta(\gamma)^{-\half}\,d\nu(\gamma) 
  \end{equation*}
  for all $a\in A$ and $\phi\in \cc(G)$.

  For each $a\in A$, let
\begin{equation}
    \label{eq:36}
    V(\gamma):=
    U_{\gamma}M_{s(\gamma)}\bigl(a\bigl(s(\gamma)\bigr)\bigr) -
    M_{r(\gamma)}\bigl(
    \alpha_{\gamma}\bigl(a\bigl(s(\gamma)\bigr)\bigr)\bigr)
    U_{\gamma}.
  \end{equation}
Then 
\begin{equation*}
   \int_{G}
      \phi(\gamma)\bip(V(\gamma)h\bigl(s(\gamma)|{k\bigl(r(\gamma)\bigr)})
      \Delta(\gamma)^{-\half}
      \,d \nu(\gamma)=0
\end{equation*}
for all $h,k\in L^{2}(\go*\HH,\mu)$ and $\phi\in\cc(G)$.
  In particular, 
    for each $h,k\in L^{2}(\go*\HH,\mu)$, there is a $\nu$-null set
    $N(h,k)$ such that $\gamma\notin N(h,k)$ implies that
    \begin{equation}\label{eq:35}
      \bip(V(\gamma)h\bigl(s(\gamma)\bigr)|{k\bigl(r(\gamma)\bigr)})=0.
    \end{equation}
    Since $L^{2}(\go*\HH,\mu)$ is separable, there is a $\nu$-null set
    $N$ such that $\gamma\notin N$ implies \eqref{eq:35} holds for all
    $h$ and $k$. In other works, $V(\gamma)=0$ for $\nu$-almost all $\gamma$.

  Therefore there is a $\nu$-null set $N(a)$ such that $\gamma\notin N(a)$
  implies that
  \begin{equation}
    \label{eq:37}
    U_{\gamma}M_{s(\gamma)}\bigl(a\bigl(s(\gamma)\bigr)\bigr) =
    M_{r(\gamma)}\bigl( 
    \alpha_{\gamma}\bigl(a\bigl(s(\gamma)\bigr)\bigr)\bigr)
    U_{\gamma} .
  \end{equation}
  Since $A$ is separable, and $a\mapsto a(u)$ is a surjective
  homomorphism of $A$ onto $A(u)$, there is a $\nu$-null set $N$ such
  that \eqref{eq:37} holds for all $a\in A$ and $\gamma\notin N$.

  It follows that $(M,\nu,\go*\HH,L)$ is covariant and that $\pi$ is
  the restriction of its integrated form to $\hoo$. The rest is easy.
\end{proof}

\section{Proof of the Equivalence Theorem}
\label{sec:proof}

The discussion to this point provides us with the main tools we need to
complete the proof of Theorem~\vref{thm-renault}.  Another key
observation is that the inner products and actions are continuous with
respect to the inductive limit topology.  Since this is slightly more
complicated in the not necessarily Hausdorff setting, we include a
statement and proof for convenience.

\begin{lemma}
  \label{lem-cts-ilt}
  The actions and inner products on the
  pre-$\cp(\B,H,\beta)\sme\acg$-\ib{} $\X_{0}:=\sacc(X; \E)$ of
  Theorem~\ref{thm-renault} are continuous in the inductive limit
  topology.  In particular, if $v_{i}\to v$ and $w_{i}\to w$ in the
  inductive limit topology on $\sacc(X;\E)$ and if $f_{i}\to f$ in the
  inductive limit topology on $\sacc(H;r^{*}\B)$, then
  \begin{enumerate}
  \item  $f_{i}\cdot w_{i}\to f\cdot w$ in the
  inductive limit topology on $\sacc(X;\E)$ and
\item $\acgip<w_{i},v_{i}>\to \acgip<w,v>$ in the inductive limit topology
  on $\sacc(G;r^{*}\A)$.
  \end{enumerate}
\end{lemma}
\begin{proof}
  By symmetry, it suffices to just check \partref1 and \partref2.  Let
  $K_{v}$, $K_{w}$ and $K_{f}$ be compact sets such that $v(x)=0$ if
  $x\notin K_{v}$, $w(x)=0$ if $x\notin K_{w}$ and $f(\eta)=0$ if
  $\eta\notin K_{f}$.  Then $f\cdot w(x)=0$ if $x\notin K_{f}\cdot
  K_{v}$.  Using Lemma~\vref{lem-ip-invar}, we see that $\|\eta\cdot
  w(x)\|=\|w(x)\|$, and thus $\|f\cdot w\|_{\infty}\le
  \|f\|_{\infty}\|w\|_{\infty}\sup_{u\in\ho}\lambda_{H}^{u}(K_{f})$.
  Now establishing \partref1 is straightforward.

To prove \partref2, notice that as in
Lemma~\vref{lem-fund-groupoid-lemma}, there is a continuous map
$\sigma:X\starr X\to G$ which induces a homeomorphism of $H\backslash
X\starr X$ onto $G$ such that $x\cdot \sigma(x,y)=y$.  (In particular,
$\sigma(y\cdot \gamma^{-1},y)=\gamma$.)  Thus $K_{r}\sigma(K_{w}\starr
K_{v})$ is compact and $\acgip<w,v>(\gamma)=0$ if $\gamma\notin
K_{r}$.  Also, there is a compact set $K_{1}$ such that
$s_{X}(K_{1})=s_{G}(K_{r})$.  Thus if the integral in \eqref{eq:5} is
nonzero, we can assume that $y\in K_{1}$.  Since the $H$-action is
proper,
\begin{equation*}
  K_{0}:=\set{\eta\in H:\eta\cdot K_{v}\cap K_{1}\not=0}
\end{equation*}
is compact.  Since the $G$-action on $\E$ is isometric,
$\|\acgip<w,v>\|_{\infty} \le
\|w\|_{\infty}\|v\|_{\infty}\sup_{u\in\ho}\lambda^{u}_{H}(K_{0})$, and
the rest is straightforward.
\end{proof}

We have already
observed in Remark~\vref{rem-well-defined} that
\eqref{eq:1}--\eqref{eq:5} are well-defined and take values in the
appropriate functions spaces.  To complete the proof, we are going to
apply \cite{rw:morita}*{Definition~3.9}.  We have also already checked
the required algebraic identities in parts
\partref1~and \partref4 of that Definition.  All that remains in order
to verify \partref1 is to show that inner products are positive.  This
and the density of the range of the inner products (a.k.a.\
part~\partref2) follow from Proposition~\vref{prop-main-approx-id},
Lemma~\ref{lem-cts-ilt} and symmetry by standard means (cf., e.g.,
\cite{wil:crossed}*{p.~115+}, or \cite{rie:pspm82} or the discussion
following Lemma~2 in~\cite{gre:am78}).

To establish the boundedness of the inner products, we need to verify
that
\begin{align}
  \label{eq:44}
  \acgip<f\cdot z,f\cdot z>&\le \|f\|_{\cp(\B,H,\beta)}^{2}\acgip<z,z>
  \quad\text{and}\\
  \bchip<z\cdot g,z\cdot g>&\le \|g\|_{\acg}^{2}\bchip<z,z>.
\end{align}
By symmetry, it is enough to prove \eqref{eq:44}.

But if $\rho$ is a state on $\acg$, then
\begin{equation*}
  \ip(\cdot  | \cdot)_{\rho}:=\rho\bigl(\acgip<\cdot,\cdot>\bigr)
\end{equation*}
makes $\sacc(X;\E)$ a pre-Hilbert space.  Let $\H_{0}$ be the dense
image of $\sacc(X;\E)$ in the Hilbert space completion $\H_{\rho}$.  The left
action of $\sacc(H;r_{H}^{*}\B)$ on $\sacc(X;\E)$ gives a
homomorphism $\pi$ of $\sacc(H;r_{H}^{*}\B)\subset
\cp(\B,H,\beta)$ into the linear operators on $\H_{0}$.  We want to check that the
requirements \partref1--\partref3
of Theorem~\ref{thm-4.2-prime} are satisfied.  

Notice that if $g_{i}\to g$ in the inductive limit topology on
$\sacc(G;r^{*}\A)$, then $\|g_{i}-g\|_{I}\to0$ and $g_{i}\to g$ in the
\cs-norm.  Thus, $\rho(g_{i})\to\rho(g)$.  If $f_{i}\to
f$ in the inductive limit topology on $\sacc(H;r^{*}\B)$, then
Lemma~\ref{lem-cts-ilt} implies that $\acgip<f_{i}\cdot w,v>\to
\acgip<f\cdot w,v>$ in the inductive limit topology.  Therefore
$\bip(\pi(f_{i})v|w)_{\rho}\to \bip(\pi(f)v|w)_{\rho}$.  This
establishes requirement~\partref2 of Theorem~\ref{thm-renault}.
Requirement~\partref1 follows in a similar way using the approximate
identity for $\sacc(H;r^{*}\B)$ as constructed in Proposition~\vref{prop-main-approx-id}. To see that \partref3 holds, we just need to observe
that
\begin{equation}\label{eq:75}
  \acgip<f\cdot w,v>=\acgip<w,f^{*}\cdot v>.
\end{equation}
We could verify \eqref{eq:75} directly via a complicated computation.
However, notice that \eqref{eq:75} holds for all $f$ in the span of
the left inner product as in the proof of
\cite{rw:morita}*{Lemma~3.7}.  However,
Proposition~\ref{prop-main-approx-id} implies that given any
$f\in\sacc(H;r^{*}\B)$, there is a net $\set{f_{i}}$ in the span of
the inner product such that $f_{i}\to f$ (and therefore $f_{i}^{*}\to
f^{*}$) in the inductive limit topology.  Then by
Lemma~\ref{lem-cts-ilt},
\begin{equation*}
  \acgip<f\cdot w,v>=\lim_{i}\acgip<f_{i}\cdot w,v>=\lim_{i}\acgip
  <w,f_{i}^{*}\cdot v>=\acgip<w,f^{*}\cdot v>.
\end{equation*}

Since the requirements of 
Theorem~\ref{thm-4.2-prime} are satisfied, it follows
that $\pi$ is bounded with respect to the
\cs-norm on $\sacc(H;r_{H}^{*}\B)$.  In particular,
\begin{equation*}
  \rho\bigl(\acgip<f\cdot z,f\cdot z>\bigr) \le
  \|f\|_{\cp(\B,H,\beta)}^{2} \rho\bigl(\acgip<z,z>\bigr).
\end{equation*}
As this holds for all $\rho$, \eqref{eq:44} follows, and this
completes the proof.

\section{Applications}
\label{sec:applications}

The equivalence theorem is a powerful tool, and we plan to make
considerable use of it in a subsequent paper on the equivariant Brauer
semigroup of a groupoid, extending the results in \cite{hrw:tams00} to
the groupoid setting.  Here we want to remark that a number of the
constructions and results in \cite{kmrw:ajm98} can be succinctly
described in terms of equivalences and the equivalence theorem.  

\subsection{Morita Equivalent Actions}
\label{sec:another-example}

Our first application, which asserts that Morita equivalent dynamical systems
induce Morita equivalent crossed products, is the natural
generalization to the setting of groupoids of the main results in
\cite{cmw:pams84} and~\cite{com:plms84}.
The key definition is lifted directly from
\cite{kmrw:ajm98}*{Definition~3.1}; the only difference is that we
allow the weaker notion of Banach bundle and dynamical system.
\begin{definition}
  \label{def-kmrw-3.1}
  Let $G$ be a \lhlc{} groupoid and suppose that $G$ acts on two
  \usc-\cs-bundles over $\go$, $\A$ and $\B$.  Then the two dynamical
  systems $(\A,G,\alpha)$ and $(\B,G,\alpha)$ are called \emph{Morita
    equivalent} if there is an $\A\sme\B$-\ib{} bimodule $\XX$ over
  $\go$ (see Example~\vref{ex-new-ib}), and a $G$-action on $\XX$ such that
  $x\mapsto V_{\gamma}(x):=\gamma\cdot x$ is an isomorphism and such that
  \begin{equation*}
    \blip\A<V_{\gamma}(x),V_{\gamma}(y)> =
    \alpha_{\gamma}\bigl(\lip\A<x,y>\bigr) \quad\text{and} \quad
    \brip\B<V_{\gamma}(x),V_{\gamma}(y)> =
    \beta_{\gamma}\bigl(\lip\A<x,y>\bigr).
  \end{equation*}
\end{definition}

We considered the equivalence relation of Morita equivalence of
dynamical systems in \cite{kmrw:ajm98}.  However, we did not consider
the corresponding crossed products.  But in the situation of
Definition~\vref{def-kmrw-3.1}, there is an equivalence between
$(\A,G,\alpha)$ and $(\B,G,\beta)$.  Then Theorem~\vref{thm-renault}
implies that the crossed products are Morita equivalent and provides a
concrete \ib.  This generalizes
\citelist{\cite{cmw:pams84}\cite{com:plms84}}.  It is instructive to
work out the details.  We let
\begin{equation*}
  \E:=r_{G}^{*}\XX=G*\XX=\set{(\gamma,x):s(\gamma)=p_{\XX}(x)}
\end{equation*}
with $p_{\E}$ given by $(\gamma,x)\mapsto \gamma$, and we view $G$ as
a $(G,G)$-equivalence.  Note that $\E_{x}$ is naturally identified
with $\X\bigl(r(\gamma)\bigr)$, which is given to be a
$A\bigl(r(\gamma)\bigr) \sme B\bigl(r(\gamma)\bigr)$-\ib.  However,
$\beta_{\gamma}^{-1}$ is an isomorphism of $B\bigl(r(\gamma)\bigr)$
onto $B\bigl(s(\gamma)\bigr)$, and so we obtain a
$A\bigl(r(\gamma)\bigr) \sme B\bigl(s( \gamma)\bigr)$-\ib{} via
composition.  Thus we have
\begin{align*}
  \blip\A<{(\gamma,x)},(\gamma,y)>&:=\lip\A<x,y> &
  \brip\B<{(\gamma,x)} , (\gamma,y)>&:=
  \beta_{\gamma}^{-1}\bigl(\rip\B<x,y>\bigr) \\
  a\cdot (\gamma,x)&:= (\gamma,a\cdot x) & (\gamma,x)\cdot b &:=
  \bigl(\gamma,x\cdot\beta_{\gamma}(b)\bigr) .
\end{align*}

We define commuting $G$-actions on the right and the left by
\begin{equation*}
  \sigma\cdot(\gamma,x):= \bigl(\sigma\gamma,V_{\sigma}(x)\bigr)
  \quad\text{and} \quad (\gamma,x)\cdot \sigma:= (\gamma\sigma,x).
\end{equation*}

Recall Definition~\vref{def-equi-dyn-sys}.  Clearly continuity and
equivariance are satisfied.  For compatibility, we check:
\begin{align*}
  \blip\A<\sigma\cdot{(\gamma,x)},\sigma\cdot(\gamma,y)> &=
  \blip\A<V_{\sigma}(x), V_{\sigma}(y)> \\
  &= \alpha_{\gamma}\bigl(\blip\A<{(\gamma,x)},(\gamma,y)>\bigr),
\end{align*}
while
\begin{align*}
  \brip\B<{(\gamma,x)\cdot\sigma},(\gamma,y)\cdot\sigma>&=
  \blip\B<V_{\sigma^{-1}\gamma^{-1}}
  (x),V_{\sigma^{-1}\gamma^{-1}}(y)> \\
  &= \beta_{\sigma}^{-1}
  \bigl(\brip\B<V_{\gamma^{-1}}(x),V_{\gamma^{-1}}(y)>\bigr) \\
  &= \beta_{\sigma}^{-1}\bigl(\brip\B<{(\gamma,x)},(\gamma,y)>\bigr).
\end{align*}
There are equally exciting computations involving the actions:
\begin{align*}
  \sigma\cdot \bigl(a\cdot (\gamma,x)\bigr) &=
  \sigma\cdot(\gamma,a\cdot x) \\
  &=\bigl(\sigma\gamma,V_{\sigma}(a\cdot x)\bigr) \\
  &= \bigl(\sigma\gamma,\alpha_{\sigma}(a)\cdot V_{\sigma}(x)\bigr) \\
  &= \alpha_{\sigma}(a)\cdot \bigl(\sigma\cdot(\gamma,x)\bigr),
\end{align*}
while
\begin{align*}
  \bigl((\gamma,x)\cdot b\bigr)\cdot\sigma &= \bigl(\gamma,
  x\cdot\beta_{\gamma}(b)\bigr) \cdot \sigma \\
  &= \bigl(\gamma\sigma,x\cdot\beta_{\gamma}(b)\bigr) \\
  &= (\gamma\sigma,x)\cdot \beta_{\sigma}^{-1}(b) \\
  &= \bigl((\gamma,x)\cdot \sigma\bigr)\cdot \beta_{\sigma}^{-1}(b).
\end{align*}
We also have to check invariance:
\begin{align*}
  \sigma\cdot\bigl((\gamma,x)\cdot b\bigr) &=
  \sigma\cdot\bigl(\gamma,x\cdot \beta_{\gamma}(b)\bigr) \\
  &= \bigl(\sigma\gamma,V_{\sigma}\bigl(x\cdot
  \beta_{\gamma}(b)\bigr)\bigr) \\
  &=
  \bigl(\sigma\gamma,V_{\gamma}(x)\cdot\beta_{\sigma\gamma}(b)\bigr)
  \\
  &= \bigl(\sigma\gamma,V_{\sigma}(x)\bigr)\cdot b \\
  &= \bigl(\sigma\cdot(\gamma,x)\bigr)\cdot b,
\end{align*}
while
\begin{align*}
  a\cdot\bigl((\gamma,x)\cdot\sigma\bigr) &= a\cdot(\gamma\sigma,x) \\
  &= (\gamma\sigma,a\cdot x) \\
  &= (\gamma,a\cdot x)\cdot\sigma \\
  &= \bigl(a\cdot(\gamma,x)\bigr)\cdot\sigma.
\end{align*}

Thus $\E$ is a $\aga\sme(\B,G,\beta)$-equivalence and $\acg$ is Morita
equivalent to $\cp(\B,G,\beta)$ via the completion of the pre-\ib{}
$\X_{0}=\sacc(G;\E)$.  Of course, each section of $\E$ is of the form
$z(\gamma)=(\gamma,\check z(\gamma))$ where $\check z:G\to\XX$ is a continuous
function satisfying the appropriate properties.  Naturally, we want to
identify $\X_{0}$ with these functions.  Then the appropriate inner
products and actions are given by
\begin{align}
  \tlip\acg<z,w>(\eta) &= \int_{G}
  \blip\A<z(\eta\gamma),V_{\eta}(w(\gamma))>
  \,d\lambda^{s(\eta)}(\gamma) \label{eq:14}\\
  f\cdot z(\gamma)&= \int_{G} f(\eta)\cdot V_{\eta}\bigl(
  z(\eta^{-1}\gamma)\bigr) \,d\lambda^{r(\gamma)}(\eta)\label{eq:15} \\
  \trip\cp(\B,G,\beta)<z,w>(\gamma) &= \int_{G}
  \beta_{\eta}\bigl(\rip\B< z(\eta^{-1}), w(\eta^{-1}\gamma)>\bigr)
  \,d\lambda^{r(\gamma)}(\eta) \text{ and}\label{eq:16}\\
  z\cdot g(\eta)&=
  z(\gamma)\cdot\beta_{\gamma}\bigl(g(\gamma^{-1}\eta)\bigr)
  \,d\lambda^{r(\eta)} (\gamma).\label{eq:17}
\end{align}
These equations are verified as follows: for \eqref{eq:14}, we have
\begin{align*}
  \tlip\acg<z,w>(\eta)&=\int_{G} \blip\A<z(\eta\gamma),\eta\cdot
  w(\gamma)> \,d\lambda^{s(\eta)}(\gamma) \\
  &=\int_{G} \blip\A<{(\eta\gamma,\bar z(\eta\gamma))},
  (\eta\gamma,V_{\eta}\bigl(\bar w(\gamma)\bigr))>
  \,d\lambda^{s(\eta)}(\gamma) \\
  &= \int_{G}\blip\A<\bar z(\eta\gamma),V_{\eta}\bigl(\bar
  w(\gamma)\bigr)> \,d\lambda^{s(\eta)}(\gamma) .
\end{align*}
And
\begin{align*}
  f\cdot z(\gamma)&= \int_{G} f(\eta)\cdot\bigl(\eta\cdot
  z(\eta^{-1}\gamma)\bigr) \,d\lambda^{r(\gamma)}(\eta) \\
  &= \int_{G}f(\eta)\cdot \bigl(\gamma,V_{\eta}\bigl(\bar
  z(\eta^{-1}\gamma) \bigr)\bigr)  \,d\lambda^{r(\gamma)}(\eta) \\
  &=\int_{G}\bigl(\gamma, f(\eta)\cdot V_{\eta}\bigl(\bar
  z(\eta^{-1}\gamma)\bigr)\bigr)\,d\lambda^{r(\gamma)}(\eta)
\end{align*}
gives \eqref{eq:15}.  Equations \eqref{eq:16}~and \eqref{eq:17} follow
from similar computations.

\subsection{Equivalence and the Basic Construction}
\label{sec:main-example}

In \cite{kmrw:ajm98}, we introduced the Brauer Group $\br(G)$ of a
groupoid second countable locally compact Hausdorff groupoid.  One of
the basic results is that if $X$ is a $(H,G)$-equivalence, then there
is a natural isomorphism $\phi^{X}$ of $\br(G)$ onto $\br(H)$.  The
map $\phi^{X}$ is defined via the ``basic construction'' which
associates a dynamical system $(\A^{X},H,\alpha^{X})$ to any given
dynamical system $(\A,G,\alpha)$
\cite{kmrw:ajm98}*{Proposition~2.15}. (In \cite{kmrw:ajm98}, we worked
with \cs-bundles rather than \usc-\cs-bundles, but the construction
is easily modified to handle the more general bundles we are working
with in this paper.)  We briefly recall the details.  The pull-back
\begin{equation*}
  s_{X}^{*}\A=\set{(x,a)\in X\times\A:s_{X}(x)=p(a)}
\end{equation*}
is a right $G$-space:
\begin{equation*}
  (x,a)\cdot \gamma=\bigl(x\cdot
  \gamma,\alpha_{\gamma}^{-1}(a)\bigr). 
\end{equation*}
Using the proof of \cite{kmrw:ajm98}*{Proposition~2.5}, we can show
 that the quotient
$\A^{X} := s_{X}^{*}\A/G$ is an \usc-\cs-bundle.  If we denote the
image of $(x,a)$ in $\A^{X}$ by $[x,a]$, then the action of $H$ is
given by
\begin{equation*}
  \alpha^{X}(\eta)[x,a]:=[\eta\cdot x,a].
\end{equation*}
Our goal in this section is the use the equivalence theorem to see
that $\A\rtimes_{\alpha}G$ is Morita equivalent to
$\A^{X}\rtimes_{\alpha^{X}}H$ (\emph{and} to exhibit the equivalence
bimodule).  As a special case, we see that the isomorphism $\phi^{X}:
\br(G)\to\br(H)$ induces a Morita equivalences of the corresponding
dynamical systems.

Here we let
\begin{equation}
  \label{eq:7}
  \E=s_{X}^{*}\A=\set{(x,a):s(x)=p_{\A}(a)}.
\end{equation}
Then $\E_{x}$ is easily identified with $A\bigl(s(x)\bigr)$, and we
give it an $A^{X}\bigl(r(x)\bigr) \sme A\bigl(s(x)\bigr)$-\ib{}
structure as follows.  First, since $x$ is given, it is not hard to
identify
\begin{equation}
  \label{eq:8}
  A^{X}\bigl(r(x)\bigr)=\set{[x\cdot \gamma,a]:s(\gamma)=p_{\A}(a)}
\end{equation}
with $A(s(x))$ via $[x\cdot\gamma,a]\mapsto \alpha_{\gamma}(a)$.  Then
the \ib{} structure is just the usual
$\sideset{_{A(s(x))}}{_{A(s(x))}}{\mathop{A\bigl(s(x)\bigr)}}$ one.
Thus
\begin{align*}
  \blip\A^{X}<{(x,a)},(x,b)>&:= [x,ab^{*}] & \brip\A<{(x,a)},(x,b)>&:=
  a^{*}b \\
  [x\cdot\gamma,b]\cdot(x,a)&:= \bigl(x,\alpha_{\gamma}(b)a\bigr) &
  (x,a)\cdot b &:= (x,ab).
\end{align*}
The $H$- and $G$-actions on $\E$ are given by
\begin{equation*}
  \eta\cdot (x,a):=(\eta\cdot x,a)\quad\text{and}\quad
  (x,a)\cdot\gamma := \bigl(x\cdot\gamma,
  \alpha_{\gamma}^{-1}(a)\bigr). 
\end{equation*}

It remains to check conditions (\ref{item:2})--(\ref{item:1}) of
Definition~\vref{def-equi-dyn-sys}.  We start with continuity.
Clearly the maps $\E*\E\to\A^{X}$, $\E*\E\to\A$ and $\E*\A\to\E$ are
continuous.  Showing that $\A^{X}*\E\to\E $ is continuous requires a
little fussing.  Suppose that $(x_{i},a_{i})\to (x_{0},a_{0})$ in $\E$
while $[y_{i},b_{i}]\to[y_{0},b_{0}]$ in $A^{X}$ with $x_{i}\cdot G=
y_{i}\cdot G$ for all $i$.  We need to see that
$[y_{i},b_{i}]\cdot(x_{i},a_{i}) \to [y_{0},b_{0}]\cdot(x_{0},a_{0})$.
It will suffice to see that a subnet has this property.\footnote{To
  show that a given net $\set{x_{i}}$ converges to $x$, it suffices to
  see that every subnet has a subnet converging to $x$.  In the case
  here, we can simply begin by replacing the given net with a subnet
  and then relabeling.  Then it does suffice to find a convergent
  subnet.}  Also, we may as well let $y_{i}=x_{i}$ for all $i$.  Then
we can pass to a subnet, and relabel, so that there are $\gamma_{i}$
such that $\bigl(x_{i}\cdot \gamma_{i},
\alpha_{\gamma_{i}}^{-1}(b_{i})\bigr)\to (x_{0},b_{0})$.  Since
$x_{i}\to x_{0}$ and since the $G$-action on $X$ is proper, we can
pass to another subnet, relabel, and assume that $\gamma_{i}\to
s(x_{0})$.  Thus $b_{i}=\alpha_{\gamma_{i}}\circ
\alpha_{\gamma_{i}}^{-1}(b_{i}) \to b_{0}$ and
\begin{equation*}
  [x_{i},b_{i}]\cdot (x_{i},a_{i})=(x_{i},b_{i}a_{i})\to
  (x_{0},b_{0}a_{0}) 
\end{equation*}
as required.

Equivariance is clear and invariance follows from some unexciting
computations.  For example,
\begin{align*}
  [x,b]\cdot\bigl((x,a)\cdot \gamma\bigr) &=
  [x,b]\cdot\bigl(x\cdot\gamma,\alpha_{\gamma}^{-1}(a)\bigr) \\
  &= [x\cdot\gamma,\alpha_{\gamma^{-1}}(b)]\cdot\bigl(x\cdot
  \gamma,\alpha_{\gamma}^{-1}(a)\bigr) \\
  &= \bigl(x\cdot\gamma,\alpha_{\gamma}^{-1}(ba)\bigr)\\
  &= \bigl([x,b]\cdot (s,a)\bigr) \cdot \gamma.
\end{align*}

To check compatibility, notice that
\begin{align*}
  \blip\A^{X}<{\eta\cdot(x,a)},\eta\cdot(x,b)>&= [\eta\cdot
  x,ab^{*}]\\
  &=\alpha^{X}_{\eta}\bigl(\blip\A^{X}<{(x,a)},(x,b)>\bigr).
\end{align*}
Similarly,
\begin{equation*}
  \brip\A<{(x,a)\cdot\gamma}, (x,b)\cdot\gamma>
  =\alpha_{\gamma}^{-1}\bigl( \brip\A<{(x,a)},(x,b)>\bigr).
\end{equation*}
The fact that the actions are compatible are easy, but we remark that
it also follows from invariance and Lemma~\vref{lem-ip-invar}:
\begin{align*}
  \eta\cdot\bigl(\lip\A^{X}<z,w>\cdot v\bigr) &= \eta\cdot
  \bigl(z\cdot
  \rip\A<w,v>\bigr) \\
  &=(\eta\cdot z)\cdot \bigl(\rip\A<\eta\cdot w,\eta\cdot v>\bigr) \\
  &= \lip\A^{X}<\eta\cdot z,\eta\cdot w>\cdot (\eta\cdot v) \\
  &= \alpha^{X}_{\eta}\bigl(\lip\A^{X}<z,w>\bigr) \cdot(\eta\cdot v).
\end{align*}
The fullness of the inner products gives
\begin{equation*}
  \eta\cdot (a\cdot v)= \alpha^{X}_{\eta}(a)\cdot(\eta\cdot v).
\end{equation*}

Before writing down the corresponding pre-\ib{} structure on
$\sacc(H;r_{H}^*\A^{X})$, a few comments about the nature of sections
of $r_{H}^*\A^{X}$ will be helpful.  First recall that $\A^{X}=X*\A/G$
and that we can identify $H$ with $X*X/G$ via $\eta\mapsto [\eta\cdot
x,x]$ (with any $x\in r^{-1}\bigl(s(x)\bigr)$).  Thus
\begin{equation*}
  r_{H}^{*}\A^{X}=\set{\bigl([x,y],[z,a]\bigr): \text{$x\cdot G=z\cdot
      G$ and $s(z)=p_{\A}(a)$}}.
\end{equation*}
If $X*X*\A=\set{(x,y,a):s(x)=s(y)=p_{\A}(a)}$, then
$X*X*\A/G=\A^{X*X}$ is a \cs-bundle over $H$ which is isomorphic to
$r_{H}^{*}\A^{X}$.  Consequently, $f\in\sacc(H;r_{H}^*\A^{X})$ must be
of the form
\begin{equation*}
  f\bigl([x,y]\bigr) =\bigl([x,y],[x,y,\tilde f(x,y)]\bigr)
\end{equation*}
for a function $\tilde f:X*X\to\A$ such that $p_{\A}\bigl(\tilde
f(x,y)\bigr) =s(x)$, such that $\tilde f(x\cdot \gamma,y\cdot
\gamma)=\alpha_{\gamma}^{-1} \bigl(\tilde f(x,y)\bigr)$ and such that
$\supp f/G$ is compact.  In fact, $\tilde f$ must also be continuous.
Let $(x_{i},y_{i})\to (x,y)$ in $X*X$.  Again, it will be enough to
see that $\tilde f(x_{i},y_{i})$ has a subnet converging to $\tilde
f(x,y)$.  Since $[x_{i},y_{i},\tilde f(x_{i},y_{i})]\to [x,y,\tilde
f(x,y)]$, we can pass to a subnet, relabel, and find $\gamma_{i}$ such
that $\bigl(x_{i}\cdot \gamma_{i},y_{i}\cdot
\gamma_{i},\alpha_{\gamma_{i}}^{-1} \bigl(\tilde
f(x_{i},y_{i})\bigr)\bigr)\to \bigl(x,y,\tilde f(x,y)\bigr)$.  Since
the $G$-action is proper, we can pass to another subnet, relabel, and
assume that $\gamma_{i} \to s(x)$.  It follows that $\tilde
f(x_{i},y_{i})\to f(x,y)$ as required.  Thus we will often identify
$f$ and $\tilde f$.  Moreover, we will view $\cp(\A^{X},H,\alpha^{X})$
as the completion of the set $C_{\alpha}(X*X;\A)$ of functions with
the above properties.

If $z\in\X_{0}:=\sacc(X;s_{X}^{*}\A)$, then $z(x)=(x,\check z(x))$ for
some continuous function $\check z:G\to\A$ such that $p_{\A}\bigl(\check
z(x)\bigr) =s(x)$.  Consequently,
\begin{align*}
  \tlip\cp(\A^{X},H,\alpha^{X})<z,w>(\eta) &= \int_{G}
  \blip\A^{X}<z(\eta\cdot x\cdot\gamma), \eta\cdot w(x\cdot
  \gamma)> \,d\lambda_{G}^{s(x)}(\gamma) \\
  &= \int_{G} [\eta\cdot x\cdot\gamma, \check z(\eta\cdot
  x\cdot\gamma)\check w (x\cdot \gamma)^{*}]
  \,d\lambda_{G}^{s(x)}(\gamma) \\
  &= \int_{G} [\eta\cdot x, \alpha_{\gamma}^{-1}\bigl( \check
  z(\eta\cdot x\cdot \gamma)\check w(x\cdot \gamma)^{*}\bigr)
  \,d\lambda_{G}^{s(x)}(\gamma) .
\end{align*}
Thus identifying $z$ and $\check z$, we have
\begin{equation}
  \label{eq:9}
  \tlip\cp(\A^{X},H,\alpha^{X})<z,w>(x,y) =\int_{G}
  \alpha_{\gamma}^{-1} \bigl( z(x\cdot \gamma)w(y\cdot
  \gamma)^{*}\bigr) \,d\lambda_{G}^{(s(x)}(\gamma)
\end{equation}
as a function on $X*X$.  To work out the left-action of $f\in
C_{\alpha}(X*X;\A)$, notice that
\begin{align*}
  f\cdot z(x)&=\int_{H} f(\eta)\cdot\bigl(\eta\cdot z(\eta^{-1}\cdot
  x)\bigr) \,d\lambda_{H}^{r(x)}(\eta) \\
  &= f(\eta)\cdot [x,\check z(\eta^{-1}\cdot x)]
  \,d\lambda_{H}^{r(x)}(\eta) \\
  &= \int_{H} [x,\eta^{-1}\cdot x,\tilde f(x,\eta^{-1}\cdot x)]\cdot
  [x,\check z(\eta^{-1}\cdot x)]
  \,d\lambda_{H}^{r(x)}(\eta) \\
  &= \int_{H} [x,\tilde f(x,\eta^{-1}\cdot x)\check z(\eta^{-1}\cdot x)]
  \,d\lambda_{H}^{r(x)}(\eta) .
\end{align*}
Thus, after identifications, the correct formula is
\begin{equation}
  \label{eq:11}
  f\cdot z(x)=\int_{H} f(x,\eta^{-1}\cdot x)z(\eta^{-1}\cdot x)
  \,d\lambda_{H}^{r(x)}(\eta) . 
\end{equation}
Similar, but less involved, considerations show that
\begin{gather}
  \label{eq:12}
  \acgip<z,w>(x,x\cdot \gamma) =\int_{H} z(\eta^{-1}\cdot x)^{*}
  \alpha_{\gamma} \bigl(w(\eta^{-1}\cdot x\cdot \gamma)\bigr) \,
  d\lambda_{H}^{r(x)}(\eta)\text{ and} \\
  z\cdot g(x) = \int_{G} \alpha_{\gamma}\bigl(z(x\cdot \gamma)\cdot
  g(\gamma^{-1} )\bigr) \,d\lambda_{G}^{s(x)}(\gamma).\label{eq:13}
\end{gather}

\appendix
\section{Radon Measures}
\label{sec:radon-measures}

In the proof of the disintegration theorem, we will need some facts
about \emph{complex Radon measures} and ``Radon'' measures on \lhlc{}
spaces that are a bit beyond the standard measure theory courses we
all teach --- although much of what we need in the Hausdorff case can
be found in authorities like
\cite{edw:functional}*{Chap.~4}. (In particular, complex Radon
  measures are defined in \cite{edw:functional}*{Definition~4.3.1},
  and in the Hausdorff case, the Radon-Nikodym Theorem we need can be
  sorted out from \cite{edw:functional}*{\S\S4.15.7--9}.)

\subsection{Radon Measures: The Hausdorff Case}
\label{sec:radon-meas-hausd}

To start with, let $X$ be a locally compact \emph{Hausdorff} space.
For simplicity, we will assume that $X$ is second countable.  A
(positive) \emph{Radon measure} on $X$ is a regular Borel measure
associated to a positive linear functional $\mu:C_{c}(X)\to\C$ via the
Riesz Representation Theorem. It is standard
practice amongst the cognoscenti to identify the measure and the
linear functional, and we will do so here --- cognoscente or not.
Additionally, we don't add the adjective ``positive'' unless we're
trying to be pedantic.  Notice that if $\mu$ is a Radon measure on
$X$, then $\mu:C_{c}(X)\to \C$ is continuous in the inductive limit
topology.  Thus we define a \emph{complex Radon measure} on $X$ to be
a linear functional $\nu:C_{c}(X)\to\C$ which is continuous in the
inductive limit topology.\footnote{As we shall see in the next
  paragraph, a complex Radon measure must be relatively bounded.
  Hence, if $X$ is compact, then $\nu$ is always bounded as a linear
  functional on $C(X)$, and we're back in the standard textbooks.}  If
$\nu$ is actually bounded with respect to the supremum norm on
$C_{c}(X)$, so that $\nu$ extends to a bounded linear functional on
$C_{0}(X)$, then $\nu$ is naturally associated to a bona fide complex
measure on $X$ (whose total variation norm coincides with the norm of
$\nu$ as a linear functional) \cite{rud:real}*{Theorem~6.19}.
However, in general, a complex Radon measure need not be bounded.
Nevertheless, we want to associate a measure of sorts (that is, a set
function) to $\nu$.  The problem is that for complex measures, it
doesn't make sense to talk about sets of infinite measure so we can't
expect to get a set function defined on unbounded sets in the general
case.

Let $\phi=\operatorname{Re} \nu$, the real linear functional on
$C_{c}(X)$ (viewed as a real vector space).  Fix $f\in\C_{c}^{+}(X)$
and consider
\begin{equation}
  \label{eq:60}
  \set{\phi(g)\in\R:|g|\le f}.
\end{equation}
If \eqref{eq:60} were not bounded, then we could find $g_{n}$ such
that $|g_{n}|\le \frac1n f$ and such that $|\phi(g_{n})|\ge n$.  This
gives us a contradiction since $g_{n}\to 0$ in the inductive limit
topology.  Consequently, $\phi$ is \emph{relatively bounded} as
defined in \cite{hr:abstract}*{Definition~B.31}, and
\cite{hr:abstract}*{Theorem~B.36} implies that $\phi=\mu_{1}-\mu_{2}$
where each $\mu_{i}$ is a positive linear functional on $C_{c}(X)$;
that is, each $\mu_{i}$ is a Radon measure.  Applying the same
analysis to the complex part of $\nu$, we find that there are Radon
measures $\mu_{i}$ such that $\nu=\mu_{1}-\mu_{2}+i(\mu_{3}-\mu_{4})$,
and for each $f\in C_{c}(X)$, we have
\begin{multline}\label{eq:62}
  \nu(f)=\int_{X}f(x)\,d\mu_{1}(x)-\int_{X}f(x)\,d\mu_{2}(x) \\
  +i\int_{X}f(x)\,d\mu_{3}(x) -i\int_{X}f(x)\,d\mu_{4}(x).
\end{multline}
Although in general the $\mu_{i}$ will not be finite
measures  --- so that it makes no
sense to talk about $\mu_{1}-\mu_{2}+i(\mu_{3}-\mu_{4})$ as a complex
measure on $X$ --- we nevertheless want a ``measure theory''
associated to $\nu$. (Since we are assuming that $X$ is second countable,
  Radon measures are necessarily $\sigma$-finite.) In particular, we can define
$\mu_{0}:=\mu_{1}+\mu_{2}+\mu_{3}+\mu_{4}$.  Then $\mu_{i}\ll\mu$ for
all $i$ and there are Borel functions $h_{i}:X\to [0,\infty)$ such
that $\mu_{i}=h_{i}\mu$.  Since the $h_{i}$ are finite-valued, we
can define a $\C$-valued Borel function by
$h=h_{1}-h_{2}+ih_{3}-ih_{4}$.  For each $f\in C_{c}(X)$, we have
\begin{equation*}
  \nu(f)=\int_{X}f(x)h(x)\,d\mu_{0}(x).
\end{equation*}
We can write $h(x)=\rho(x)p(x)$ for a nonnegative Borel function $p$
and a unimodular Borel function $\rho$.  Replacing $p\mu_{0}$ by
$\mu$, we then have
\begin{equation}
  \nu(f)=\int_{X}f(x)\rho(x)\,d\mu(x)\quad\text{for all $f\in
    C_{c}(X)$.} \label{eq:63}
\end{equation}
If, for example, $X$ is compact, then it is well-known that the
measure $\mu$ appearing in \eqref{eq:63} is unique, and that $\rho$ is
determined $\mu$-almost everywhere.  If $X$ is second countable, and
therefore $\sigma$-compact, then we see that $\mu$ and $\rho$ satisfy
the same uniqueness conditions.  As in the compact case, we will write
$|\nu|$ for $\mu$ and call $|\nu|$ the \emph{total variation} of
$\nu$.

Since Radon measures are finite on compact subsets, we can certainly
make perfectly good sense out of $\nu(f)$ for any $f\in \bb_{c}(X)$
--- that is, for any bounded Borel function $f$ which vanishes outside
a compact set --- simply by using \eqref{eq:63}.  (In fact, we can
make sense out of $\nu(f)$ whenever $f\in \L^{1}(|\nu|)$.)  In
particular, if $B$ is a pre-compact\footnote{We say that a set is
  \emph{pre-compact} if it is contained in a compact subset.
  Alternatively, if $X$ is Hausdorff, $B$ is pre-compact if its
  closure is compact.} Borel set in $X$, then we will happily write
$\nu(B)$ for $\nu(\charfcn B)$.  We say that a Borel set (possibly not
pre-compact) is \emph{locally $\nu$-null} if $\nu(B\cap K)=0$ for all
compact sets $K\subset X$.

We will also need a version of the Radon-Nikodym Theorem for our
complex Radon measures.  Specifically, we suppose that $\mu$ is a
Radon measure and that $\nu$ is a complex Radon measure such that
$\nu\ll \mu$ --- that is, $\mu(B)=0$ implies $B$ is locally
$\nu$-null.  If $\nu\ll\mu$ and if $\mu(E)=0$, then for each Borel set
$F\subset E$, we have
\begin{equation*}
  \int_{F} \rho(x)\,d|\nu|(x)=0.
\end{equation*}
It follows that $\rho(x)=0$ for $|\nu|$-almost all $x\in E$.  Since
$|\rho(x)|=1\not=0$ for all $x$, we must have $|\nu|(E)=0$.  That is
$\nu\ll\mu$ if and only if $|\nu|\ll\mu$.  Therefore there is a Borel
function $\phi:X\to[0,\infty)$ such that
\begin{equation*}
  \nu(f)=\int_{X}f(x)\rho(x)\,d|\nu|(x) = \int_{X}
  f(x)\rho(x)\phi(x)\,d\mu(x) ,
\end{equation*}
and we call $\frac{d\nu}{d\mu}:=\phi\rho$ the Radon-Nikodym derivative
of $\nu$ with respect to $\mu$.  Of course, $\frac{d\nu}{d\mu}$ is
determined $\mu$-almost everywhere.

\subsection{Radon Measures on Locally Hausdorff, Locally Compact
  Spaces}
\label{sec:radon-meas-locally}

Now we want to consider functionals on $\cc(X)$ where $X$ is a locally
Hausdorff, locally compact space.  The situation is more complicated
because we will not be able to invoke
\cite{hr:abstract}*{Theorem~B.36} since the vector space $\cc(X)$
\emph{need not} have the property that $f\in\cc(X)$ implies
$|f|\in\cc(X)$, and hence $\cc(X)$ need not be a lattice. This
troubling possibility was illustrated in
Example~\ref{ex-not-closed-abs}.

Consider a second countable locally Hausdorff, locally compact space
$X$.  As in the Hausdorff case, a Radon measure on $X$ starts life as
a linear functional $\mu:\cc(X)\to\C$ which is positive in the usual
sense: $f\ge 0$ should imply that $\mu(f)\ge0$.  To produced a bona
fide Borel measure on $X$ corresponding to $\mu$, we will need the
following straightforward observation.

\begin{lemma}
  \label{lem-basic-mu}
  Suppose that $(X,\M)$ is a Borel space, that $\set{U_{i}}$ is a
  cover of $X$ by Borel sets and that $\mu_{i}$ are Borel measures on
  $U_{i}$ such that if $B$ is a Borel set in $U_{i}\cap U_{j}$, then
  $\mu_{i}(B)=\mu_{j}(B)$.  Then there is a Borel measure $\mu$ on $X$
  such that $\mu\restr{U_{i}}=\mu_{i}$ for all $i$.

  Furthermore, if $\set{U_{j}'}$ and $\mu_{j}'$ is another such family
  of measures resulting in a Borel measure $\mu'$, and if the
  $\mu_{i}$ and $\mu_{j}'$ agree on overlaps as above, then
  $\mu=\mu'$.
\end{lemma}
\begin{proof}
  [Sketch of the Proof] As usual, we can find pairwise disjoint Borel
  sets $B_{i}\subset U_{i}$ such that for each $n$,
  $\bigcup_{i=1}^{n}B_{i}=\bigcup_{i=1}^{n}U_{i}$.  Then we define
  $\mu$ by
  \begin{equation*}
    \mu(B):=\sum_{i=1}^{n}\mu_{i}(B\cap B_{i}).
  \end{equation*}
  Suppose that $B$ is the countable \emph{disjoint} union
  $\bigcup_{k=1}^{\infty} E_{k}$.  Then, since the $\mu_{i}$ are each
  countably additive,
  \begin{align*}
    \mu(B) &=\sum_{i=1}^{n} \mu_{i}(B\cap B_{i})
    =\sum_{i=1}^{\infty}\sum_{k=1}^{\infty} \mu_{i}(E_{k}\cap B_{i}) =
    \sum_{k=1}^{\infty} \sum_{i=1}^{\infty} \mu_{i}(E_{k}\cap B_{i})
    \\
    &= \sum_{k=1}^{\infty}\mu(E_{k}).
  \end{align*}
  Therefore $\mu$ is a measure.

  If $B\subset U_{k}$, then
  \begin{align*}
    \mu(B)&=\sum\mu_{i}(B\cap B_{i}) \\
    \intertext{which, since $B_{i}\cap U_{k}=B_{i}\cap
      \bigcup_{j=1}^{k} B_{j}$ and since the $B_{j}$ are pairwise
      disjoint, is}
    &= \sum_{i=1}^{k} \mu_{i}(B\cap B_{i}) \\
    \intertext{which, since $B\cap B_{j}\subset U_{j}\cap U_{k}$ is}
    &= \sum_{i=1}^{k}\mu_{k}(B\cap B_{j}) 
   = \mu_{k}(B\cap \bigcup_{j=1}^{k}B_{j}) \\
    &=\mu_{k}(B).
  \end{align*}
  Thus $\mu\restr{U_{k}}=\mu_{k}$ as claimed.

  The proof of uniqueness is straightforward.
\end{proof}

If $\mu$ is a Radon measure on $\cc(X)$, we can let $\set{U_{i}}$ be a
countable open cover of $X$ by Hausdorff open sets.  We can let
$\mu_{i}:= \mu\restr{C_{c}(U_{i})}$.  Then the $\mu_{i}$ are measures
as in Lemma~\ref{lem-basic-mu}, and there is a measure $\bar\mu$ on $X$
such that $\bar\mu\restr{U_{i}}=\mu_{i}$.  If $f\in\cc(X)$, then by
\cite{khoska:jram02}*{Lemma~1.3}, we can write $f=\sum f_{i}$, where
each $f_{i}\in C_{c}(U_{i})$ and only finitely many $f_{i}$ are
nonzero.  Then
\begin{align*}
  \mu(f)&=\sum\mu_{i}(f_{i}) \\&= \sum \int_{X} f_{i}(x)\,d\mu_{i}(X)
  \\
  &=\sum \int_{X} f_{i}(x)\,d\bar\mu(x) \\
  &= \int_{X} f(x)\,d\bar\mu(x).
\end{align*}
Moreover, $\bar\mu$ does not depend on the cover $\set{U_{i}}$.  In
the sequel, we will drop the ``bar'' and write simply ``$\mu$'' for
both the linear functional and the measure as in the Hausdorff case.

Suppose that $\nu$ and $\mu$ are Radon measures on $\cc(X)$ and that
we use
the same letters for the associated measures on $X$.  As expected, we
write $\nu\ll \mu$ if $\mu(E)=0$ implies $\nu(E)=0$.  Let
$\set{U_{i}}$ be a countable cover of $X$ by Hausdorff open sets, and
let $\nu_{i}$ and $\mu_{i}$ be the associated (honest) Radon measures
on $U_{i}$.  Clearly we have $\nu_{i}\ll\mu_{i}$ and we can let
$\rho_{i}=\frac{d\nu_{i}}{d\mu_{i}}$ be the Radon-Nikodym derivative.
The usual uniqueness theorems imply that $\rho_{i}=\rho_{j}$
$\mu$-almost everywhere on $U_{i}\cap U_{j}$.  A standard argument, as
in the proof of Lemma~\ref{lem-basic-mu}, implies that there is a
Borel function $\rho:X\to[0,\infty)$ such that $\rho=\rho_{i}$
$\mu$-almost everywhere on $U_{i}$.  Then if $f=\sum f_{i}\in\cc(X)$,
we have
\begin{align*}
  \nu(f)&=\sum\nu_{i}(f_{i}) =\sum\int_{X} f_{i}(x)\rho_{i}(x)\,d\mu_{i}(x) \\
  &=\sum \int_{X}f_{i}(x)\rho(x)\,d\mu(x)=\int_{X}f(x)\rho(x)\,d\mu(x) \\
  &=\mu(f\rho).
\end{align*}
Naturally, we call $\rho$ the Radon-Nikodym derivative of $\nu$ with
respect to $\mu$.

By a \emph{complex Radon measure} on $\cc(X)$, we mean a linear
functional $\nu:\cc(X)\to\C$ which is continuous in the inductive
limit topology.  Since $\cc(X)$ is not a lattice, the usual proofs
that $\nu$ decomposes into a linear combination of (positive) Radon
measures fail (for example, the proof of 
\cite{edw:functional}*{Theorem~4.3.2} requires that
$\min(f,g)\in\cc(X)$ when $f,g\in\cc(X)$, and
Example~\ref{ex-not-closed-abs} shows this need not be the case), and
we have been unable to supply a ``non-Hausdorff'' proof.  Nevertheless,
we can employ the techniques of Lemma~\ref{lem-basic-mu} to build what
we need from an open cover $\set{U_{i}}$ of $X$ by Hausdorff subsets.
By restriction, we get complex Radon measures $\nu_{i}$ on
$C_{c}(U_{i})$.  As above there are essentially unique unimodular
functions $\rho_{i}$ such that
\begin{equation*}
  \nu_{i}(f)=\int_{X} f(x)\rho_{i}(x)\,d|\nu_{i}|(x)\quad\text{for all
    $f\in C_{c}(U_{i})$.}
\end{equation*}
Standard uniqueness arguments imply that $|\nu_{i}|(B)=|\nu_{j}|(B)$
for Borel sets $B\subset U_{i}\cap U_{j}$.  We can let $|\nu|$ be the
corresponding measure on $X$. 
 Then $\rho_{i}(x)=\rho_{j}(x)$ for
$|\nu|$-almost every $x\in B$, and we can define a Borel function
$\rho:X\to\T$ such that $\rho(x)=\rho_{i}(x)$ for $|\nu|$-almost $x\in
U_{i}$.  The measure $|\nu|$ and the $|\nu|$-equivalence class of
$\rho$ are independent of $\set{U_{i}}$, and
\begin{equation*}
  \nu(f)=\int_{X}f(x)\rho(x)\,d|\nu|(x)\quad\text{for all
    $f\in\cc(X)$.} 
\end{equation*}

Suppose that $\mu$ is a (positive) Radon measure on $\cc(X)$ and that
$\nu$ is a complex Radon measure on $\cc(X)$.  As expected, we write
$\nu\ll\mu$ if every $\mu$-null set is locally $\nu$-null.  Let
$\set{U_{i}}$ be as above.  Clearly $\nu_{i}\ll\mu_{i}$ and therefore
$|\nu_{i}|\ll \mu_{i}$.  It follows that $|\nu|\ll\mu$.  Arguing as
above, there is a $\C$-valued Borel function $\rho$ that acts as a
Radon-Nikodym derivative for $\nu$ with respect to $\mu$; that is,
\begin{equation}
  \label{eq:64}
  \nu(f)=\int_{X}f(x)\rho(x)\,d\mu(x)\quad\text{for all $f\in\cc(X)$.}
\end{equation}
Using \eqref{eq:64} and the continuity of $\nu$, it is not hard to see
that $|\nu|$ is continuous in the inductive limit topology and
therefore a Radon measure.

\section{Proof of the Disintegration Theorem}
\label{sec:muhlys-proof-theorem}

In this section, we want to give a proof of Renault's disintegration
theorem (Theorem~\vref{thm-ren-4.2}).  Let $L$, $\H$, $\H_{0}$ and
$\hoo$ be as in the statement of Theorem~\ref{thm-ren-4.2}.  In
particular, if $\operatorname{Lin}(\H_{0})$ is the collection of
linear operators on the vector space $\H_{0}$, then
$L:\cc(G)\to\operatorname{Lin}(\H_{0})$ is a homomorphism satisfying
conditions \partref1, \partref2 and \partref3 of
Theorem~\ref{thm-ren-4.2}.  If $\H_{0}'$ is a dense subspace of a
Hilbert space $\H'$, then we say that
$L':\cc(G)\to\operatorname{Lin}(\H_{0}')$ is \emph{equivalent} to $L$
is there is a unitary $U:\H\to\H'$ intertwining $L$ and $L'$ as well
as the dense subspaces $\H_{0}$ and $\H_{0}'$.

The first step in the proof will be to produce the measure $\mu$ that
appears in the direct integral in the disintegration.  This is
straightforward and is done in the next proposition.  The real work
will be to show that the measure is quasi-invariant.

\begin{prop}
  \label{prop-measure-class}
  Suppose that $L:\cc(G)\to\operatorname{Lin}(\H_{0})$ is as above.
  Then there is a representation $M:C_{0}(\go)\to B(\H)$ such that for
  all $h\in C_{0}(\go)$, $f\in\cc(G)$ and $\xi\in\H_{0}$ we have
  \begin{equation}
    \label{eq:24p}
    M(h)L(f)\xi=L\bigl((h\circ r)\cdot f\bigr)\xi.
  \end{equation}
  In particular, after replacing $L$ by an equivalent representation,
  we may assume that $\H=L^{2}(\go*\VV,\mu)$ for a Borel Hilbert
  bundle $\go*\VV$ and a finite Radon measure $\mu$ on $\go$ such that
  \begin{equation*}
    M(h)\xi(u)=h(u)\xi(u)\quad\text{for all $h\in C_{0}(\go)$ and
      $\xi\in L^{2}(\go*\VV,\mu)$.}
  \end{equation*}
\end{prop}
\begin{proof}
  We can easily make sense of $(h\circ r)\cdot f$ for $h\in
  C_{0}(\go)^{\sim}$.\footnote{As usual, if $A$ is \cs-algebra, then
    $\widetilde A$ is equal to $A$ if $1\in A$ and $A$ with a unit
    adjoined otherwise.}  Furthermore, we can compute that
  \begin{equation*}
    \bip(L((h\circ r)\cdot f)\xi|{L(g)\eta})=\bip(L(f)\xi|{L((\bar
      h\circ r)\cdot g)\eta}).
  \end{equation*}
  Then, if $k\in C_{0}(\go)^{\sim}$ is such that
  \begin{equation*}
    \|h\|_{\infty}^{2}1-|h|^{2}=|k|^{2},
  \end{equation*}
  we can compute that
  \begin{align*}
    \|h\|_{\infty}^{2}\Bigl\|\sum_{i=1}^{n} L(f_{i})\xi_{i}\Bigr\|^{2}
    -{}& \Bigl\|\sum _{i=1}^{n} L\bigl((h\circ r)\cdot
    f_{i}\bigr)\xi_{i}\Bigr\|^{2}  \\
    &=\sum_{ij} \bip(L\bigl(\bigl((\|h\|_{\infty}^{2}1-|h|^{2})\circ
    r\bigr)\cdot f_{i}\bigr)\xi_{i} | {L(f_{j}})\xi_{j}
    )\\
    &= \Bigl\|\sum_{i=1}^{n} L\bigl((k\circ r)\cdot f_{i}\bigr)
    \xi_{i}
    \Bigr\|^{2} \\
    &\ge 0
  \end{align*}
  Since $\hoo$ is dense in $\H$, it follows that there is a
  well-defined bounded operator $M(h)$ on all of $\H$ satisfying
  \eqref{eq:24p}.  It is not hard to see that $M$ is a
  $*$-homomorphism.  To see that $M$ is a representation, by
  convention, we must also see that $M$ is nondegenerate.  But if
  $f\in C_{c}(V)\subset\cc(G)$, then $r\bigl(\supp_{V}f\bigr)$ is
  compact in $\go$.  Hence there is a $h\in C_{0}(\go)$ such that
  $M(h)f=f$.  From this, it is straightforward to see that $M$ is
  nondegenerate and therefore a representation.

  Since $M$ is a representation of $C_{0}(\go)$, it is equivalent to a
  multiplication representation on $L^{2}(\go*\VV,\mu)$ for an
  appropriate Borel Hilbert bundle $\go*\VV$ and finite Radon measure
  $\mu$ --- for example, see \cite{wil:crossed}*{Example~F.25}.  The
  second assertion follows, and this completes the proof.
\end{proof}

\begin{lemma}
  \label{lem-dense}
  If $\hoo'$ is a dense subspace of $\hoo$, then
  \begin{equation*}
    \operatorname{span}\set{L(f)\xi:\text{$f\in \cc(G)$ and $\xi\in\hoo'$}}
  \end{equation*}
  is dense in $\H$.
\end{lemma}

\begin{proof}
  In view of Proposition~\vref{prop-main-approx-id}, there is a
  self-adjoint approximate identity $\set{e_{i}}$ for $\cc(G)$ in the
  inductive limit topology.  Then if $L(f)\xi\in\hoo$, we see that
  \begin{multline*}
    \|L(e_{i})L(f)\xi- L(f)\xi\|^{2} =\\
    \bip(L(f^{*}*e_{i}*e_{i}*f)\xi|\xi)
    -2\operatorname{Re}\bip(L(f^{*}*e_{i}*f)\xi|\xi) +
    \bip(L(f^{*}*f)\xi|\xi),
  \end{multline*}
  which tends to zero since $L$ is continuous in the inductive limit
  topology (by part~\partref2 of Theorem~\ref{thm-ren-4.2}).  It
  follows that $\hoo' \subset
  \overline{\operatorname{span}}\set{L(f)\xi:\text{$\xi\in\hoo'$ and
      $f\in\cc(G)$}}$.  Since $\hoo'$ is dense, the result follows.
\end{proof}

The key to Renault's proof, which we are following here, is to realize
$\H$ as the completion of (a quotient of) the algebraic tensor product
$\cc(G)\atensor \H_{0}$ which has a natural fibring over $\go$.
\begin{lemma}
  \label{lem-alg-tensor}
  Then there is a positive sesquilinear form $\rip<\cdot,\cdot>$ on
  $\ccgatho$ such that
  \begin{equation}
    \label{eq:31p}
    \rip<f\tensor \xi,g\tensor\eta>=\bip(L(g^{*}*f)\xi|\eta).
  \end{equation}
  Furthermore, the Hilbert space completion $\mathcal{K}$ of
  $\ccgatho$ is isomorphic to $\H$.  In fact, if $[f\tensor\xi]$ is
  the class of $f\tensor \xi$ in $\mathcal{K}$, then
  $[f\tensor\xi]\mapsto L(f)\xi$ is well-defined and induces an
  isomorphism of $\mathcal{K}$ with $\H$ which maps the quotient
  $\ccgatho/\N$, where $\N$ is the subspace
  $\N=\set{\sum_{i}f_{i}\tensor \xi_{i}:\sum_{i}L(f_{i})\xi_{i}=0}$ of
  vectors in $\ccgatho$ of length zero, onto $\hoo$.
\end{lemma}
\begin{proof}
  Using the universal properties of the \emph{algebraic} tensor
  product, as in the proof of \cite{rw:morita}*{Proposition~2.64} for
  example, it is not hard to see that there is a unique sesquilinear
  form on $\ccgatho$ satisfying \eqref{eq:31p}.\footnote{For fixed $g$
    and $\eta$, the left-hand side of \eqref{eq:31p} is bilinear in $f$
    and $\xi$.  Therefore, by the universal properties of the
    algebraic tensor product, \eqref{eq:31p} defines linear map
    $m(g,\eta):\ccgatho\to\C$.  Then $(g,\eta)\mapsto
    \overline{m(g,\eta)}$ is a bilinear map into the space
    $\operatorname{CL}(\ccgatho)$ of conjugate linear functionals on
    $\ccgatho$.  Then we get a linear map $N:\ccgatho\to
    \operatorname{CL}(\ccgatho)$.  We can then define
    $\rip<\alpha,\beta>:=\overline{N(\beta)(\alpha)}$.  Clearly
    $\alpha\mapsto \rip<\alpha,\beta>$ is linear and it is not hard to
    check that $\overline{\rip<\alpha,\beta>}=\rip<\beta,\alpha>$.}
  Thus to see that $\rip<\cdot,\cdot>$ is a pre-inner product, we just
  have to see that it is positive.  But
  \begin{equation}\label{eq:32p}
    \begin{split}
      \brip<\sum_{i}f_{i}\tensor\xi_{i},\sum_{i} f_{i}\tensor\xi_{i}>
      &= \sum_{ij}\bip(L(f_{j}^{*}*f_{i})\xi_{i}|\xi_{j}) \\
      &= \sum_{ij} \bip(L(f_{i})\xi_{i}|{L(f_{j})\xi_{i}}) \\
      &=\bigl\|\sum_{i}L(f_{i})\xi_{i}\bigr\|^{2}.
    \end{split}
  \end{equation}
  As in \cite{rw:morita}*{Lemma~2.16}, $\rip<\cdot,\cdot>$ defines an
  inner-product on $\ccgatho/\N$, and $[f_{i}\tensor\xi]\mapsto
  L(f_{i})\xi$ is well-defined in view of \eqref{eq:32p}.  Since this
  map has range $\hoo$ and since $\hoo$ is dense in $\H$ by
  part~\partref1 of Theorem~\ref{thm-ren-4.2}, the map extends to an
  isomorphism of $\mathcal K$ onto $\H$ as claimed.
\end{proof}

From here on, using Lemma~\ref{lem-alg-tensor}, we will normally
identify $\H$ with $\K$, and $\hoo$ with $\ccgatho/\N$.  Thus we will
interpret $[f\tensor\xi]$ as a vector in $\hoo\subset\H_{0}\subset
\H$.  Then we have
\begin{align}
  \label{eq:33}
  L(g)[f\tensor\xi]&=[g*f\tensor \xi]\quad\text{and}\\
  M(h)[f\tensor\xi]&=[(h\circ r)\cdot f\tensor\xi],\label{eq:34p}
\end{align}
where $M$ is the representation of $C_{0}(\go)$ defined in
Proposition~\vref{prop-measure-class}, $g\in\cc(G)$ and $h\in
C_{0}(\go)$.

\begin{remark}
  \label{rem-extend-m}
  In view of Proposition~\vref{prop-measure-class}, $M$ extends to a
  $*$-homomorphism of $\bb_{c}(G)$ into $B(\H)$ such that $M(h)=0$ if
  $h(u)=0$ for $\mu$-almost all $u$ (where $\mu$ is the measure
  defined in that proposition).  However, at this point, we can not
  assert that \eqref{eq:34p} holds for any $h\notin C_{0}(\go)$.
\end{remark}

Showing that $\mu$ is quasi-invariant requires that we extend
equations \eqref{eq:33} and \eqref{eq:34p} to a larger class of
functions.  This can't be done without also enlarging the domain of
definition of $L$.  This is problematic as we don't as yet know that
each $L(f)$ is bounded in any sense, and $\H_{0}$ is not complete.
We'll introduce only those functions we absolutely need.

\begin{definition}
  \label{def-bboc}
  Suppose that $V$ is an open Hausdorff set in $G$.  Let $\bboc(V)$ be
  the collection of bounded Borel functions on $V$ which are the
  pointwise limit of a uniformly bounded sequence $\set{f_{n}}\subset
  C_{c}(V)$ such that there is a compact set $K\subset V$ such that
  $\supp f_{n}\subset K$ for large $n$.  We let $\bbocg$ be the vector
  space spanned by the $\bboc(V)$ for all $V\subset G$ open and
  Hausdorff.
\end{definition}

It is important to note that $\bbocg$ is not a very robust class of
functions on $G$.  In particular, it is not closed under the type of
convergence used in its definition.  Nevertheless, its elements are
all integrable with respect to any Radon measure on $G$, and the
following lemma is an easy consequence of the dominated convergence
theorem applied to the total variation measure.
\begin{lemma}
  \label{lem-grm-ext}
  Suppose that $\sigma$ is a complex Radon measure on $\cc(G)$ such
  that
  \begin{equation}
    \label{eq:65}
    \sigma(f)=\int_{G} f(\gamma)\rho(\gamma)\,d|\sigma|(\gamma)
  \end{equation}
  for a unimodular function $\rho$ and total variation $|\sigma|$ (see
  Appendix~\ref{sec:radon-meas-locally}).  Then $\sigma $ extends to a
  linear functional on $\bbocg$ such that \eqref{eq:65} holds and such
  that if $\set{f_{n}}$ is a uniformly bounded sequence in $\bbocg$
  converging pointwise to $f\in \bbocg$ with supports eventually
  contained in a fixed compact set, then $\sigma(f_{n})\to\sigma(f)$.
\end{lemma}
\begin{proof}[Sketch of the Proof]
  Since $|\sigma|$ is a Radon measure, \eqref{eq:65} makes good sense
  for any $f\in\bbocg$.  Thus $\sigma$ extends as claimed.  The rest
  is an easy consequence of the dominated convergence theorem applied
  to $|\sigma|$.
\end{proof}

\begin{cor}
  \label{cor-bbocg-conv}
  If $f,g\in\bbocg$, then
  \begin{equation*}
    f*g(\gamma):=\int_{G}
    f(\eta)g(\eta^{-1}\gamma)\,d\lambda^{r(\gamma)}( \eta)
  \end{equation*}
  defines an element $f*g$ of $\bbocg$.
\end{cor}
\begin{proof}
  As in the proof of Proposition~\vref{prop-star-alg}, we can assume
  that there are Hausdorff open sets $U$ and $V$ such that $UV$ is
  Hausdorff and such that $f\in \bboc(U)$ while $g\in \bboc(V)$.  Let
  $\set{f_{n}}$ and $\set{g_{n}}$ be uniformly bounded sequences in
  $C_{c}(U)$ and $C_{c}(V)$, respectively, such that $f_{n}\to f$ and
  $g_{n}\to g$ pointwise and with supports contained in a fixed
  compact set.  Since
  \begin{equation*}
    |f_{n}*g_{n}(\gamma)|\le
    \|f_{n}\|_{\infty}\|g_{n}\|_{\infty}\sup_{u\in\go}
    \lambda^{u}\bigl((\supp f_{n})(\supp g_{n})\bigr),
  \end{equation*}
  it follows that $\set{f_{n}*g_{n}}$ is a uniformly bounded sequence
  in $C_{c}(UV)$, all of whose supports are in a fixed compact set,
  converging pointwise to $f*g$.  Thus $f*g\in
  \bboc(UV)\subset\bbocg$.
\end{proof}

In view of the continuity assumption on $L$, we can define a complex Radon
measure $L_{\xi,\eta}$ on $\cc(G)$ via
\begin{equation*}
  L_{\xi,\eta}(f):=\bip(L(f)\xi|\eta)
\end{equation*}
for each $\xi$ and $\eta$ in $\H_{0}$.
Keep in mind that we can extend $L_{\xi,\eta}$ to a linear functional
on all of $\bbocg$.

\begin{lemma}
  \label{lem-paul2}
  There is a positive sesquilinear form on $\bbocgatho$, extending
  that on $\ccgatho$, such that
  \begin{equation*}
    \rip<f\tensor \xi,g\tensor
    \eta>=L_{\xi,\eta}(g^{*}*f)\quad\text{for all $f,g\in\cc(G)$ and
      $\xi,\eta\in \H_{0}$.}
  \end{equation*}
  In particular, if
  \begin{equation*}
    \N_{b}:=\set{\sum_{i}f_{i}\tensor\xi\in\ccgatho:\brip<\sum_{i}f_{i}\tensor \xi, \sum_{i}f_{i}\tensor\xi_{i}>=0}
  \end{equation*}
  is the subspace of vectors of zero length, then the quotient
  $\bbocgatho/\N_{b}$ can be identified with a subspace of $\H$
  containing $\hoo:=\ccgatho/\N$.
\end{lemma}
\begin{remark}
  \label{rem-notation}
  As before, we will write $[f\tensor\xi]$ for the class of
  $f\tensor\xi$ in the quotient $\bbocgatho/\N_{b}\subset \H$.
\end{remark}

\begin{proof}
  Just as in Lemma~\vref{lem-alg-tensor}, there is a well-defined
  sesquilinear form on $\bbocgatho$ extending that on $\gcbatho$.
  (Note that the right-hand side of \eqref{eq:31p} can be rewritten as
  $L_{\xi,\eta}(g^{*}*f)$.) In particular, we have
  \begin{equation*}
    \brip<\sum_{i} f_{i}\tensor\xi_{i},\sum_{j} g_{j}\tensor\eta_{j}> =
    \sum_{ij} L_{\xi_{i},\eta_{j}}(g_{j}^{*}*f_{i}).
  \end{equation*}
  We need to see that the form is positive.  Let $\alpha:=\sum_{i}
  f_{i}\tensor \xi_{i}$, and let $\set{f_{i,n}}$ be a uniformly
  bounded sequence in $\cc(G)$ converging pointwise to $f_{i}$ with
  all the supports contained in a fixed compact set.  Then for each
  $i$ and $j$, $f_{j,n}^{*}*f_{i,n}\to f_{j}^{*}*f_{i}$ in the
  appropriate sense.  In particular, Lemma~\vref{lem-grm-ext} implies
  that
  \begin{align*}
    \rip<\alpha,\alpha>&= \sum_{ij}
    L_{\xi_{i},\xi_{j}}(f_{j}^{*}*f_{i}) \\
    &= \lim_{n} \sum_{ij} L_{\xi_{i},\xi_{j}}(f_{j,n}^{*}*f_{i,n}) \\
    &= \lim_{n} \rip<\alpha_{n},\alpha_{n}>,
  \end{align*}
  where $\alpha_{n}:=\sum_{i}f_{i,n}\tensor \xi_{i}$.
  Since$\rip<\cdot,\cdot>$ is positive on $\ccgatho$ by
  Lemma~\vref{lem-alg-tensor}, we have $\rip<\alpha_{n},
  \alpha_{n}>\ge0$ and we've shown that $\rip<\cdot,\cdot>$ is still
  positive on $\bbocgatho$.

  Clearly the map sending the class $f\tensor\xi+\N$ to
  $f\tensor\xi+\N_{b}$ is isometric and therefore extends to an
  isometric embedding of $\H$ into the Hilbert space completion
  $\H_{b}$ of $\bbocgatho$ with respect to $\rip<\cdot,\cdot>$.
  However if $g\tensor\xi \in\bbocgatho$ and if $\set{g_{n}}$ is a
  sequence in $\cc(G)$ such that $g_{n}\to g$ in the usual way, then
  \begin{equation*}
    \|(g_{n}\tensor\xi+\N_{b})-(g\tensor\xi+\N_{b})\|^{2} =
    L_{\xi,\xi}(g_{n}^{*}*g_{n} -g_{n}^{*}*g -g*g_{n}^{*}+g^{*}*g),
  \end{equation*}
  and this tends to zero by Lemma~\vref{lem-grm-ext}.  Thus the image
  of $\H$ in $\H_{b}$ is all of $\H_{b}$.  Consequently, we can
  identify the completion of $\bbocgatho$ with $\H$ and
  $\bbocgatho/\N_{b}$ with a subspace of $\H$ containing $\hoo$.
\end{proof}

The ``extra'' vectors provided by $\bbocgatho/\N_{b}$ are just enough
to allow us to use a bit of general nonsense about unbounded operators
to extend the domain of each $L(f)$.  More precisely, for $f\in
\cc(G)$, we can view $L(f)$ as an operator in $\H$ with domain
$\D(L(f))=\hoo$.  Then using part~\partref3 of
Theorem~\vref{thm-ren-4.2}, we see that
\begin{equation*}
  L(f^{*})\subset L(f)^{*}.
\end{equation*}
This implies that $L(f)^{*}$ is a densely defined operator.  Hence
$L(f)$ is closable \cite{con:course}*{Proposition~X.1.6}.
Consequently, the closure of the graph of $L(f)$ in $\H\times\H$ is
the graph of the closure $\overline{L(f)}$ of $L(f)$
\cite{con:course}*{Proposition~X.1.4}.

Suppose that $g\in\bbocg$.  Let $\set{g_{n}}$ be a uniformly bounded
sequence in $\cc(G)$ all with supports in a fixed compact set such
that $g_{n}\to g$ pointwise.  Then
\begin{equation}\label{eq:38}
  \|[g_{n}\tensor\xi]-[g\tensor\xi]\|^{2} =
  L_{\xi,\xi}(g_{n}^{*}*g_{n}-g^{*}*g_{n} -g_{n}^{*}*g +g*g). 
\end{equation}
However $\set{g_{n}^{*}*g_{n}-g^{*}*g_{n} -g_{n}^{*}*g +g*g}$ is
uniformly bounded and converges pointwise to zero.  Since the supports
are all contained in a fixed compact set, the left-hand side of
\eqref{eq:38} tends to zero by Lemma~\vref{lem-grm-ext}.  Similarly,
\begin{equation*}
  \|[f*g_{n}\tensor\xi]-[f*g\tensor\xi]\|^{2}\to 0.
\end{equation*}
If follows that
\begin{equation*}
  \bigl([g_{n}\tensor\xi,L(f)[g_{n}\tensor \xi]\bigr) \to
  \bigl([g\tensor\xi],[f*g\tensor\xi]\bigr) 
\end{equation*}
in $(\bbocgatho/\N_{b})\times (\bbocgatho/\N_{b})\subset\H\times\H$.
Therefore $[g\tensor\xi]\in\D\bigl(\overline{L(f)}\bigr)$ and
$\overline{L(f)}[g\tensor\xi]= [f*g\tensor\xi]$.  We have proved the
following.

\begin{lemma}
  \label{lem-ext-lb}
  For each $f\in \cc(G)$, $L(f)$ is a closable operator in $\H$ with
  domain $\D(L(f))=\hoo=\ccgatho/\N$.  Furthermore $\bbocgatho/\N_{b}$
  belongs to $\D\bigl(\overline{L(f)}\bigr)$, and
  \begin{equation*}
    \overline{L(f)}[g\tensor\xi]=[f*g\tensor\xi]\quad\text{for all
      $f\in \cc(G)$, $g\in\bbocg$ and $\xi\in\H_{0}$.}
  \end{equation*}
\end{lemma}

Now can extend $L$ a bit.

\begin{lemma}
  \label{lem-lb}
  For each $f\in\bbocg$, there is a well-defined operator $L_{b}(f)\in
  \operatorname{Lin}(\bbocgatho)/N_{b})$ such that
  \begin{equation}
    \label{eq:39}
    L_{b}(f)[g\tensor\xi]=[f*g\tensor\xi].
  \end{equation}
  If $f\in\cc(G)$, then $L_{b}(f)\subset\overline{L(f)}$.
\end{lemma}
\begin{proof}
  To see that \eqref{eq:39} determines a well-defined operator, we
  need to see that
  \begin{equation}
    \label{eq:40}
    \sum_{i}[g_{i}\tensor\xi_{i}]=0\quad\text{implies}\quad \sum_{i}
    [f*g_{i}\tensor \xi_{i}]=0.
  \end{equation}
  However,
  \begin{equation}\label{eq:41p}
    \bigl\|\sum_{i} [f*g_{i}\tensor\xi_{i}]\bigr\|^{2}=\sum_{ij}
    L_{\xi_{i},\xi_{j}} (g_{j}^{*}*f^{*}*f*g_{i}).
  \end{equation}
  Since $f\in\bbocg$, we can approximate the right-hand side of
  \eqref{eq:41p} by sums of the form
  \begin{equation}
    \label{eq:42}
    \sum_{ij} L_{\xi_{i},\xi_{j}} (g_{j}^{*}*h^{*}*h*g_{i}),
  \end{equation}
  where $h\in\cc(G)$.  But \eqref{eq:42} equals
  \begin{equation*}
    \bigl\| \overline{L(h)} \sum_{i}[g_{i}\tensor\xi_{i}]\bigr\|^{2}
  \end{equation*}
  which is zero if the left-hand side of \eqref{eq:40} is zero.  Hence
  the right-hand side of \eqref{eq:40} is also zero and $L_{b}(f)$ is
  well-defined.

  If $f\in\cc(G)$, then $L_{b}(f)\subset \overline{L(f)}$ by
  Lemma~\vref{lem-ext-lb}.
\end{proof}

The previous gymnastics have allowed us to produce some additional
vectors in $\H$ and to extend slightly the domain of $L$.  The next
lemma provides the technical assurances that, despite the subtle
definitions above, our new operators act via the formulas we expect.

\begin{lemma}
  \label{lem-paul-tech}
  Suppose that $f\in\bbocg$ and that $k$ is a bounded Borel function
  on $\go$ which is the pointwise limit of a uniformly bounded
  sequence from $C_{0}(\go)$.  Then for all $g,h\in\cc(G)$ and
  $\xi,\eta\in\H_{0}$, we have the following.
  \begin{align}
    \bip(L_{b}(f)[g\tensor\xi]|[h\tensor\eta]) &=
    \bip([f*g\tensor\xi]|[ h\tensor\eta]) \tag{a}\\
    &= L_{\xi,\eta}(h^{*}*f*g) \notag\\
    &= L_{[g\tensor\xi],[h\tensor\eta]}(f)\notag \\[1ex]
    \big(M(k)[g\tensor\xi]|[h\tensor\eta]) &=
    L_{\xi,\eta}(h^{*}*((k\circ r)\cdot g)) \tag{b}\\
    &= \bip([(k\circ r)\cdot g\tensor \xi]|[h\tensor\eta])\notag \\
    &= \bip(M(k)L(g)\xi|L(h)\eta)\notag\\[1ex]
    \bip(M(k)L_{b}(f)[g\tensor\xi]|[h\tensor\eta]) &=
    \bip(L_{b}((k\circ r)\cdot f)[g\tensor\xi]|[h\tensor\eta]).\tag{c}
  \end{align}
\end{lemma}
\begin{proof}
  We start with \partref1.  The first equality is just the definition
  of $L_{b}(f)$.  The second follows from the definition of the inner
  product on $\bbocgatho/\N_{b}$.  If $f$ is in $\cc(G)$, then the
  third equation holds just by untangling the definition of the
  complex Radon measure $L_{\xi,\eta}$ and using the continuity in the
  inductive limit topology. The third equality holds for $f\in\bbocg$
  by applying the continuity assertion in Lemma~\vref{lem-grm-ext}.

  Part~\partref2 is proved similarly.  The first equation holds if
  $k\in C_{0}(\go)$ by definition of $M(k)$ and $L_{\xi,\eta}$.  If
  $\set{k_{n}}\subset C_{0}(\go)$ is a bounded sequence converging
  pointwise to $k$, then $M(k_{k})\to M(k)$ in the weak operator
  topology by the dominated convergence theorem.  On the other hand
  $h^{*}*(k_{n}\circ r)\cdot g\to h^{*}*(k\circ r)g$ in the required
  way.  Thus $L_{\xi,\eta}(h^{*}*(k_{n}\circ r)\cdot g) \to
  L_{\xi,\eta}(h^{*}*(k\circ r)g)$ by Lemma~\vref{lem-grm-ext}.  Thus
  the first equality is valid.  The second equality is clear if $k\in
  C_{0}(\go)$ and passes to the limit as above.  The third equality is
  simply our identification of $[g\tensor\xi]$ with $L(g)\xi$ as in
  Lemma~\vref{lem-alg-tensor}.

  For part~\partref3, first note that if $f_{n}\to f$ and $k_{n}\to k$
  are uniformly bounded sequences converging pointwise with supports
  in fixed compact sets independent of $n$, then $(k\circ r)\cdot
  f=\lim_{n} (k_{n}\circ r)\cdot f_{n}$.  It follows that $(k\circ
  r)\cdot f\in\bbocg$.  Also, $[f\tensor\xi]=\lim[f_{n}\tensor \xi]$,
  and since $M(k)$ is bounded, part~\partref2 implies that
  \begin{align*}
    M(k)[f\tensor\xi]&=\lim_{n}M(x)[f_{n}\tensor\xi] \\
    &= \lim_{n}[(k\circ r)\cdot f_{n}\tensor\xi]\\
    &= [(k\circ r)\cdot f\tensor \xi].
  \end{align*}
  Since it is not hard to verify that $M(k)^{*}[f\tensor\xi]= (\bar
  k\circ r)\cdot f\tensor \xi]$, we can compute that
  \begin{align*}
    \bip(M(k)L_{b}(f)[g\tensor\xi]|[h\tensor\eta]) &=
    \bip([f*g\tensor\xi] | (\bar k\circ r)\cdot h\tensor\eta]) \\
    &= \bip([k\circ r)\cdot(f*g)\tensor\xi]|[h\tensor\eta]) \\
    &=\bip([((k\circ r)\cdot f)*g\tensor \xi]|[h\tensor\eta]) \\
    &= \bip(L_{b}((k\circ r)\cdot
    f)[g\tensor\xi]|[h\tensor\eta]).\qed
  \end{align*}
\renewcommand\qed{}
\end{proof}

\begin{prop}
  \label{prop-quasi-inv}
  Let $\mu$ be the Radon measure on $\go$ associated to $L$ by
  Proposition~\vref{prop-measure-class}.  Then $\mu$ is
  quasi-invariant.
\end{prop}
\begin{proof}
  We need to show that measures $\nu$ and $\nu^{-1}$ (defined in
  \eqref{eq:10}~and \eqref{eq:18}, respectively) are equivalent.
  Therefore, we have to show that if $A$ is pre-compact in $G$, then
  $\nu(A)=0$ if and only if $\nu(A^{-1})=0$.  Since $(A^{-1})^{-1}=A$,
  it's enough to show that $A$ $\nu$-null implies that $A^{-1}$ is
  too.  Further, we can assume that $A\subset V$, where $V$ is open
  and Hausdorff.  Since $\nu\restr V$ is regular, we may as well
  assume that $A$ is a $G_{\delta}$-set so that $ f :=\charfcn A$ is
  in $\bboc(V)\subset \bbocg$.  Let $\tilde f (x)= f (x^{-1})$.  We
  need to show that $\tilde f (x)=0$ for $\nu$-almost all $x$.  Since
  $A$ is a $G_{\delta}$, we can find a sequence $\set{ f _{n}}\subset
  C_{c}^{+}(V)$ such that $ f _{n}\searrow f $.

  If $k $ is any function in $\cc(G)$, then $k f \bar k =|k |^{2} f
  \in\bbocg$ and vanishes $\nu$-almost everywhere.  By the monotone
  convergence theorem,
  \begin{equation*}
    \lambda(k  f \bar k )(u) 
    :=\int_{G}|k (\gamma)|^{2} f (\gamma)\,d\lambda^{u}(\gamma) 
  \end{equation*}
  defines a function in $\bboc(\go)$ which is equal to $0$ for
  $\mu$-almost all $u$.  In particular, $M(\lambda(k f \bar k ))=0$.

  It then follows from part~\partref2 of Lemma~\vref{lem-paul-tech}
  that
  \begin{equation}
    \label{eq:47}
    0=\bip(M\bigl(\lambda(k  f \bar k )\bigr)L(g)\xi|L(g)\xi) =
    L_{\xi,\xi} (g^{*}*\bigl(\lambda(k  f \bar k )\circ r)\cdot
    g\bigr) 
  \end{equation}
  for all $g\in\cc(G)$ and $\xi\in\H_{0}$.  On the other hand, if
  \eqref{eq:47} holds for all $g,k\in\cc(G)$ and $\xi\in\H_{0}$, then
  we must have $M\bigl(\lambda(k f \bar k )\bigr)=0$ for all $k
  \in\cc(G)$.  Since $ f (\gamma)\ge0$ everywhere, this forces $|k
  (\gamma)|^{2} f (\gamma)=0$ for $\nu$-almost all $\gamma$.  Since $k
  $ is arbitrary, we conclude that $ f (\gamma)=0$ for $\nu$-almost
  all $\gamma$.  Therefore it will suffice to show that
  \begin{equation}
    \label{eq:48p}
    L_{\xi,\xi}\bigl(g^{*}*\bigl(\lambda(k \tilde f \bar k )\circ
    r\bigr)\cdot g\bigr)=0  \quad\text{for all $g,k\in\cc(G)$ and
      $\xi\in\H_{0}$,}
  \end{equation}
  where we have replaced $ f $ with $\tilde f $ in the right-hand side
  of \eqref{eq:47}.  First, we compute that with $ f $ in
  \eqref{eq:47} we have
  \begin{align}
    g^{*}*\bigl(\lambda(k f &\bar k )\circ r\bigr)\cdot g(\sigma) =
    \int_{G}\overline{g(\gamma^{-1})}\bigl(\lambda(k f \bar k )\circ
    r\bigr)
    \cdot g(\gamma^{-1}\sigma) \,d\lambda^{r(\sigma)}(\gamma)\notag \\
    &= \int_{G } \overline{g(\gamma^{-1})} \lambda(k f \bar k
    )\bigl(s(\gamma)\bigr)
    g(\gamma^{-1}\sigma) \,d\lambda^{r(\sigma)}(\gamma)\notag\\
    &=\int_{G}\int_{G} \overline{g(\gamma^{-1})} \overline{k (\eta)} f
    (\eta)k (\eta) g(\gamma^{-1}\sigma)
    \,d\lambda^{s(\gamma)}(\eta)\,d\lambda^{r(\sigma)}(\gamma)
    \notag\\
    \intertext{which, after sending $\eta\mapsto \gamma^{-1}\eta$ and
      using left-invariance of the Haar system, is} &=\int_{G}\int_{G}
    \overline{g(\gamma^{-1})} \overline{k (\gamma^{-1}\eta)} f
    (\gamma^{-1}\eta) k (\gamma^{-1}\eta) g(\gamma^{-1}\sigma)
    \,d\lambda^{r(\sigma)}(\eta)\,d\lambda^{r(\sigma)}(\gamma) \notag\\
    \intertext{which, after defining $F(\gamma,\eta):=k
      (\gamma^{-1}\eta)g(\gamma^{-1})$ and $ f \cdot F(\gamma,\eta):=
      f (\gamma^{-1}\eta)F(\gamma,\eta)$ for $(\gamma,\eta)\in G\starr
      G$, is} &=\int_{G}\int_{G}\overline{F(\gamma,\eta)} f \cdot
    F(\sigma^{-1}\gamma,\sigma^{-1}\eta) \,d\lambda^{r(\sigma)}(\eta)
    \, d\lambda^{r(\sigma)}(\gamma).\label{eq:49p}
  \end{align}

  We will have to look at integrals of the form \eqref{eq:49p} in some
  detail.  First note that if $U$ and $V$ are Hausdorff open sets in
  $G$, then $U\starr V$ is a Hausdorff open set in $G\starr G$.  Thus
  if $g,k\in \cc(G)$, then
  $F(\gamma,\eta):=k(\gamma^{-1}\eta)g(\gamma^{-1})$ defines an
  element $F\in \cc(G\starr G)$.\footnote{For example, we can assume
    that $k\in C_{c}(U)$ and $g\in C_{c}(V)$ with $U$ and $V$ both
    open and Hausdorff.  A partition of unity argument as in the proof
    of Proposition~\vref{prop-star-alg}, allows us to assume that $VU$
    is Hausdorff.  Then observe that $\supp F\subset VU\starr
    V^{-1}$.}

  \begin{lemma}
    \label{lem-sp-F-one}
    Suppose that $F_{1},F_{2}\in\cc(G\starr G)$.  Then
    \begin{equation*}
      \sigma\mapsto
      \int_{G}\int_{G}\overline{F_{1}(\gamma,\eta)}
   F_{2}(\sigma^{-1}\gamma,\sigma^{-1}\eta)\,
      d\lambda^{r(\sigma)}(\eta) 
      \, d\lambda^{r(\sigma)}(\gamma)
    \end{equation*}
    defines an element of $\cc(G)$ which we denote by
    $\overline{F}_{1}*_{\lambda*\lambda}F_{2}$.
  \end{lemma}
  \begin{proof}
    We can take $F_{i}\in C_{c}(U_{i}\starr V_{i})$ with each $U_{i}$
    and $V_{i}$ open and Hausdorff.  As in the proof of
    Proposition~\vref{prop-star-alg}, we can assume that
    $U_{1}U_{2}^{-1}$ and $V_{1}V_{2}^{-1}$ are Hausdorff.  Note that
    \begin{equation*}
      \|\overline{F}_{1}*_{\lambda*\lambda}F_{2}\|_{\infty}\le
      \|F_{1}\|_{\infty}
      \|F_{2}\|_{\infty}\sup_{u\in\go}\lambda^{u}(K_{1})\lambda^{u}(K_{2})
    \end{equation*}
    whenever $\supp F_{1}\subset K_{1}\starr K_{2}$.  Thus to see that
    $\overline{F}_{1}*_{\lambda*\lambda} F_{2}\in
    C_{c}(U_{1}U_{2}^{-1}\cap V_{1}V_{2}^{-1})$, it will suffice to
    consider only those $F_{i}$ is dense subspaces of
    $C_{c}(U_{1}\starr V_{1})$ and $C_{c}(U_{2}\starr V_{2})$.  In
    particular, we can assume that each $F_{i}$ is of the form
    $F_{i}(\gamma,\eta)= k_{i}(\eta)g(\gamma^{-1})$.  But then
    \begin{equation*}
      \overline{F}_{1}*_{\lambda*\lambda}F_{2}(\sigma)=
      \overline{k}_{1}*\tilde k_{2}(\sigma) g_{1}^{*}*g_{2}(\sigma),  
    \end{equation*}
    and the result follows.
  \end{proof}

  \begin{lemma}
    \label{lem-denseb}
    Functions of the form
    \begin{equation}
      \label{eq:66}
      (\gamma,\eta)\mapsto
      k(\gamma^{-1}\eta)g(\gamma^{-1})\quad\text{with $k,g\in\cc(G)$}
    \end{equation}
    span a dense subspace of $\cc(G\starr G)$ in the \ilt.
  \end{lemma}
  \begin{proof}
    We have already noted that functions of the form given in
    \eqref{eq:66} determine elements $\theta_{k,g}$ in $\cc(G\starr
    G)$.  Furthermore, arguing as in the proof of
    Proposition~\vref{prop-star-alg}, it will suffice to show that we
    can approximate functions $\theta\in C_{c}(U\starr V)$ with $U$
    and $V$ open, Hausdorff and such that $UV$ is Hausdorff.  Then the
    span of functions $\theta_{k,g}$ with $k\in C_{c}(UV)$ and $g\in
    C_{c}(V^{-1})$ is dense in $C_{c}(U\starr V)$ in the \ilt{} by the
    Stone-Weierstrass Theorem.
  \end{proof}

  Let $\AA\subset\sa_{c}(G\starr G;\iota^{*}\B)$ be the dense subspace
  of functions of the form considered in Lemma~\vref{lem-denseb}.  We
  continue to write $ f $ for the characteristic function of our fixed
  pre-compact, $\nu$-null set.  Then we know from \eqref{eq:47} that
  \begin{equation}
    \label{eq:50p}
    L_{\xi,\xi}\bigl(\overline{F}*_{\lambda*\lambda}( f \cdot
    F)\bigr)=0\quad\text{for all $F\in\AA$.}
  \end{equation}
  It is not hard to check that, if $ f '\in\bboc(G)$, then
  $F*_{\lambda*\lambda}( f '\cdot F)\in \bbocg$ and that if $F_{n}\to
  F$ in the \ilt{} in $\cc(G\starr G)$, then
  $\set{\overline{F}_{n}*_{\lambda*\lambda}( f '\cdot F)}$ is
  uniformly bounded and converges pointwise to
  $\overline{F}*_{\lambda*\lambda}( f '\cdot F)$.  In particular the
  continuity of the $L_{\xi,\xi}$ (see Lemma~\vref{lem-grm-ext})
  implies that \eqref{eq:50p} holds for all $F\in\cc(G\starr G)$.  But
  if we define $\tilde F(x,y):=F(y,x)$, then we see from the
  definition that
  \begin{equation*}
    \overline{\tilde F}*_{\lambda*\lambda}( f \cdot \tilde F)=
    \overline{F}*_{\lambda*\lambda}(\tilde  f \cdot F)),
  \end{equation*}
  where we recall that $\tilde f (x):= f (x^{-1})$.  Thus
  \begin{equation*}
    L_{\xi,\xi}\bigl(\overline{F}*_{\lambda*\lambda}(\tilde  f \cdot F))\bigr)
    =0 \quad\text{for all $F\in\cc(G\starr G)$}.
  \end{equation*}
  Since the above holds in particular for $F\in\AA$, this implies
  \eqref{eq:48p}, and completes the proof.
\end{proof}

To define the Borel Hilbert bundle we need, we need to see that the
complex Radon measures $L_{\xi,\eta}$ defined above are absolutely
continuous with respect to the measure $\nu$.  In order to prove this,
we need to restrict $\xi$ and $\eta$ to lie in $\hoo$, and to employ
Lemma~\vref{lem-paul-tech}. 

\begin{lemma}
  \label{lem-abs-cont}
  Let $a$ and $b$ be vectors in $\hoo$ (identified with
  $\ccgatho/\N$).  Let $L_{a,b}$ be the complex Radon measure given by
  \begin{equation*}
    L_{a,b}(f):=\bip(L(f)a|b).
  \end{equation*}
  Then $L_{a,b}$ is absolutely continuous with respect to the measure
  $\nu$ defined in \eqref{eq:10}.\footnote{Absolute continuity of
    complex Radon measures on \lhlc{} spaces is discussed in
    Appendix~\ref{sec:radon-meas-locally}.}
\end{lemma}

\begin{proof}
  It is enough to show that if $M$ is a pre-compact $\nu$-null set and
  if $f:=\charfcn M$, then $L_{a,b}(f)=0$.  We can also assume that
  $M\subset V$ where $V$ is a Hausdorff open set.  Since $\nu\restr V$
  is a Radon measure, and therefore regular, we may as well assume
  that $M$ is a $G_{\delta}$-set.  Then $f\in \bboc(V)\subset \bbocg$.

  On the other hand,
  \begin{equation*}
    0=\int_{\go}\int_{G}  f (\gamma)\,d\lambda^{u}(\gamma)\,d\mu(u),
  \end{equation*}
  so there is a $\mu$-null set $N\subset\go$ such that
  $\lambda^{u}(M\cap G^{u})=0$ if $u\notin N$.  As above, we can
  assume that $N$ is a $G_{\delta}$ set.  Then for any $g\in\cc(G)$,
  we have
  \begin{equation*}
    f*g(\gamma)=\int_{G}f(\eta)g(\eta^{-1}\gamma)\,d\lambda^{r(\gamma)}(\eta)=0
  \end{equation*}
  whenever $r(\gamma)\notin N$.  Since $\supp
  \lambda^{r(\gamma)}=G^{r(\gamma)}$, it follows that for all
  $\gamma\in G$ (without exception),
  \begin{equation}\label{eq:46p}
    f*g(\gamma)= \charfcn N\bigl(r(\gamma)\bigr) f*g(\gamma) =
    \bigl((\charfcn N\circ 
    r)\cdot f\bigr)*g(\gamma).  
  \end{equation}

  Since $a,b\in\hoo$, it suffices to consider $a=[g\tensor\xi]$ and
  $b=[h\tensor \eta]$ (with $g,h\in\cc(G)$ and $\xi,\eta\in\H_{0})$.
  Note that $f$ and $\charfcn N$ satisfy the hypotheses of
  Lemma~\vref{lem-paul-tech}.  Therefore, by part~\partref1 of that
  lemma,
  \begin{align*}
    L_{[g\tensor\xi,h\tensor\eta}(f)&=\bip([f*g\tensor
    \xi]|[h\tensor\eta]) \\
    \intertext{which, by \eqref{eq:46p}, is}
    &= \bip([((\charfcn N\circ r)\cdot f)*g\tensor\xi]|[h\tensor\eta]) \\
    \intertext{which, by part~\partref1 of Lemma~\ref{lem-paul-tech},
      is} &= \bip(L_{b}((\charfcn N\circ r)\cdot
    f)[g\tensor\xi]|[h\tensor\eta])
    \\
    \intertext{which, by part~\partref3 of Lemma~\ref{lem-paul-tech},
      is}
    &= \bip(M(\charfcn N)L_{b}(f)[g\tensor\xi]|[h\tensor\eta]).
  \end{align*}
Since $M(\charfcn N)=0$, the last inner product is zero
  as desired.  This completes the proof.
\end{proof}

 Since the measures $\nu$ and $\nu_{0}$ are equivalent, for each
$\xi,\eta\in\hoo$, we can, in view of Lemma~\vref{lem-abs-cont}, let
$\rho_{\xi,\eta}$ be the Radon-Nikodym derivative of $L_{\xi,\eta}$
with respect to $\nu_{0}$.  Then for each $\xi,\eta\in\hoo$, we have
\begin{align*}
  \bip(L(f)\xi|\eta)&=
  L_{\xi,\eta}(f) = \int_{G} f(x)
  \,d L_{\xi,\eta}(x)
  \\
  &= \int_{G}f(x) \rho_{\xi,\eta}(x)
  \Delta(x)^{-\half} \,d\nu(x)
  \\
  &=\int_{\go}\int_{G}f(x)
  \rho_{\xi,\eta}(x) \Delta(x)^{-\half}\,d\lambda^{u}(x)\,d\mu(u) .
\end{align*}

Our next computation serves to motivate the construction in
Lemma~\vref{lem-tensor-fix}.  If $\xi,\eta\in\hoo$, then we can 
apply Lemma~\vref{lem-abs-cont} and compute that
\begin{align*}
  \bip(L(f)\xi|L&(g)\eta)= \bip(L(g^{*}*f)\xi|\eta) =L_{\xi,\eta}(g^{*}*f) \\
  &= \int_{\go}\int_{G}g^{*}*f(\gamma) \rho_{\xi,\eta}(\gamma)
  \Delta(\gamma)^{-\half} \,d\lambda^{u}(\gamma)\,d\mu(u) \\
  &= \int_{\go}\int_{G}\int_{G}
  \overline{g(\eta^{-1})}f(\eta^{-1}\gamma)
  \rho_{\xi,\eta}(\gamma)\Delta(\gamma)^{-\half}
  \,d\lambda^{u}(u)\,d\lambda^{u}(\gamma)\,d\mu(u)
  \\
  \intertext{which, by Fubini and sending $\gamma\mapsto \eta\gamma$, is} &=
  \int_{\go}\int_{G}\int_{G}
  \overline{g(\eta^{-1})}f(\gamma) \rho_{\xi,\eta}(\eta\gamma)
  \Delta(\eta\gamma)^{-\half} \, d\lambda^{s(\eta)}(\gamma) \,d\lambda^{u}(\eta)\,d\mu(u)
  \\
  \intertext{which, after sending $\eta\mapsto \eta^{-1}$, and using the
    symmetry of $\nu_{0}$, is} &= \int_{\go}\int_{G}\int_{G}
  \overline{g(\eta)}f(\gamma) \rho_{\xi,\eta}(\eta^{-1}\gamma)
  \Delta(\eta)^{-\half} \Delta(\gamma)^{-\half} \\
  &\hskip2.5in \,d\lambda^{u}(\gamma)\,d\lambda^{u}(\eta)\, d\mu(u).
\end{align*}

Since it is not clear to what extent $\rho_{\xi,\eta}$ is a
sesquilinear function of $(\xi,\eta)$, we fix once and for all a
countable orthonormal basis $\set{\spxi_{i}}$ for
$\hoo$.  (Actually, any countable linearly independent set whose span
is dense in $\hoo$ will do.)  We let
\begin{equation*}
  \hoop:=\operatorname{span}\set{\spxi_{i}}.
\end{equation*}
To make the subsequent formulas a bit easier to read, we will write
$\rho_{ij}$ in place of the Radon-Nikodym derivative
$\rho_{\spxi_{i},\spxi_{j}}$.  The linear independence
of the $\spxi_{i}$ guarantees that each $\alpha\in\cc(G)\atensor\hoo'$
can 
be written \emph{uniquely} as
\begin{equation*}
  \alpha=\sum_{i} f_{i}\tensor\spxi_{i}
\end{equation*}
where all but finitely many $f_{i}$ are zero.

\begin{lemma}
  \label{lem-tensor-fix}
  For each $u\in\go$, there is a sesquilinear form $\ipu<\cdot,\cdot>$
  on $\ccgathoop$ such that
  \begin{equation}
    \label{eq:52}
    \ipu<f\tensor\spxi_{i},g\tensor \spxi_{j}> = \int_{G}\int_{G}
    \overline{g(\eta)}f(\gamma) \rho_{ij}(\eta^{-1}\gamma)
    \Delta(\eta\gamma)^{-\half} \,d\lambda^{u}(\gamma)\,d\lambda^{u}(\eta).
  \end{equation}
  Furthermore, there is a $\mu$-conull set $F\subset \go$ such that
  $\ipu<\cdot,\cdot>$ is a pre-inner product for all $u\in F$.
\end{lemma}
\begin{remark}
  \label{rem-tensor-prob}
  As mentioned earlier, we fixed the $\spxi_{i}$ because it isn't
  clear that the right-hand side of \eqref{eq:52} is linear in
  $\spxi_{i}$ or conjugate linear in $\spxi_{j}$.
\end{remark}

\begin{proof}
  Given $\alpha=\sum_{i}f_{i}\tensor \spxi_{i}$ and $\beta=\sum_{j}
  g_{j}\tensor \spxi_{j}$, we get a well-defined form via the definition
  \begin{equation*}
    \ipu<\alpha,\beta> = 
    \sum_{ij} \int_{G}\int_{G}
    \overline{g_{j}(\eta)}f_{i}(\gamma) 
    \rho_{ij}(\eta^{-1}\gamma) \Delta(\eta\gamma)^{-\half}
    \,d\lambda^{u}(\gamma)\,d\lambda^{u}(\eta).  
  \end{equation*}
  This clearly satisfies \eqref{eq:52}, and is
  linear in $\alpha$ and conjugate linear in $\beta$.  It only remains
  to provide a conull Borel set $F$ such that $\ipu<\cdot,\cdot>$ is
  positive for all $u\in F$.

  However, \eqref{eq:52} was inspired by the calculation preceding the
  lemma.  Hence if $\alpha:=\sum_{i}f_{i}\tensor\spxi_{i}$, then
  \begin{equation}\label{eq:55p}
    \begin{split}
      \Bigl\|\sum L(f_{i})\spxi_{i}\Bigr\|^{2}&=
      \sum_{ij} \bip(L(f_{i})\spxi_{i}|{L(f_{j})\spxi_{j}}) \\
      &= \sum_{ij} \bip(L(f_{j}^{*}*f_{i})\spxi_{i}|\spxi_{j}) \\
      &= \sum_{ij} \int_{\go}\int_{G} \int_{G} 
      \overline{f_{j}(\eta)}f_{i}(\gamma)  \rho_{ij}(\eta^{-1}\gamma)
      \Delta(\gamma\eta)^{-\half}\\ 
      &\hskip2in
      \,d\lambda^{u}(\gamma) \, d\lambda^{u}(\eta) \,d\mu(u) \\
      &= \sum_{ij} \int_{\go} \ipu<f_{i}\tensor \spxi_{i},f_{j}\tensor
      \spxi_{j}>\,d\mu(u) \\
      &= \int_{\go} \ipu<\alpha,\alpha>\,d\mu(u).
    \end{split}
  \end{equation}
  Thus, for $\mu$-almost all $u$, we have $\ipu<\alpha,\alpha>\ge0$.
  The difficulty is that the exceptional null set depends on $\alpha$.
  However, 
  there is a sequence $\set{f_{i}}\subset\cc(G)$ which is dense in
  $\cc(G)$ in the \ilt.
  Let $\AA$ be the rational vector space spanned by the countable set
  $\set{f_{i}\tensor\spxi_{j}}_{i,j}$.  Since $\AA$ is countable,
  there is a $\mu$-conull set $F$ such that $\ipu<\cdot,\cdot>$ is a
  positive $\mathbf{Q}$-sesquilinear form on $\AA$.  However, if
  $g_{i}\to g$ and $h_{i}\to h$ in the \ilt{} in $\cc(G)$, then, since
  $\lambda^{u}\times\lambda^{u}$ is a Radon measure on $G^{u}\times G^{u}$, we
  have $\ipu<g_{i}\tensor\spxi_{j}, h_{i}\tensor \spxi_{k} > \to
  \ipu<g\tensor\spxi_{j},h\tensor\spxi_{k}>$.  It follows that for all
  $u\in F$, $\ipu<\cdot,\cdot>$ is a positive sesquilinear form (over
  $\C$) on the complex vector space generated by
  \begin{equation*}
    \set{f\tensor \spxi_{i}:f\in\cc(G)}.
  \end{equation*}
  However, as that is all $\ccgathoop$, the proof is complete.
\end{proof}

We need the following technical result which is a rather specialized
 version of the
Tietze Extension Theorem for \lhlc{} spaces.

\begin{lemma}
  \label{lem-tietze}
  Suppose that $g\in C_{c}(G^{u})$ for some $u\in\go$.  Then there is
  a $f\in \cc(G)$ such that $f\restr{G^{u}}=g$.
\end{lemma}
\begin{proof}
  There are Hausdorff open sets $V_{1},\dots,V_{n}$ such that $\supp
  g\subset \bigcup V_{i}$. Then, using a partition of unity, we can
  find
 $g_{i}\in C_{c}(G^{u})$ such
  that $\supp g_{i}\subset V_{i}$ and such that $\sum g_{i}=g$.  By
  the Tietze Extension Theorem, there are $f_{i}\in C_{c}(V_{i})$ such
  that $f_{i}\restr{G^{u}}=g_{i}$.  Then $f:=\sum f_{i}$ does the
  job. 
\end{proof}

Note that for \emph{any} $u\in\go$, the value of
$\ipu<f\tensor\spxi_{i},g \tensor\spxi_{j}>$ depends only on $f\restr
{G^{u}}$ and $g\restr{G^{u}}$.  Furthermore, using our specialized
Tietze Extension result above, we can view $\ipu<\cdot,\cdot>$ as a
sesquilinear form on $C_{c}(G^{u})$.  (Clearly, since $G^{u}$ is
Hausdorff, each $f\in\cc(G)$ determines an element of $C_{c}(G^{u})$.
We need Lemma~\vref{lem-tietze} to know that every function in
$C_{c}(G^{u})$ arises in this fashion.)  In particular, if $f\in
\cc(G)$ and $\sigma\in G$, then we let $\hatpi(\sigma)f$ be any
element of $\cc(G)$ such that
\begin{equation*}
  (\hatpi(\sigma)f)(\gamma)=\Delta(\sigma)^{\half}f(\sigma^{-1}\gamma)\quad\text{for all 
  $\gamma\in G^{r(\sigma)}$.}
\end{equation*}
Of course, $\hatpi(\sigma)f$ is only well-defined on $G^{r(\sigma)}$.

The next computation is critical to what follows.  We have
\begin{align*}
  \bipb r(\sigma)<\hatpi(\sigma)f\tensor\spxi_{i},g\tensor\spxi_{j}&> = \int_{G}\int_{G}
  \overline{g(\eta)}f(\sigma^{-1}\gamma)
  \rho_{ij} (\gamma^{-1}\eta) \Delta(\sigma^{-1}\gamma\eta)^{-\half}\\
  &\hskip2.25in \, d\lambda^{r(\sigma)}(\eta) \,
  d\lambda^{r(\sigma)} (\gamma) \\
  \intertext{which, after sending $\gamma\mapsto \sigma\gamma$, is} &=
  \int_{G}\int_{G} \overline{g(\eta)}f(\gamma)
  \rho_{ij}(\gamma^{-1}\sigma^{-1}\eta) \Delta(\gamma\eta)^{-\half} \\
  &\hskip2.25in
  \,d\lambda^{r(\sigma)}(\eta)\,d\lambda^{s(\sigma)}(\gamma) \\
  \intertext{which, after sending $\eta\mapsto \sigma\eta$, is} &=
  \int_{G}\int_{G} \overline{g(\sigma\eta)}f(\gamma)
  \rho_{ij}(\gamma^{-1}\eta) \Delta(\sigma)^{-\half}\Delta(\gamma\eta)^{-\half}\\
  &\hskip2.25in \,
  d\lambda^{s(\sigma)}(\eta) \,d\lambda^{s(\sigma)}(\gamma) \\
  &= \int_{G}\int_{G}
  \bigl(\hatpi(\sigma^{-1})g\bigr)(\eta)f(\gamma)
   \rho_{ij}(\gamma^{-1}\eta) \Delta(\gamma\eta)^{-\half} \,\\
  &\hskip2.25in d\lambda^{s(\sigma)}(\eta) \,d\lambda^{s(\sigma)}(\gamma) \\
  &= \bipb s(\sigma)<f\tensor\spxi_{i},\hatpi(\sigma^{-1})g\tensor\eta>.
\end{align*}

Recall that $G$ acts continuously on the left of
$\go$: $\gamma\cdot s(\gamma)=r(\gamma)$.  In particular, if $C$ is compact
in $G$ and if $K$ is compact in $\go$, then
\begin{equation*}
  C\cdot K=\set{\gamma\cdot u:(\gamma,u)\in G^{(2)}\cap(C\times K)}
\end{equation*}
is compact.  If $U\subset \go$, then we say that $U$ is
\emph{saturated} if $U$ is $G$-invariant.  More simply, $U$ is
saturated if $s(x)\in U$ implies $r(x)$ is in $U$.  If $V\subset\go$,
then its \emph{saturation} is the set $[V]=G\cdot V$ which is the
smallest saturated set containing $V$.  

The next result is a key technical step in our proof and takes the
place of the Ramsay selection theorems
(\cite{ram:jfa82}*{Theorem~3.2} and \cite{ram:am71}*{Theorem~5.1})
used in Muhly's and Renault's proof.
\begin{lemma}
  \label{lem-saturated}
  We can choose the $\mu$-conull Borel set $F\subset\go$ in
  Lemma~\ref{lem-tensor-fix} to be saturated for the $G$-action on
  $\go$. 
\end{lemma}
\begin{proof}
  Let $F$ be the Borel set from Lemma~\vref{lem-tensor-fix}.  We want
  to see that $\ipv<\cdot,\cdot>$ is positive for all $v$ in the
  saturation of $F$.  To this end, suppose that $u\in F$ and that
  $\sigma\in G$ is such that $s(\sigma)=u$ and $r(\sigma)=v$.  Then
  \begin{equation*}
    \gamma\mapsto \Delta(\sigma)^{\half}  f(\sigma^{-1}\gamma)
  \end{equation*}
is in $C_{c}(G^{v})$, and such functions span a dense
subspace of $C_{c}(G^{v})$ in the \ilt.  Moreover, as we observed at the end
of the proof of Lemma~\vref{lem-tensor-fix},
\begin{equation*}
  \ipv<f_{i}\tensor\spxi_{j},g_{i}\tensor\spxi_{k}>\to
  \ipv<f\tensor\spxi_{j} , g\tensor\spxi_{k}>
\end{equation*}
provided $f_{i}\to f$ and $g_{i}\to g$ in the \ilt{} in
$C_{c}(G^{v})$. 
Therefore, to show that $\ipv<\cdot,\cdot>$ is positive, it will
suffice to check on vectors of the form $\alpha:=\sum_{i}
\hatpi(\sigma)(f_{i})\tensor \spxi_{i}$.  Then using the key calculation
preceding Lemma~\ref{lem-saturated}, we have
\begin{equation} \label{eq:67}
\begin{split}
  \ipv<\alpha,\alpha>&=\sum_{ij}
  \bipu<\hatpi(\sigma^{-1}\sigma)f_{i}\tensor\spxi_{i} ,
  f_{j}\tensor\spxi_{j}>. \\
&= \sum_{ij} \bipu<f_{i}\tensor\spxi_{i} ,
  f_{j}\tensor\spxi_{j}>. \\
&=\bipu<\sum_{i} f_{i}\tensor\spxi_{i}, \sum_{i}
f_{i}\tensor\spxi_{i}>
\end{split}
\end{equation}
which is positive since $u\in F$.

It only remains to verify that the saturation of $F$ is Borel.  Since
$\mu$ is a Radon measure --- and therefore regular --- we can shrink $F$ a bit, if necessary, and
assume it is $\sigma$-compact.  Say $F=\bigcup K_{n}$.  On the other
hand, $G$ is second countable and therefore $\sigma$-compact.  If
$G=\bigcup C_{m}$, then $[F]=\bigcup C_{m}\cdot K_{n}$.  Since each
$C_{m}\cdot K_{n}$ is compact, $[F]$ is $\sigma$-compact and therefore
Borel.  This completes the proof.
\end{proof}

From here on, we will assume that $F$ is saturated.  In view of
Lemma~\vref{lem-tensor-fix}, for each $u\in F$ we can define $\H(u)$
to be the Hilbert space completion of $\ccgathoop$ with respect to
$\ipu<\cdot,\cdot>$.  We will denote the image of $f\tensor\spxi_{i}$
in $\H(u)$ by $f\tensor_{u}\spxi_{i}$.  Since the complement of $F$ is
$\mu$-null and also saturated, what we do off $F$ has little
consequence.  In particular, $G$ is the disjoint union of $G\restr F$
and the $\nu$-null set $G\restr{\go\setminus F}$.\footnote{The
  saturation of $F$ is critical to what follows.  If $F$ is not
  saturated, then in general $G$ is not the union of $G\restr F$ and $
  G\restr{\go\setminus F}$.  But as $F$ \emph{is} saturated, note that
  a homomorphism $\phi:G\restr F\to H$ can be trivially extended to a
  homomorphism on all of $G$ by letting $\phi$ be suitably trivial on
  $G\restr{\go\setminus F}$.  This is certainly not the case if $F$ is
  not saturated.}  Nevertheless, for the sake of nicety, we let
$\mathcal{V}$ be a Hilbert space with an orthonormal basis
$\set{e_{ij}}$ doubly indexed by the same index sets as for
$\set{f_{i}}$ and $\set{\spxi_{j}}$, and set $\H(u)=\mathcal{V}$ if
$u\in\go\setminus F$.  We then let
\begin{equation*}
  \go*\HH=\set{(u,h):\text{$u\in F$ and $h\in\H(u)$}},
\end{equation*}
and define $\Phi_{ij}:F\to F*\HH$ by 
\begin{equation*}
  \Phi_{ij}(u):=
  \begin{cases}
     f_{i}\tensor_{u}
\spxi_{j}&\text{if $u\in F$ and} \\
e_{ij}&\text{if $u\notin F$.}
  \end{cases}
\end{equation*}
(Technically, $\Phi_{ij}(u)=(u, f_{i}\tensor_{u}\spxi_{j})$ --- at least
for $u\in F$ --- but we have agreed to obscure this subtlety.)  Then
\cite{wil:crossed}*{Proposition~F.8} implies that we can make $\go*\HH$
into a Borel Hilbert bundle over $\go$ in such a way that the
$\set{\Phi_{ij}}$ form a fundamental sequence (see
\cite{wil:crossed}*{Definition~F.1}).  Note that if
$f\tensor\spxi_{i}\in \ccgathoop$ and if
$\Phi(u):=f\tensor_{u}\spxi_{i}$, then
\begin{equation*}
  u\mapsto \bipu<\Phi(u),\Phi_{ij}(u)>
\end{equation*}
is Borel on $F$.\footnote{We can define $\Phi(u)$ to be zero off $F$.
  We are going to continue to pay as little attention as possible to the null
  complement of $F$ in the
  sequel.}  It follows that $\Phi$ is a Borel section of $\go*\HH$ and
defines a class in
$L^{2}(\go*\HH,\mu)$.

Furthermore, \eqref{eq:67} shows that for each $\sigma\in G\restr F$,
there is a unitary $U_{\sigma}:\H\bigl(s(\sigma)\bigr) \to
\H\bigl(r(\sigma)\bigr)$ characterized by
\begin{equation*}
  U_{\sigma}(f\tensor_{s(\sigma)}\spxi_{i})
=\hatpi(\sigma)f\tensor_{r(\sigma)}\spxi_{i}. 
\end{equation*}
If $\sigma\notin G\restr F$, then
$\H\bigl(s(\sigma)\bigr)=\H\bigl(r(\sigma)\bigr)=\mathcal V$, and we
can let $U_{\sigma}$ be the identity operator.

\begin{lemma}
  \label{lem-pi}
  The map $\hat U$ from $G$ to $ \operatorname{Iso}(\go*\HH)$ defined by $\hat
  U(\sigma):= \bigl(r(\sigma), U_{\sigma}, s(\sigma)\bigr)$ is a Borel
  homomorphism.  Hence $(\mu,\go*\HH,\hat U)$ is a unitary 
  representation of $G$ on $L^{2}(\go*\HH,\mu)$.
\end{lemma}
\begin{proof}
  If $\sigma\in G\restr F$, then
  \begin{multline*}
    \bip(U_{\sigma}\Phi_{ij}\bigl(s(\sigma)\bigr)
    |{\Phi_{kl}\bigl(r(\sigma)\bigr) }) = \\
\int_{G}\int_{G} \overline{f_{k}(\eta)}f_{i}(\sigma^{-1}\gamma)
 \rho_{jl}(\eta^{-1}\gamma)\Delta(\sigma^{-1}\gamma\eta)^{-\half} \,
d\lambda^{r(\sigma)}\, d\lambda^{r(\sigma)}.
  \end{multline*}
  Thus $\sigma\mapsto \bip(U_{\sigma}\Phi_{ij}\bigl(s(\sigma)\bigr)
  |{\Phi_{kl}\bigl(r(\sigma)\bigr) })$ is Borel on $F$ by Fubini's
  Theorem.  Since it is clearly Borel on the complement of $F$,
  $\hat U$ is Borel.  The algebraic properties are
  straightforward.  For example, assuming that $\gamma\in
G^{r(\sigma)}$, we have on the one hand,
\begin{equation*}
  \bigl(\hatpi(\sigma\eta)f\bigr)(\gamma)=\Delta(\sigma\eta)^{\half} 
  f\bigl( (\sigma\eta)^{-1}\gamma\bigr),
\end{equation*}
while
\begin{align*}
  \bigl(\hatpi(\sigma)\hatpi(\eta)f\bigr)(\gamma)&=
  \Delta(\sigma)^{\half}\bigl(\hatpi(\eta)f\bigr)(\sigma^{-1}\gamma\bigr) \\
&= \Delta(\sigma\eta)^{\half}  f\bigl(\eta^{-1}\sigma^{-1}\gamma\bigr).
\end{align*}
It follows that $\hat U$ is multiplicative on
$G\restr F$.  Of course, it is clearly multiplicative on the
complement (which is $G\restr{\go\setminus F}$ since $F$ is
saturated). 
\end{proof}

\begin{lemma}
  \label{lem-V}
  Each $f\tensor\spxi_{i}\in\ccgathoop$ determines a Borel section
  $\Phi(u):=f\tensor_{u}\spxi_{i}$ in $\L^{2}(\go*\HH,\mu)$ whose class in
  $L^{2}(\go*\HH,\mu)$ depends only on the class of
  $[f\tensor\spxi_{i}]\in \ccgathoop/\N\subset \ccgatho/\N=\hoo$.  Furthermore,
  there is a unitary isomorphism $V$ of $\H$ onto $L^{2}(\go*\HH,\mu)$
  such that $V(L(f)\spxi_{i})=[\Phi]$.
\end{lemma}
\begin{proof}
  We have already seen that $\Phi$ is in $\L^{2}(F*\HH,\mu)$.  More
  generally, the computation \eqref{eq:55p} in the proof of
  Lemma~\vref{lem-tensor-fix} shows that if $\alpha=\sum_{i}
  f_{i}\tensor\spxi_{i}$ and $\Psi(u):=\sum_{i}
  f_{i}\tensor_{u}\spxi_{i}$, then
  \begin{equation*}
    \|\Psi\|_{2}^{2}=\Bigl\|\sum_{i} L(f_{i})\spxi_{i}\Bigr\|^{2}
  \end{equation*}
  Thus there is a well defined isometric map $V$ as in the statement
  of lemma mapping $\operatorname{span}\set{L(f)\spxi_{i}:f\in\cc(G)}$
  onto a dense subspace of $L^{2}(F*\HH,\mu)$.  Since $\hoop$ is dense
  in $\hoo$, and therefore in $\H$, the result follows by
  Lemma~\vref{lem-dense}.
\end{proof}

The proof of Theorem~\vref{thm-ren-4.2} now follows almost
immediately from the next proposition.

\begin{prop}
  \label{prop-main}
  The unitary $V$ defined in Lemma~\vref{lem-V} intertwines $L$ with a
  representation $L'$ which in the integrated form of the unitary
  representation $(\mu,\go*\HH,U)$ from Lemma~\vref{lem-pi}.
\end{prop}
\begin{proof}
  We have $L'(f_{1})=VL(f_{1})V^{*}$.  On the one hand,
  \begin{multline*}
    \bip(L(f_{1})[f\tensor\spxi_{i}]|[g\tensor \spxi_{j}])_{\H}
    = \\ \bip(VL(f_{1})[f\tensor\spxi_{i}]|V[g\tensor\spxi_{j}]) =\\
    \bip(L'(f_{1})V[f\tensor\spxi_{1}]|V[g\tensor\spxi_{j}]).
  \end{multline*}
  But the left-hand side is
  \begin{align*}
    \bip(L(f_{1}&*f)\spxi_{i}|L(g)\spxi_{j}) =
    L_{\spxi_{i},\spxi_{j}}(g^{*}*f_{1}*f) \\
    &= \int_{\go}\int_{G}\int_{G}
    \overline{g(\eta)}f_{1}*f(\gamma)
    \rho_{ij}(\eta^{-1}\gamma)\Delta(\eta\gamma)^{-\half}
    \,d\lambda^{u}(\gamma)\,d\lambda^{u}(\eta)\,d\mu(u) \\
    &= \int_{\go}\int_{G}\int_{G}\int_{G}
    \overline{g(\eta)}f_{1}(\sigma)f(\sigma^{-1}\gamma)
    \rho_{ij}(\eta^{-1}\gamma) \Delta(\eta\gamma)^{-\half} \\
    &\hskip2in
    \,d\lambda^{u}(\sigma)
    \,d\lambda^{u}(\gamma)\,d\lambda^{u}(\eta)\,d\mu(u) \\ 
    &= \int_{\go}\int_{G}\int_{G}\int_{G}f_{1}(\sigma) \overline{g(\eta)}
    \bigl(\hatpi(\sigma)f\bigr)(\gamma)
    \rho_{ij}(\eta^{-1}\gamma) \Delta(\eta\gamma)^{-\half}
    \Delta(\sigma)^{-\half}\\ 
    &\hskip2in
    \,d\lambda^{u}(\sigma)
    \,d\lambda^{u}(\gamma)\,d\lambda^{u}(\eta)\,d\mu(u) \\ 
    &= \int_{F}\int_{G}f_{1}(\sigma) \bipu<\hatpi(\sigma)f\tensor\spxi_{i},
    g\tensor\spxi_{j}> \Delta(\sigma)^{-\half} \,d\lambda^{u}(\sigma)\,d\mu(u) \\
    &= \int_{F}\int_{G} f_{1}(\gamma)\bipu<U_{\gamma}
    (f\tensor_{s(\gamma)}\spxi_{i}), (g\tensor_{u} \spxi_{j}>
    \Delta(\sigma)^{-\half}
    \,d\lambda^{u}(\sigma)\,d\mu(u) \\
    &= \int_{G}f_{1}(\sigma)\bipb
    r(\sigma)<U_{\sigma}V[f\tensor\spxi_{i}]\bigl(s(\sigma)\bigr) ,
    V[g\tensor\spxi_{j}]\bigl(r(\sigma)\bigr)> \Delta(\sigma)^{-\half}
    \,d\nu(\sigma).
  \end{align*}
  Thus $L'$ is the integrated form as claimed.
\end{proof}

%

\begin{bibdiv}
\begin{biblist}

\bib{anaren:amenable00}{book}{
      author={Anantharaman-Delaroche, Claire},
      author={Renault, Jean},
       title={Amenable groupoids},
      series={Monographies de L'Enseignement Math\'ematique [Monographs of
  L'Enseignement Math\'ematique]},
   publisher={L'Enseignement Math\'ematique},
     address={Geneva},
        date={2000},
      volume={36},
        ISBN={2-940264-01-5},
        note={With a foreword by Georges Skandalis and Appendix B by E.
  Germain},
      review={\MR{MR1799683 (2001m:22005)}},
}

\bib{bla:bsmf96}{article}{
      author={Blanchard, {\'E}tienne},
       title={D\'eformations de {$C\sp *$}-alg\`ebres de {H}opf},
        date={1996},
        ISSN={0037-9484},
     journal={Bull. Soc. Math. France},
      volume={124},
      number={1},
       pages={141\ndash 215},
      review={\MR{97f:46092}},
}

\bib{bou:generalc}{book}{
      author={Bourbaki, Nicolas},
       title={General topology. {C}hapters 1--4},
   publisher={Springer-Verlag},
     address={Berlin},
        date={1989},
        ISBN={3-540-19374-X},
        note={Translated from the French, Reprint of the 1966 edition},
      review={\MR{90a:54001a}},
}

\bib{com:plms84}{article}{
      author={Combes, Fran{\c c}ois},
       title={Crossed products and {M}orita equivalence},
        date={1984},
        ISSN={0024-6115},
     journal={Proc. London Math. Soc. (3)},
      volume={49},
      number={2},
       pages={289\ndash 306},
      review={\MR{86c:46081}},
}

\bib{con:lnm79}{incollection}{
      author={Connes, Alain},
       title={Sur la th\'eorie non commutative de l'int\'egration},
        date={1979},
   booktitle={Alg\`ebres d'op\'erateurs ({S}\'em., les {P}lans-sur-{B}ex,
  1978)},
      series={Lecture Notes in Math.},
      volume={725},
   publisher={Springer},
     address={Berlin},
       pages={19\ndash 143},
      review={\MR{MR548112 (81g:46090)}},
}

\bib{con:pspm80}{incollection}{
      author={Connes, Alain},
       title={A survey of foliations and operator algebras},
        date={1982},
   booktitle={Operator algebras and applications, {P}art {I} ({K}ingston,
  {O}nt., 1980)},
      series={Proc. Sympos. Pure Math.},
      volume={38},
   publisher={Amer. Math. Soc.},
     address={Providence, R.I.},
       pages={521\ndash 628},
      review={\MR{MR679730 (84m:58140)}},
}

\bib{con:course}{book}{
      author={Conway, John~B.},
       title={A course in functional analysis},
      series={Graduate texts in mathematics},
   publisher={Springer-Verlag},
     address={New York},
        date={1985},
      volume={96},
}

\bib{cmw:pams84}{article}{
      author={Curto, Ra{\'u}l~E.},
      author={Muhly, Paul~S.},
      author={Williams, Dana~P.},
       title={Cross products of strongly {M}orita equivalent {$C\sp{\ast}
  $}-algebras},
        date={1984},
        ISSN={0002-9939},
     journal={Proc. Amer. Math. Soc.},
      volume={90},
      number={4},
       pages={528\ndash 530},
      review={\MR{85i:46083}},
}

\bib{dg:banach}{book}{
      author={Dupr{\'e}, Maurice~J.},
      author={Gillette, Richard~M.},
       title={Banach bundles, {B}anach modules and automorphisms of
  ${C}^*$-algebras},
   publisher={Pitman (Advanced Publishing Program)},
     address={Boston, MA},
        date={1983},
      volume={92},
        ISBN={0-273-08626-X},
      review={\MR{85j:46127}},
}

\bib{echwil:jot01}{article}{
      author={Echterhoff, Siegfried},
      author={Williams, Dana~P.},
       title={Locally inner actions on {$C_0(X)$}-algebras},
        date={2001},
        ISSN={0379-4024},
     journal={J. Operator Theory},
      volume={45},
      number={1},
       pages={131\ndash 160},
      review={\MR{1 823 065}},
}

\bib{edw:functional}{book}{
      author={Edwards, Robert~E.},
       title={Functional analysis. {T}heory and applications},
   publisher={Holt},
     address={Rinehart and Winston, New York},
        date={1965},
      review={\MR{MR0221256 (36 \#4308)}},
}

\bib{fd:representations1}{book}{
      author={Fell, James M.~G.},
      author={Doran, Robert~S.},
       title={Representations of {$*$}-algebras, locally compact groups, and
  {B}anach {$*$}-algebraic bundles. {V}ol. 1},
      series={Pure and Applied Mathematics},
   publisher={Academic Press Inc.},
     address={Boston, MA},
        date={1988},
      volume={125},
        ISBN={0-12-252721-6},
        note={Basic representation theory of groups and algebras},
      review={\MR{90c:46001}},
}

\bib{gre:am78}{article}{
      author={Green, Philip},
       title={The local structure of twisted covariance algebras},
        date={1978},
     journal={Acta Math.},
      volume={140},
       pages={191\ndash 250},
}

\bib{hah:tams78}{article}{
      author={Hahn, Peter},
       title={Haar measure for measure groupoids},
        date={1978},
        ISSN={0002-9947},
     journal={Trans. Amer. Math. Soc.},
      volume={242},
       pages={1\ndash 33},
      review={\MR{MR496796 (82a:28012)}},
}

\bib{hr:abstract}{book}{
      author={Hewitt, Edwin},
      author={Ross, Kenneth~A.},
       title={Abstract harmonic analysis. {V}ol. {I}: {S}tructure of
  topological groups. {I}ntegration theory, group representations},
      series={Die Grundlehren der mathematischen Wissenschaften, Band 115},
   publisher={Springer-Verlag},
     address={New York},
        date={1963},
      review={\MR{MR0156915 (28 \#158)}},
}

\bib{hof:74}{unpublished}{
      author={Hofmann, Karl~Heinrich},
       title={Banach bundles},
        date={1974},
        note={Darmstadt Notes},
}

\bib{hof:lnim77}{incollection}{
      author={Hofmann, Karl~Heinrich},
       title={Bundles and sheaves are equivalent in the category of {B}anach
  spaces},
        date={1977},
   booktitle={{$K$}-theory and operator algebras ({P}roc. {C}onf., {U}niv.
  {G}eorgia, {A}thens, {G}a., 1975)},
      series={Lecture Notes in Math},
      volume={575},
   publisher={Springer},
     address={Berlin},
       pages={53\ndash 69},
      review={\MR{58 \#7117}},
}

\bib{hofkei:lnm79}{incollection}{
      author={Hofmann, Karl~Heinrich},
      author={Keimel, Klaus},
       title={Sheaf-theoretical concepts in analysis: bundles and sheaves of
  {B}anach spaces, {B}anach {$C(X)$}-modules},
        date={1979},
   booktitle={Applications of sheaves ({P}roc. {R}es. {S}ympos. {A}ppl. {S}heaf
  {T}heory to {L}ogic, {A}lgebra and {A}nal., {U}niv. {D}urham, {D}urham,
  1977)},
      series={Lecture Notes in Math.},
      volume={753},
   publisher={Springer},
     address={Berlin},
       pages={415\ndash 441},
      review={\MR{MR555553 (81f:46085)}},
}

\bib{hrw:tams00}{article}{
      author={Huef, Astrid~an},
      author={Raeburn, Iain},
      author={Williams, Dana~P.},
       title={An equivariant {B}rauer semigroup and the symmetric imprimitivity
  theorem},
        date={2000},
        ISSN={0002-9947},
     journal={Trans. Amer. Math. Soc.},
      volume={352},
      number={10},
       pages={4759\ndash 4787},
      review={\MR{2001b:46107}},
}

\bib{kas:im88}{article}{
      author={Kasparov, Gennadi~G.},
       title={Equivariant {$KK$}-theory and the {N}ovikov conjecture},
        date={1988},
        ISSN={0020-9910},
     journal={Invent. Math.},
      volume={91},
      number={1},
       pages={147\ndash 201},
      review={\MR{88j:58123}},
}

\bib{khoska:jram02}{article}{
      author={Khoshkam, Mahmood},
      author={Skandalis, Georges},
       title={Regular representation of groupoid {$C\sp *$}-algebras and
  applications to inverse semigroups},
        date={2002},
        ISSN={0075-4102},
     journal={J. Reine Angew. Math.},
      volume={546},
       pages={47\ndash 72},
      review={\MR{MR1900993 (2003f:46084)}},
}

\bib{khoska:jot04}{article}{
      author={Khoshkam, Mahmood},
      author={Skandalis, Georges},
       title={Crossed products of {$C\sp *$}-algebras by groupoids and inverse
  semigroups},
        date={2004},
        ISSN={0379-4024},
     journal={J. Operator Theory},
      volume={51},
      number={2},
       pages={255\ndash 279},
      review={\MR{MR2074181 (2005f:46122)}},
}

\bib{kmrw:ajm98}{article}{
      author={Kumjian, Alexander},
      author={Muhly, Paul~S.},
      author={Renault, Jean~N.},
      author={Williams, Dana~P.},
       title={The {B}rauer group of a locally compact groupoid},
        date={1998},
        ISSN={0002-9327},
     journal={Amer. J. Math.},
      volume={120},
      number={5},
       pages={901\ndash 954},
      review={\MR{2000b:46122}},
}

\bib{gal:94}{thesis}{
      author={Le~Gall, Pierre-Yves},
       title={Th\'eorie de {K}asparov \'equivariante et groupo\"\i des},
        type={Th\`ese de Doctorat},
        date={1994},
}

\bib{gal:kt99}{article}{
      author={Le~Gall, Pierre-Yves},
       title={Th\'eorie de {K}asparov \'equivariante et groupo\"\i des. {I}},
        date={1999},
        ISSN={0920-3036},
     journal={$K$-Theory},
      volume={16},
      number={4},
       pages={361\ndash 390},
      review={\MR{2000f:19006}},
}

\bib{gal:cm01}{incollection}{
      author={Le~Gall, Pierre-Yves},
       title={Groupoid {$C\sp *$}-algebras and operator {$K$}-theory},
        date={2001},
   booktitle={Groupoids in analysis, geometry, and physics ({B}oulder, {CO},
  1999)},
      series={Contemp. Math.},
      volume={282},
   publisher={Amer. Math. Soc.},
     address={Providence, RI},
       pages={137\ndash 145},
      review={\MR{2002h:22005}},
}

\bib{muh:cbms}{techreport}{
      author={Muhly, Paul~S.},
       title={Coordinates in operator algebra},
 institution={CMBS Conference Lecture Notes (Texas Tech 1990)},
        date={1999},
        note={In continuous preparation},
}

\bib{mrw:jot87}{article}{
      author={Muhly, Paul~S.},
      author={Renault, Jean~N.},
      author={Williams, Dana~P.},
       title={Equivalence and isomorphism for groupoid {$C^*$}-algebras},
        date={1987},
        ISSN={0379-4024},
     journal={J. Operator Theory},
      volume={17},
      number={1},
       pages={3\ndash 22},
      review={\MR{88h:46123}},
}

\bib{muhwil:ms90}{article}{
      author={Muhly, Paul~S.},
      author={Williams, Dana~P.},
       title={Continuous trace groupoid {\cs}-algebras},
        date={1990},
     journal={Math. Scand.},
      volume={66},
       pages={231\ndash 241},
}

\bib{muhwil:plms395}{article}{
      author={Muhly, Paul~S.},
      author={Williams, Dana~P.},
       title={Groupoid cohomology and the {D}ixmier-{D}ouady class},
        date={1995},
     journal={Proc. London Math. Soc. (3)},
       pages={109\ndash 134},
}

\bib{nil:iumj96}{article}{
      author={Nilsen, May},
       title={{$C^*$}-bundles and {$C_0(X)$}-algebras},
        date={1996},
        ISSN={0022-2518},
     journal={Indiana Univ. Math. J.},
      volume={45},
      number={2},
       pages={463\ndash 477},
      review={\MR{98e:46075}},
}

\bib{pat:groupoids99}{book}{
      author={Paterson, Alan L.~T.},
       title={Groupoids, inverse semigroups, and their operator algebras},
      series={Progress in Mathematics},
   publisher={Birkh\"auser Boston Inc.},
     address={Boston, MA},
        date={1999},
      volume={170},
        ISBN={0-8176-4051-7},
      review={\MR{MR1724106 (2001a:22003)}},
}

\bib{rae:ma88}{article}{
      author={Raeburn, Iain},
       title={Induced {$C\sp *$}-algebras and a symmetric imprimitivity
  theorem},
        date={1988},
        ISSN={0025-5831},
     journal={Math. Ann.},
      volume={280},
      number={3},
       pages={369\ndash 387},
      review={\MR{90k:46144}},
}

\bib{raewil:tams85}{article}{
      author={Raeburn, Iain},
      author={Williams, Dana~P.},
       title={Pull-backs of {$C\sp \ast$}-algebras and crossed products by
  certain diagonal actions},
        date={1985},
        ISSN={0002-9947},
     journal={Trans. Amer. Math. Soc.},
      volume={287},
      number={2},
       pages={755\ndash 777},
      review={\MR{86m:46054}},
}

\bib{rw:morita}{book}{
      author={Raeburn, Iain},
      author={Williams, Dana~P.},
       title={Morita equivalence and continuous-trace {$C^*$}-algebras},
      series={Mathematical Surveys and Monographs},
   publisher={American Mathematical Society},
     address={Providence, RI},
        date={1998},
      volume={60},
        ISBN={0-8218-0860-5},
      review={\MR{2000c:46108}},
}

\bib{ram:am71}{article}{
      author={Ramsay, Arlan},
       title={Virtual groups and group actions},
        date={1971},
     journal={Advances in Math.},
      volume={6},
       pages={253\ndash 322 (1971)},
      review={\MR{43 \#7590}},
}

\bib{ram:am76}{article}{
      author={Ramsay, Arlan},
       title={Nontransitive quasi-orbits in {M}ackey's analysis of group
  extensions},
        date={1976},
        ISSN={0001-5962},
     journal={Acta Math.},
      volume={137},
      number={1},
       pages={17\ndash 48},
      review={\MR{MR0460531 (57 \#524)}},
}

\bib{ram:jfa82}{article}{
      author={Ramsay, Arlan},
       title={Topologies on measured groupoids},
        date={1982},
     journal={J. Funct. Anal.},
      volume={47},
       pages={314\ndash 343},
}

\bib{ren:groupoid}{book}{
      author={Renault, Jean},
       title={A groupoid approach to {\cs}-algebras},
      series={Lecture Notes in Mathematics},
   publisher={Springer-Verlag},
     address={New York},
        date={1980},
      volume={793},
}

\bib{ren:jot87}{article}{
      author={Renault, Jean},
       title={Repr\'esentations des produits crois\'es d'alg\`ebres de
  groupo\"\i des},
        date={1987},
     journal={J. Operator Theory},
      volume={18},
       pages={67\ndash 97},
}

\bib{rie:pspm82}{incollection}{
      author={Rieffel, Marc~A.},
       title={Applications of strong {M}orita equivalence to transformation
  group {$C^*$}-algebras},
        date={1982},
   booktitle={Operator algebras and applications, {P}art {I} ({K}ingston,
  {O}nt., 1980)},
      series={Proc. Sympos. Pure Math.},
      volume={38},
   publisher={Amer. Math. Soc.},
     address={Providence, R.I.},
       pages={299\ndash 310},
      review={\MR{84k:46046}},
}

\bib{rud:real}{book}{
      author={Rudin, Walter},
       title={Real and complex analysis},
   publisher={McGraw-Hill},
     address={New York},
        date={1987},
}

\bib{shu:private88}{misc}{
      author={Shultz, Frederic~W.},
       title={Personal correspondence},
        date={1988},
}

\bib{tu:doc04}{article}{
      author={Tu, Jean-Louis},
       title={Non-{H}ausdorff groupoids, proper actions and {$K$}-theory},
        date={2004},
        ISSN={1431-0635},
     journal={Doc. Math.},
      volume={9},
       pages={565\ndash 597 (electronic)},
      review={\MR{MR2117427 (2005h:22004)}},
}

\bib{wil:jfa81}{article}{
      author={Williams, Dana~P.},
       title={Transformation group {\cs}-algebras with continuous trace},
        date={1981},
     journal={J. Funct. Anal.},
      volume={41},
       pages={40\ndash 76},
}

\bib{wil:crossed}{book}{
      author={Williams, Dana~P.},
       title={Crossed products of {$C{\sp \ast}$}-algebras},
      series={Mathematical Surveys and Monographs},
   publisher={American Mathematical Society},
     address={Providence, RI},
        date={2007},
      volume={134},
        ISBN={978-0-8218-4242-3; 0-8218-4242-0},
      review={\MR{MR2288954 (2007m:46003)}},
}

\end{biblist}
\end{bibdiv}

\def\noopsort#1{}\def\cprime{$'$} \def\sp{^}

\end{document}